\documentclass[12pt,twoside,leqno,openany]{amsart}

\usepackage{amssymb,amsbsy,amsmath,amsfonts,amssymb,amscd,times,
graphics,color,xypic,footmisc,fancyhdr,multicol,fancybox,
graphicx,mathrsfs,rotating,ifthen,wasysym}
\usepackage[all]{xy}

\usepackage[T1]{fontenc}
\newcommand{\C}{\mathbb{C}}

\newcommand{\LL}{\mathbb{L}}
\newcommand{\N}{\mathbb{N}}
\newcommand{\R}{\mathbb{R}}

\newcommand{\red}{\textcolor{red}}

\newcommand{\HEAD}[2]{%
\pagestyle{fancy}
\fancyhead[RO]{\tiny\sf\thepage}
\fancyhead[CO]{{\tiny\sf #1}}
\fancyhead[LE]{\tiny\sf\thepage}
\fancyhead[CE]{{\tiny\sf #2}}
\fancyfoot{}}

\theoremstyle{definition}

\newcommand{\codim}{\text{\footnotesize\sf codim}}

\newcommand{\constant}{\text{\scriptsize\sf constant}}

\newcommand{\CRdim}{\text{\footnotesize\sf CRdim}}

\renewcommand{\det}{\text{\footnotesize\sf det}}

\renewcommand{\dim}{\text{\footnotesize\sf dim}}

\newcommand{\dist}{\text{\footnotesize\sf dist}}

\newcommand{\explain}[1]{\text{\scriptsize\sf [#1]}}

\newcommand{\function}{\text{\footnotesize\sf function}}

\newcommand{\genrank}{\text{\footnotesize\sf genrank}}

\newcommand{\Id}{\text{\scriptsize\sf Id}}

\renewcommand{\Im}{\text{\footnotesize\sf Im}}

\newcommand{\isqrt}{{\scriptstyle{\sqrt{-1}}}}

\newcommand{\Jac}{\text{\footnotesize\sf Jac}}

\renewcommand{\lim}{\text{\footnotesize\sf lim}}

\newcommand{\mathmotsf}[1]{\text{\footnotesize\sf #1}}

\renewcommand{\mod}{\text{\footnotesize\sf mod}}

\newcommand{\NN}{\text{\scriptsize\sc n}}

\newcommand{\NNN}{\text{\sc n}}

\newcommand{\radius}{\text{\scriptsize\sf radius}}

\newcommand{\rank}{\text{\footnotesize\sf rank}}

\renewcommand{\Re}{\text{\footnotesize\sf Re}}

\newcommand{\same}{\text{\footnotesize\sf same}}

\newcommand{\something}{\text{\scriptsize\sf something}}

\newcommand{\Span}{\text{\footnotesize\sf Span}}

\newcommand{\vf}{\vfill\end{document}}

\newcommand{\zero}[1]{\underline{#1}_{\red{\circ}}}

\let\mathcal\mathscr

\begin{document}

$\:$

\bigskip\bigskip

\begin{center}

{\Large\bf Equivalences of $5$-dimensional CR-manifolds}

\medskip

{\Large\bf II: General classes 
$\text{\sf I}$, 
$\text{\sf II}$, 
$\text{\sf III}_{\text{\sf 1}}$,
$\text{\sf III}_{\text{\sf 2}}$,
$\text{\sf IV}_{\text{\sf 1}}$,
$\text{\sf IV}_{\text{\sf 2}}$
}

\end{center}

\medskip

\begin{center}
Jo\"el {\sc Merker}, 
Samuel {\sc Pocchiola},
and Masoud {\sc Sabzevari}
\end{center}

\bigskip

\begin{center}
\begin{minipage}[t]{10.25cm}
\baselineskip =0.32cm 
{\scriptsize
{\bf Abstract.}
For later use in subsequent upcoming arxiv.org prepublications, basic
foundational material on local, smooth or real analytic, CR-generic
submanifolds of complex Euclidean spaces is developed from scratch,
with strong emphasis on the interplay between extrinsic and intrinsic
aspects, a constructive option that commands to perform
computational syntheses in 
coordinates. Mainly, one finds a self-contained
proof of the existence of precisely six general classes
$\text{\sf I}$, 
$\text{\sf II}$, 
$\text{\sf III}_{\text{\sf 1}}$,
$\text{\sf III}_{\text{\sf 2}}$,
$\text{\sf IV}_{\text{\sf 1}}$,
$\text{\sf IV}_{\text{\sf 2}}$
of nondegenerate general CR manifolds up to dimension 5,
class $\text{\sf III}_{\text{\sf 2}}$ being unobserved untill now.
}
\end{minipage}
\end{center}

\medskip

\begin{center}
\begin{minipage}[t]{11.75cm}
\baselineskip =0.35cm {\scriptsize

\centerline{\bf Table of contents (excerpt)}

\medskip

{\bf \ref{yet-new-class}.~New general class of {\bf 5}-dimensional 
CR manifolds $M^5 \subset \C^4$
\dotfill~\pageref{yet-new-class}.}

{\bf \ref{freeman-n-2}.~Levi kernel and Freeman form nondegeneracies 
in CR dimension $n = {\bf 2}$
\dotfill~\pageref{freeman-n-2}.}

{\bf \ref{summary-general-classes}.~Classes 
$\text{\sf I}$, 
$\text{\sf II}$, 
$\text{\sf III}_{\text{\sf 1}}$,
$\text{\sf III}_{\text{\sf 2}}$,
$\text{\sf IV}_{\text{\sf 1}}$,
$\text{\sf IV}_{\text{\sf 2}}$,
of 
$M^3 \subset \C^2$, $M^4 \subset \C^3$, $M^5 \subset \C^4$, $M^5 \subset \C^3$
\dotfill~\pageref{summary-general-classes}.}

}\end{minipage}
\end{center}

%%%%%%%%%%%%%%%%%%%%%%%%%%%%%%%%%%%%%%%%%%%%%%%%%%%%%%%%%%%%%%%%%%%%%

\bigskip

\section{\sf Real analytic ($\mathcal{ C}^\omega$)
submanifolds of $\C^\NN$:
\\
Zariski-generic features}
\label{real-analytic-submanifolds}
\HEAD{\ref{real-analytic-submanifolds}.~Real 
analytic ($\mathcal{ C}^\omega$)
submanifolds of $\C^\NN$:
Zariski-generic features}{
Jo\"el {\sc Merker} (Paris-Sud), 
Samuel {\sc Pocchiola} (Paris-Sud), 
Masoud {\sc Sabzevari} (Shahrekord)}

\medskip

\medskip\noindent{\bf Smoothness class assumptions.} In what
follows:
\[
\mathmotsf{smoothness classes}
:=
\left\{
\aligned
\mathcal{C}^{1,2,3,4,5,6,7,8,9,\dots},
\\
\mathcal{C}^\infty,
\\
\mathcal{C}^\omega.
\endaligned
\right.
\]

\medskip\noindent{\bf Complex Euclidean space.}
On $\C^\NN = \R^{2\NN}$, take $\text{\sc n}$ complex coordinates:
\[
\big({\sf z}_1,\dots,{\sf z}_\NN\big)
=
\big(
{\sf x}_1+\isqrt\,{\sf y}_1,\dots,
{\sf x}_\NN+\isqrt\,{\sf y}_\NN\big).
\]
On the real tangent bundle:
\[
T\R^{2\NN} 
\cong
T^{\sf real}\C^\NN,
\]
for which a natural frame is constituted by the $2\NNN$ vector fields:
\[
\frac{\partial}{\partial {\sf x}_1},\,\,
\dots\dots,\,\,
\frac{\partial}{\partial {\sf x}_\NN},\,\,\,\,
\frac{\partial}{\partial {\sf y}_1},\,\,
\dots\dots,\,\,
\frac{\partial}{\partial {\sf y}_\NN},\,\,
\]
the complex structure $J$ acts by definition as:
\[
\aligned
J
\bigg(
\frac{\partial}{\partial {\sf x}_k}
\bigg)
&
:=
\frac{\partial}{\partial{\sf y}_k}
\ \ \ \ \ \ \ \ \ \ \ \ \ \ \ \
{\scriptstyle{(k\,=\,1\,\cdots\,\NN)}},
\\
J
\bigg(
\frac{\partial}{\partial{\sf y}_k}
\bigg)
&
:=
-\,
\frac{\partial}{\partial{\sf x}_k}
\ \ \ \ \ \ \ \ \ \ \ \ \
{\scriptstyle{(k\,=\,1\,\cdots\,\NN)}},
\endaligned
\]
and it is an invertible automorphism of $T\R^{ 2\NN}$ satisfying:
\[
J^2
=
-\,\Id.
\]

\medskip\noindent{\bf CR submanifolds of $\C^\NN$.}
Consider a $\mathcal{ C}^\kappa$ ($\kappa \geqslant 1$), or
$\mathcal{ C}^\infty$, or $\mathcal{ C}^\omega$ {\em real}
submanifold:
\[
M
\subset
\C^\NN.
\]
At each point $q \in M$, one may view {\em extrinsically}:
\[
T_qM
\subset
T_q\C^\NN
=
T_q\R^{2\NN},
\]
so that it is meaningful to consider the vector subspace:
\[
J\big(T_qM\big)
\subset
T_q\R^{2\NN},
\]
which is of the same dimension as $T_q M$.

\medskip\noindent{\bf First elementary fact.}
{\em Any $\mathcal{ C}^\omega$ connected real submanifold:}
\[
M 
\subset
\C^\NN
\]
{\em is {\sl Cauchy-Riemann} (CR) on a certain 
Zariski open subset:}
\[
M\backslash\Sigma,
\] 
{\em in the sense that:}
\[
M\backslash\Sigma
\ni
p
\,\longmapsto\,
\dim_\R 
\big(
T_p M\cap J(T_p M)\big)
\in
\N
\] 
{\em has constant value there.\qed}

\medskip\noindent{\bf Definition.}
At various points $p \in M$, introduce the {\sl
complex-tangent} subspaces:
\[
T_p^cM
:=
T_pM
\cap
J(T_pM)
\ \ \ \ \ \ \ \ \ \ \ \ \
{\scriptstyle{(p\,\in\,M)}}
\]
of the tangent spaces $T_pM$.

\medskip
Applying then the:
\[
\text{\sl Lie-Cartan Principle of Relocalization},
\]
one disregards the non-CR locus $\Sigma$ and one
assumes that $M$ is CR at every point:
\[
\dim_\R
\big(T_pM\cap J(T_pM)\big)
\,=\,
\mathmotsf{constant},
\]
so that:
\[
\bigcup_{p\in M}\,
T_p^cM
\,\subset\,
\bigcup_{p\in M}\,
T_pM
\]
makes up a true {\em subbundle}:
\[
T^cM
\subset
TM.
\] 

Furthermore, setting:
\[
\big(a+\isqrt\,b\big)
\cdot
L_p
\,:=\,
a\,L_p
+
\isqrt\,
J\big(L_p\big),
\]
one naturally equips all the: 
\[
T_p^cM
\,\ni\,
L_p
\]
with {\em complex vector space structures}, whence:
\[
\dim_\R
\big(T_p^cM\big)
\,\in\,2\,\N.
\]

Similarly, the $J$-invariance\,\,---\,\,use $J^2 = \Id$\,\,---\,\,of:
\[
T_pM\cap J(T_pM)
\]
yields:
\[
\dim_\R
\big(T_pM+J(T_pM)\big)
\,\in\,2\,\N.
\]

\smallskip

When $M$ is CR (at every point), 
the dimension formula:
\[
\dim_\R\,\big(E+F\big)
=
\dim_\R\,E
+
\dim_\R\,F
-
\dim_\R\big(E\cap F\big)
\]
for vector subspaces $E, F$ of a certain ambient 
vector space then gives:
\[
\aligned
\underbrace{\rank_\R
\big(TM+J(TM)\big)}_{\in\,2\,\N}
&
=
\rank_\R\big(TM\big)
+
\rank_\R\big(J(TM)\big)
-
\underbrace{\rank_\R
\big(TM\cap J(TM)\big)}_{\in\,2\,\N}
\\
&
=
{\sf constant}
\\
&
=:
2\,\NNN_M
\\
&
\leqslant\,
2\,\NNN
\,=\,
\rank_\R\big(T\C^\NN\big).
\endaligned
\]

\medskip\noindent{\bf Second elementary fact.}
{\em Any $\mathcal{ C}^\omega$ connected
CR submanifold:} 
\[
M 
\subset 
\C^\NN
\]
{\em is contained in a unique thin {\em complex-analytic} strip-submanifold:}
\[
M^{i_c} 
\supset 
M
\]
{\em stretched along $M$ of complex dimension:}
\[
\aligned
\dim_\C\,M^{i_c}
&
\,=\,
\NNN_M
\\
&
\,=\,
\rank_\C\big(TM+J( TM)\big),
\endaligned
\]
{\em inside which $M$ is {\sl CR-generic}:}
\[
T_pM+J(T_pM)
=
T_pM^{i_c}
\ \ \ \ \ \ \ \ \ \ \ \ \
{\scriptstyle{(\forall\,\,p\,\in\,M)}}.
\qed
\] 

\medskip\noindent{\bf Consequence for the
biholomorphic equivalence problem.}
{\em After replacing $\C^\NN$ by:}
\[
M^{i_c}
\cong 
\C^{\NN_M}
\ \ \ \ \ \ \ \ \ \ \ \ \
{\scriptstyle{(\NN_M\,=\,{\sf rank}_\C(TM+J(TM)))}}, 
\]
{\em there is no restriction to study only {\em CR-generic} 
$\mathcal{ C}^\omega$ submanifolds $M \subset \C^\NN$.\qed}

\medskip

From now on, therefore:
\[
M 
\subset 
\C^\NN
\ \ \ \ \ \ \ \ \ \ \ \ \ \ \ \ \ \ \ \ \
\text{\rm with}
\ \ \ \ \
TM+J(TM) 
=
T\C^\NN
\big\vert_M
\]
will {\em always} be CR-generic. Introduce:
\[
c
:=
\codim_\R\,M.
\]

\smallskip

The dimension formula, again, then yields:
\[
\aligned
2\,\NNN
=
\rank_\R\big(T\C^\NN\big)
&
=
\rank_\R
\big(TM+J(TM)\big)
\\
&
=
\underbrace{\rank_\R\big(TM\big)}_{
2\,\NN\,-\,c}
+
\underbrace{\rank_\R\big(J(TM )\big)}_{
2\,\NN\,-\,c}
-
\rank_\R\big(TM\cap J(TM )\big),
\endaligned
\]
so that:
\[
\rank_\R
\big(
TM\cap J( TM)
\big)
=
2\,\NNN-2\,c.
\]

\medskip\noindent{\bf Definition-Property.}
The {\sl CR dimension} of a CR submanifold:
\[
M
\subset
\C^\NN
\]
is the rank as a $\C$-vector bundle of:
\[
\aligned
\CRdim\,M
&
\overset{\sf def}{\,:=\,}
\rank_\C\big(T^cM\big)
\\
&
\,=\,
{\textstyle{\frac{1}{2}}}\,
\rank_\R\big(TM\cap J( TM)\big)
\endaligned
\]
and when $M$ is CR-generic, one has:
\[
\aligned
\CRdim\,M
&
=
\NNN
-
\codim_\R\,M
\\
&
=
\NNN
-
c.
\qed
\endaligned
\]

\smallskip
Regularly, the CR dimension will be denoted with the letter:
\[
n
:=
\CRdim\,M.
\]
An application of the Implicit Function Theorem yields the known:

\medskip\noindent{\bf Proposition.}
{\em When $M \subset \C^\NN$ is CR-generic with:}
\[
\aligned
c
&
=
\codim_\R\,M,
\\
n
&
=
\CRdim\,M
=
\NNN-c,
\endaligned
\]
{\em then at every point:}
\[
p\in M
\]
{\em and for every choice of centered affine coordinates:}
\[
\big(z_1,\dots,z_n,w_1,\dots,w_c\big)
=
\Big(
x_1+\isqrt\,y_1,\dots,x_n+\isqrt\,y_n,\,
u_1+\isqrt\,v_1,\dots,u_c+\isqrt\,v_c
\Big)
\]
{\em in which the tangent space is straightened:}
\[
T_pM
=
\big\{
0
=
v_1
=\cdots=
v_c
\big\},
\]
{\em there exist $c$ graphing functions:} 
\[
\varphi_1,\,\dots,\,\varphi_c
\]
{\em which locally represent $M$ as:}
\[
\left[
\aligned
v_1
&
=
\varphi_1
\big(x_1,\dots,x_n,y_1,\dots,y_n,u_1,\dots,u_c\big),
\\
\cdots
&
\cdots\cdots\cdots\cdots\cdots\cdots\cdots\cdots\cdots\cdots\cdots\cdots
\cdot\cdot
\\ 
v_c
&
=
\varphi_c
\big(x_1,\dots,x_n,y_1,\dots,y_n,u_1,\dots,u_c\big),
\endaligned\right.
\]
{\em which are defined in some open neighborhood of the
origin in:}
\[
\R^n\times\R^n\times\R^c,
\]
{\em which satisfy:}
\[
\aligned
0
&
=
\varphi_1(0)
=
d\varphi_1(0),
\\
\cdots
&
\cdots\cdots\cdots\cdots\cdots\cdots
\\
0
&
=
\varphi_c(0)
=
d\varphi_c(0),
\endaligned
\]
{\em and which are of regularity:}
\[
\mathcal{C}^{1,2,3,4,5,6,7,8,9,\dots}
\ \ \ \ \
\text{\rm or}
\ \ \ \ \
\mathcal{C}^\infty
\ \ \ \ \
\text{\rm or}
\ \ \ \ \
\mathcal{C}^\omega,
\]
{\em according to which $M$ is assumed to belong to the
corresponding smoothness category.\qed}

%%%%%%%%%%%%%%%%%%%%%%%%%%%%%%%%%%%%%%%%%%%%%%%%%%%%%%%%%%%%%%%%%%%%%

\bigskip

\section{\sf CR-generic submanifolds up to dimension $5$}
\label{CR-up-5}
\HEAD{\ref{CR-up-5}.~CR-generic submanifolds up to dimension $5$}{
Jo\"el {\sc Merker} (Paris-Sud), 
Samuel {\sc Pocchiola} (Paris-Sud), 
Masoud {\sc Sabzevari} (Shahrekord)}

\medskip

Let therefore $M \subset \C^\NN$ be a connected
CR-generic $\mathcal{ C}^\omega$ submanifold with:
\[
\aligned
c
&
=
\codim_\R\,M,
\\
n
&
=
\CRdim\,M,
\\
2n+c
&\,
=
\dim_\R\,M.
\endaligned
\]

\smallskip
Recall the goal is to reach:
\[
\dim_\R\,M
=
2n+c
\leqslant
{\bf 5}.
\]
The two cases:
\[
\aligned
c
&
=
0,
\\
n
&
=
0
\endaligned
\]
are not interesting in CR geometry, for: 
\[
\aligned
M
\cong
\C^n,
\\
M
\cong
\R^c,
\endaligned
\]
respectively. Hence one assumes:
\[
\aligned
c
&
\geqslant 
1
\\
n
&
\geqslant
1.
\endaligned
\]

\medskip\noindent{\bf Possible CR dimensions and real codimensions:}
\[
\boxed{
\aligned
&
2n+c={\bf 3}
\ \ \ \ \ \
\Longrightarrow
\ \
\Big\{
n={\bf 1},\ \ \ \ \
c={\bf 1},
\\
&
2n+c={\bf 4}
\ \ \ \ \ \
\Longrightarrow
\ \
\Big\{
n={\bf 1},\ \ \ \ \
c={\bf 2},
\\
&
2n+c={\bf 5}
\ \ \ \ \ \
\Longrightarrow
\ \
\left\{
\aligned
&
n={\bf 1},\ \ \ \ \
c={\bf 3},
\\
&
n={\bf 2},\ \ \ \ \
c={\bf 1}.
\endaligned\right.
\endaligned
}
\]

In order to distinguish these cases, dimensions must be emphasized:
\[
M^{2n+c}
\,\subset\,
\C^{n+c},
\]
which gives four cases:
\[
\aligned
&
M^3
\,\subset\,\C^2,
\\
&
M^4
\,\subset\,\C^3,
\\
&
M^5
\,\subset\,
\aligned
&
\C^4,
\\
&
\C^3.
\endaligned
\endaligned
\]

In local coordinates: 
\[
\big(
z_1,\dots,z_n,w_1,\dots,w_n
\big)
=
\Big(
x_1+\isqrt\,y_1,\dots,x_n+\isqrt\,y_n,\,
u_1+\isqrt\,v_1,\dots,u_c+\isqrt\,v_c
\Big),
\]
one represents:
\[
\aligned
&
M^3\subset\C^2
\colon
\ \ \ \ \ \,
\Big[
\,\,
v
=
\varphi(x,y,u),
\\
&
M^4\subset\C^3
\colon
\ \ \ \ \
\left[
\aligned
\,\,
v_1
&
=
\varphi_1(x,y,u_1,u_2),
\\
\,\,
v_2
&
=
\varphi_2(x,y,u_1,u_2),
\endaligned\right.
\\
& 
M^5\subset\C^4
\colon
\ \ \ \ \
\left[
\aligned
\,\,
v_1
&
=
\varphi_1(x,y,u_1,u_2,u_3),
\\
\,\,
v_2
&
=
\varphi_2(x,y,u_1,u_2,u_3), 
\\
\,\,
v_3
&
=
\varphi_3(x,y,u_1,u_2,u_3),
\endaligned\right.
\\
&
M^5\subset\C^3
\colon
\ \ \ \ \ \
\Big[
\,\,
v
=
\varphi(x_1,y_1,x_2,y_2,u),
\endaligned
\]
erasing lower indices when either $n = 1$ or $c = 1$.

%%%%%%%%%%%%%%%%%%%%%%%%%%%%%%%%%%%%%%%%%%%%%%%%%%%%%%%%%%%%%%%%%%%%%

\bigskip

\section{\sf Action of local biholomorphisms
\\
on $T^{1,0}\C^\NN$, on $T^{0,1}\C^\NN$, 
on $T^{1,0}M^{2n+c}$, on $T^{0,1}M^{2n+c}$}
\label{action-local-biholomorphisms}
\HEAD{\ref{action-local-biholomorphisms}.~Action of local biholomorphisms
on $T^{1,0}\C^\NN$, on $T^{0,1}\C^\NN$, 
on $T^{1,0}M^{2n+c}$, on $T^{0,1}M^{2n+c}$}{
Jo\"el {\sc Merker} (Paris-Sud), 
Samuel {\sc Pocchiola} (Paris-Sud), 
Masoud {\sc Sabzevari} (Shahrekord)}

\medskip

\medskip\noindent{\bf Local biholomorphisms.}
On $\C^\NN = \R^{2\NN}$, take $\text{\sc n}$ complex coordinates:
\[
\big({\sf z}_1,\dots,{\sf z}_\NN\big)
=
\big(
{\sf x}_1+\isqrt\,{\sf y}_1,\dots,
{\sf x}_\NN+\isqrt\,{\sf y}_\NN\big),
\]
that will sometimes be abbreviated as:
\[
{\sf z}_\bullet
=
{\sf x}_\bullet
+
\isqrt\,{\sf y}_\bullet.
\]
On $\C' = {\R'}^2$, take $1$ complex coordinate:
\[
{\sf z}'
=
{\sf x}'
+
\isqrt\,{\sf y}'.
\]

Consider an open subset:
\[
{\sf U}
\subset
\C^\NN,
\]
and a $\mathcal{ C}^\kappa$ ($\kappa \geqslant 1$),
or $\mathcal{ C}^\infty$, or $\mathcal{ C}^\omega$ map:
\[
\aligned
h\colon
\ \ \ \ \ \ \ \ \ \ \ \ \ \ \ \ \ \ \ \ \ \ \ \ \ \ \ \ \ \ \ \ \ \ \ \ 
{\sf U}
&
\,\longrightarrow\,
\C'
\\
\big({\sf x}_1,{\sf y}_1,\dots,{\sf x}_\NN,{\sf y}_\NN\big)
&
\,\longmapsto\,
h\big({\sf x}_1,{\sf y}_1,\dots,{\sf x}_\NN,{\sf y}_\NN\big)
\\
&
\ \ \ \ \
=
f\big({\sf x}_1,{\sf y}_1,\dots,{\sf x}_\NN,{\sf y}_\NN\big)
+
\isqrt\,
g\big({\sf x}_1,{\sf y}_1,\dots,{\sf x}_\NN,{\sf y}_\NN\big)
\endaligned
\]
decomposed in real and imaginary parts:
\[
h
=
f+\isqrt\,g.
\]

\medskip
Introduce the $\NNN$ antiholomorphic vector field derivations:
\[
\frac{\partial}{\partial\overline{\sf z}_1}
\overset{\sf def}{\,:=\,}
\frac{1}{2}\,\frac{\partial}{\partial{\sf x}_1}
+
\frac{\isqrt}{2}\,\frac{\partial}{\partial{\sf y}_1},\,\,
\dots\dots\dots,\,\,
\frac{\partial}{\partial\overline{\sf z}_\NN}
\overset{\sf def}{\,:=\,}
\frac{1}{2}\,\frac{\partial}{\partial{\sf x}_\NN}
+
\frac{\isqrt}{2}\,\frac{\partial}{\partial{\sf y}_\NN}.
\]

\medskip\noindent{\bf Definition.}
A 
$\mathcal{ C}^\kappa$ ($\kappa \geqslant 1$), or $\mathcal{ C}^\infty$, 
or $\mathcal{ C}^\omega$
map $h \colon {\sf U} \longrightarrow \C'$
is {\sl holomorphic} when:
\[
0
\equiv
\frac{\partial h}{\partial\overline{\sf z}_1}
\equiv\cdots\equiv
\frac{\partial h}{\partial\overline{\sf z}_\NN}.
\]

\medskip

Equivalently:
\[
\aligned
0
&
\equiv
\bigg(
\frac{\partial}{\partial{\sf x}_l}
+
\isqrt\,\frac{\partial}{\partial{\sf y}_l}
\bigg)
\Big(
f
+
\isqrt\,g
\Big)
\\
&
\equiv
f_{{\sf x}_l}
-
g_{{\sf y}_l}
+
\isqrt\,
\big(
f_{{\sf y}_l}
+
g_{{\sf x}_l}
\big).
\endaligned
\]

\medskip\noindent{\bf Consequence.}
{\em The map $h = f + \isqrt\, g$ is holomorphic
if and only if:}

\[
\aligned
0
&
\equiv
f_{{\sf x}_l}
-
g_{{\sf y}_l}
\ \ \ \ \ \ \ \ \ \ \ \ \
{\scriptstyle{(l\,=\,1\,\cdots\,\NN)}},
\\
0
&
\equiv
f_{{\sf y}_l}
+
g_{{\sf x}_l}
\ \ \ \ \ \ \ \ \ \ \ \ \
{\scriptstyle{(l\,=\,1\,\cdots\,\NN)}},
\endaligned
\]
{\em these being called} {\sl Cauchy-Riemann equations}.

\medskip\noindent{\bf Known fundamental property.}
{\em Every holomorphic function $h$ on some open subset:}
\[
{\sf U}
\subset
\C^\NN
\]
{\em is locally expandable in converging power series:}
\[
h\big({\sf z}_1,\dots,{\sf z}_\NN\big)
=
\sum_{\alpha_1\in\N}\cdots\sum_{\alpha_\NN\in\N}\,
\underbrace{h_{\alpha_1,\dots,\alpha_\NN}}_{\in\,\C}\,
\big({\sf z}_1-{\sf z}_{01}\big)^{\alpha_1}\,
\cdots\cdots
\big({\sf z}_\NN-{\sf z}_{0\NN}\big)^{\alpha_\NN}
\]
{\em in some sufficiently small neighborhood:}
\[
\big\{
\vert{\sf z}_1-{\sf z}_{10}\vert
<
\rho_0,
\,\dots\dots,\,
\vert{\sf z}_\NN-{\sf z}_{\NN 0}\vert
<
\rho_0
\big\}
\]
{\em of any point:}
\[
{\sf z}_{0\bullet}
=
\big({\sf z}_{01},\dots,{\sf z}_{0n}\big)
\,\in\,
{\sf U},
\]
{\em for some $\rho_0$ with:}
\[
0
<
\rho_0
\leqslant
\dist\,
\big({\sf z}_{0\bullet},\,{\sf boundary}({\sf U})\big),
\]
{\em that is to say the coefficient enjoy a Cauchy-type
estimate:}
\[
\big\vert
h_{\alpha_1,\dots,\alpha_\NN}
\big\vert
\,\leqslant\,
\constant\,
\bigg(
\frac{1}{\radius}
\bigg)^{\alpha_1+\cdots+\alpha_\NN},
\]
{\em for some two positive constants:}
\[
\constant\,>\,0,
\ \ \ \ \ \ \ \ \ \ \ \ \ \ \ \ \ \ \ \ \ \ \ \ \ \
\radius\,>\,0.
\qed
\]

\medskip
In all what follows, no attention will be paid to making
any occurence of the constant $\radius$ close to any true
radius of convergence, just its positivity will matter.
Also, one will sometimes abbreviate:
\[
h({\sf z}_\bullet)
=
\sum_{\alpha_\bullet\in\N^\NN}\,
h_{\alpha_\bullet}\,
\big({\sf z}_\bullet-{\sf z}_{0\bullet}\big)^{\alpha_\bullet}.
\]

\medskip

Importantly, when one conjugates a holomorphic function,
the conjugation instantly distributes onto its converging 
power series:
\[
\overline{h({\sf z}_\bullet)}
=
\sum_{\alpha_\bullet\in\N^\NN}\,
\overline{h_{\alpha_\bullet}}\,
\big(\overline{\sf z}_\bullet
-
\overline{\sf z}_{0\bullet}\big)^{\alpha_\bullet},
\]
so that one can introduce:
\[
\overline{h}
\big({\sf z}_\bullet\big)
:=
\sum_{\alpha_\bullet\in\N^\NN}\,
\overline{h_{\alpha_\bullet}}\,
\big({\sf z}_\bullet-{\sf z}_{0\bullet}\big)^{\alpha_\bullet},
\]
by conjugating {\em only} the coefficients.

\medskip\noindent{\bf Transfer of vector fields.}
Given a $\mathcal{ C}^\kappa$ ($\kappa \geqslant 1$), 
or $\mathcal{ C}^\infty$, 
or $\mathcal{ C}^\omega$ map $h = (f, g)$:
\[
\C^\NN
=
\R^{2\NN}
\,\longrightarrow\,
{\R'}^2
=
{\C'},
\]
written out as:
\[
\big({\sf x}_1,{\sf y}_1,
\dots,{\sf x}_\NN,{\sf y}_\NN\big)
\,\longmapsto\,
\Big(
f\big({\sf x}_1,{\sf y}_1,
\dots,{\sf x}_\NN,{\sf y}_\NN\big),\,\,
g\big({\sf x}_1,{\sf y}_1,
\dots,{\sf x}_\NN,{\sf y}_\NN\big)
\Big),
\]
its $2 \times 2\text{\sc n}$ {\sl Jacobian matrix}:
\[
\Jac_\R(h)
=
\Jac_\R(f,g)
=
\left(\!
\begin{array}{ccccc}
f_{{\sf x}_1} & f_{{\sf y}_1} & \cdots &
f_{{\sf x}_\NN} & f_{{\sf y}_\NN}
\\ 
g_{{\sf x}_1} & g_{{\sf y}_1} & \cdots &
g_{{\sf x}_\NN} & g_{{\sf y}_\NN}
\end{array}
\!\right)
\]
expresses the rank of $h$ at various points. 

Moreover, $\Jac_\R ( f, g)$ enables one to 
{\em transfer} tangent vectors:
\[
h_*
=
(f,g)_*
\colon\ \ \
T{\sf U}
\,\longrightarrow\,
T\C',
\]
with the understanding that tangent vectors identify
with {\em derivations}.

The best way to see this is to look first at the transfer
of functions by composition:
\[
\xymatrix{
{\sf U} \ar@{.>}[dr]_F \ar[r]^h & {\sf U}' \ar[d]^{F'}
\\
& \R''.
}
\]
Assume:
\[
h({\sf U})
\subset
{\sf U}'
\subset
\C'.
\]
Then to every real-valued function:
\[
F'
\colon\ \ \
{\sf U}'
\,\longrightarrow\,
\R'',
\]
one associates:
\[
F
:=
F'\circ h,
\]
namely:
\[
\aligned
F
\big({\sf x}_1,{\sf y}_1,
\dots,{\sf x}_\NN,{\sf y}_\NN\big)
&
=
F'\big(
{\sf x}',{\sf y}'\big)
\circ
h\big({\sf x}_1,{\sf y}_1,
\dots,{\sf x}_\NN,{\sf y}_\NN\big)
\\
&
=
F'
\Big(
f
\big({\sf x}_1,{\sf y}_1,
\dots,{\sf x}_\NN,{\sf y}_\NN\big),\,\,
g
\big({\sf x}_1,{\sf y}_1,
\dots,{\sf x}_\NN,{\sf y}_\NN\big)
\Big).
\endaligned
\]
Applying the chain rule, one gets:
\[
\aligned
\frac{\partial F}{\partial{\sf x}_l}
&
=
f_{{\sf x}_l}\,\frac{\partial F'}{\partial{\sf x}'}
+
g_{{\sf x}_l}\,\frac{\partial F'}{\partial{\sf y}'}
\ \ \ \ \ \ \ \ \ \ \ \ \
{\scriptstyle{(l\,=\,1\,\cdots\,\NN)}},
\\
\frac{\partial F}{\partial{\sf y}_l}
&
=
f_{{\sf y}_l}\,\frac{\partial F'}{\partial{\sf x}'}
+
g_{{\sf y}_l}\,\frac{\partial F'}{\partial{\sf y}'}
\ \ \ \ \ \ \ \ \ \ \ \ \
{\scriptstyle{(l\,=\,1\,\cdots\,\NN)}},
\endaligned
\]
identically for:
\[
\big({\sf x}_1,{\sf y}_1,
\dots,{\sf x}_\NN,{\sf y}_\NN\big)
\,\in\,
{\sf U},
\]
without writing arguments.

This means that $h = (f, g)$ pushes forward vector fields as:
\[
\aligned
(f,g)_*
\bigg(
\frac{\partial}{\partial{\sf x}_l}
\bigg)
=
f_{{\sf x}_l}\,\frac{\partial}{\partial{\sf x}'}
+
g_{{\sf x}_l}\,\frac{\partial}{\partial{\sf y}'}
\ \ \ \ \ \ \ \ \ \ \ \ \
{\scriptstyle{(l\,=\,1\,\cdots\,\NN)}},
\\
(f,g)_*
\bigg(
\frac{\partial}{\partial{\sf y}_l}
\bigg)
=
f_{{\sf y}_l}\,\frac{\partial}{\partial{\sf x}'}
+
g_{{\sf y}_l}\,\frac{\partial}{\partial{\sf y}'}
\ \ \ \ \ \ \ \ \ \ \ \ \
{\scriptstyle{(l\,=\,1\,\cdots\,\NN)}}.
\endaligned
\]

In the right-hand sides, the coefficients:
\[
\aligned
&
f_{{\sf x}_l},
\ \ \ \ \ \ \ \ \ \
g_{{\sf x}_l},
\\
&
f_{{\sf y}_l},
\ \ \ \ \ \ \ \ \ \
g_{{\sf y}_l},
\endaligned
\]
live in the source space ${\sf U}$, while the fields:
\[
\frac{\partial}{\partial{\sf x}'},
\ \ \ \ \ \ \ \ \ \
\frac{\partial}{\partial{\sf y}'},
\]
live in the target space $\C'$.

To remedy this imperfection, it is better to deal with 
{\em equidimensional} source and target spaces.

Let therefore:
\[
\big({\sf z}_1',\dots,{\sf z}_\NN'\big)
=
\big({\sf x}_1',{\sf y}_1',\dots,{\sf x}_\NN',{\sf y}_\NN'\big)
\]
be complex coordinates on a target space:
\[
{\C'}^\NN
=
{\R'}^{2\NN}
\]
having the same dimension.

When advisable, abbreviate the coordinates as:
\[
\aligned
\big({\sf x}_\bullet,{\sf y}_\bullet\big)
\ \ 
&
\text{\rm on}\ \
\C^\NN
=
\R^{2\NN},
\\
\big({\sf x}_\bullet',{\sf y}_\bullet'\big)
\ \ 
&
\text{\rm on}\ \
{\C'}^\NN
=
{\R'}^{2\NN}.
\endaligned
\]

Consider an open subset:
\[
{\sf U}
\subset
\C^\NN,
\]
and a $\mathcal{ C}^\kappa$ ($\kappa \geqslant 1$), 
or $\mathcal{ C}^\infty$, 
or $\mathcal{ C}^\omega$ map:
\[
\aligned
h\colon
\ \ \ \ \ \ \ \ \ \ \ \ \ \ \ \ \ \ \ \ \ \ \ \ \ \ \ \ \ \ \ \ \ \ \ \ 
{\sf U}
&
\,\longrightarrow\,
{\C'}^\NN
=
{\R'}^{2\NN}
\\
\big({\sf x}_1,{\sf y}_1,\dots,{\sf x}_\NN,{\sf y}_\NN\big)
&
\,\longmapsto\,
\Big(
h_1\big({\sf x}_\bullet,{\sf y}_\bullet\big),
\dots,
h_\NN\big({\sf x}_\bullet,{\sf y}_\bullet\big)
\Big)
\\
&
\ \ \ \ \
=:
\Big(
f_1\big({\sf x}_\bullet,{\sf y}_\bullet\big),
g_1\big({\sf x}_\bullet,{\sf y}_\bullet\big),
\dots,
f_\NN\big({\sf x}_\bullet,{\sf y}_\bullet\big),
g_\NN\big({\sf x}_\bullet,{\sf y}_\bullet\big)
\Big)
\endaligned
\]
decomposed in real and imaginary parts:
\[
h_1
=
f_1+\isqrt\,g_1,
\,\,\dots\dots,\,\,
h_\NN
=
f_\NN+\isqrt\,g_\NN.
\]

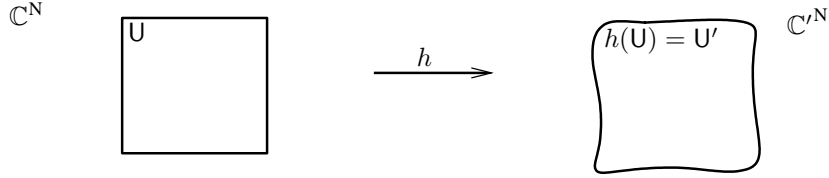
\begin{center}
\begin{picture}(0,0)%
\includegraphics{h-U-polydiscs.pstex}%
\end{picture}%
\setlength{\unitlength}{4144sp}%
\begingroup\makeatletter\ifx\SetFigFont\undefined%
\gdef\SetFigFont#1#2#3#4#5{%
  \reset@font\fontsize{#1}{#2pt}%
  \fontfamily{#3}\fontseries{#4}\fontshape{#5}%
  \selectfont}%
\fi\endgroup%
\begin{picture}(4680,1028)(887,-2319)
\put(3339,-1656){\makebox(0,0)[lb]{\smash{{\SetFigFont{10}{12.0}{\familydefault}{\mddefault}{\updefault}{\color[rgb]{0,0,0}$h$}%
}}}}
\put(1612,-1507){\makebox(0,0)[lb]{\smash{{\SetFigFont{10}{12.0}{\familydefault}{\mddefault}{\updefault}{\color[rgb]{0,0,0}${\sf U}$}%
}}}}
\put(5552,-1470){\makebox(0,0)[lb]{\smash{{\SetFigFont{10}{12.0}{\familydefault}{\mddefault}{\updefault}{\color[rgb]{0,0,0}${\C'}^\NN$}%
}}}}
\put(4459,-1546){\makebox(0,0)[lb]{\smash{{\SetFigFont{10}{12.0}{\familydefault}{\mddefault}{\updefault}{\color[rgb]{0,0,0}$h({\sf U})={\sf U}'$}%
}}}}
\put(902,-1414){\makebox(0,0)[lb]{\smash{{\SetFigFont{10}{12.0}{\familydefault}{\mddefault}{\updefault}{\color[rgb]{0,0,0}$\C^\NN$}%
}}}}
\end{picture}%

\end{center}

In fact, one shall assume that onto its image:
\[
{\sf U}
\,\longrightarrow\,
h({\sf U})
=:
{\sf U}',
\]
the map is a 
$\mathcal{ C}^\kappa$ ($\kappa \geqslant 1$), or $\mathcal{ C}^\infty$, 
or $\mathcal{ C}^\omega$
diffeomorphism.
This means: bijective, homeomorphic, and
with nowhere vanishing Jacobian determinant:
\[
\det\,\Jac_\R(f,g)
=
\det\,
\left(\!
\begin{array}{ccccc}
f_{1,{\sf x}_1} & f_{1,{\sf y}_1} & \cdots &
f_{1,{\sf x}_\NN} & f_{1,{\sf y}_\NN}
\\ 
g_{1,{\sf x}_1} & g_{1,{\sf y}_1} & \cdots &
g_{1,{\sf x}_\NN} & g_{1,{\sf y}_\NN}
\\
\vdots & \vdots & \ddots & \vdots & \vdots
\\
f_{\NN,{\sf x}_1} & f_{\NN,{\sf y}_1} & \cdots &
f_{\NN,{\sf x}_\NN} & f_{\NN,{\sf y}_\NN}
\\ 
g_{\NN,{\sf x}_1} & g_{\NN,{\sf y}_1} & \cdots &
g_{\NN,{\sf x}_\NN} & g_{\NN,{\sf y}_\NN}
\end{array}
\!\right).
\]

\medskip\noindent{\bf Lemma.}
{\em A $\mathcal{ C}^\kappa$ ($\kappa \geqslant 1$), 
or $\mathcal{ C}^\infty$, 
or $\mathcal{ C}^\omega$
diffeomorphism:}
\[
\R^{2\NN}
\supset
{\sf U}
\overset{\sim}{\,\longrightarrow\,}
{\sf U}'
\subset{\R'}^{2\NN}
\]
{\em written as:}
\[
\aligned
({\sf x}_\bullet,{\sf y}_\bullet)
&
\,\longmapsto\,
\big(
f_\bullet({\sf x}_\bullet,{\sf y}_\bullet),\
g_\bullet({\sf x}_\bullet,{\sf y}_\bullet)
\big)
\\
&
\ \ \ \ \
=
({\sf x}_\bullet',\,{\sf y}_\bullet')
\endaligned
\]
{\em is a biholomorphism if and only if (definition):}
\[
h_k
:=
f_k
+
\isqrt\,g_k
\ \ \ \ \ \ \ \ \ \ \ \ \
{\scriptstyle{(1\,\leqslant\,k\,\leqslant\,\NN)}}
\]
{\em are holomorphic with respect to ${\sf z}_1, \dots, {\sf z}_\NN$,
or equivalently (Cauchy-Riemann equations):}
\[
\boxed{\,\,
\aligned
0
&
\equiv
f_{k,{\sf x}_l}
-
g_{k,{\sf y}_l}
\ \ \ \ \ \ \ \ \ \ \ \ \
{\scriptstyle{(1\,\leqslant\,k,\,l\,\leqslant\,\NN)}},\,\,
\\
0
&
\equiv
f_{k,{\sf y}_l}
+
g_{k,{\sf x}_l}
\ \ \ \ \ \ \ \ \ \ \ \ \
{\scriptstyle{(1\,\leqslant\,k,\,l\,\leqslant\,\NN)}}.
\endaligned}
\qed
\]

\medskip

Since the $\text{\sc n}$ holomorphic
components $h_1, \dots, h_\NN$ are then locally
expandable in converging power series
in $({\sf z}_1, \dots, {\sf z}_\NN)$ and
{\em do not depend on $\big( \overline{\sf z}_1,
\dots, \overline{\sf z}_\NN\big)$}, 
one can also introduce the {\sl holomorphic Jacobian matrix}:
\[
\Jac_\C(h)
:=
\left(\!
\begin{array}{ccc}
h_{1,{\sf z}_1} & \cdots & h_{1, {\sf z}_\NN}
\\
\vdots & \ddots & \vdots
\\
h_{\NN,{\sf z}_1} & \cdots & h_{\NN, {\sf z}_\NN}
\end{array}
\!\right).
\]

By assumption:
\[
\aligned
0
&
\neq
\det\,\Jac_\C(h)
\\
&
=
\det\,
\left(\!
\begin{array}{ccc}
h_{1,{\sf z}_1} & \cdots & h_{1,{\sf z}_\NN}
\\
\vdots & \ddots & \cdots
\\
h_{\NN,{\sf z}_1} & \cdots & h_{\NN,{\sf z}_\NN}
\end{array}
\!\right),
\endaligned
\]
at every point $q \in {\sf U}_p$, or equivalently:
\[
\aligned
0
&
\neq
\det\,\Jac_\R(f,g)
\\
&
=
\det\,
\left(\!
\begin{array}{ccccc}
f_{1,{\sf x}_1} & f_{1,{\sf y}_1} & \cdots &
f_{1,{\sf x}_\NN} & f_{1,{\sf y}_\NN}
\\ 
g_{1,{\sf x}_1} & g_{1,{\sf y}_1} & \cdots &
g_{1,{\sf x}_\NN} & g_{1,{\sf y}_\NN}
\\
\vdots & \vdots & \ddots & \vdots & \vdots
\\
f_{\NN,{\sf x}_1} & f_{\NN,{\sf y}_1} & \cdots &
f_{\NN,{\sf x}_\NN} & f_{\NN,{\sf y}_\NN}
\\ 
g_{\NN,{\sf x}_1} & g_{\NN,{\sf y}_1} & \cdots &
g_{\NN,{\sf x}_\NN} & g_{\NN,{\sf y}_\NN}
\end{array}
\!\right),
\endaligned
\]
because of an:

\medskip\noindent{\bf Exercise.}
\[
\det\,\Jac_\R(f,g)
=
\big\vert
\det\,\Jac_\C(h)
\big\vert^2.
\qed
\]

However, even when $h$ is a biholomorphism, it is preferable
to work mainly with its real Jacobian matrix.

Coming back to a 
$\mathcal{ C}^\kappa$ ($\kappa \geqslant 1$), or $\mathcal{ C}^\infty$, 
or $\mathcal{ C}^\omega$
diffeomorphism $(f, g)$, the
push-forwards of basic vector fields are (exercise):
\[
\aligned
(f,g)_*
\bigg(
\frac{\partial}{\partial{\sf x}_l}
\bigg)
&
=
\sum_{k=1}^\NN\,
\bigg(
f_{k,{\sf x}_l}\,
\frac{\partial}{\partial{\sf x}_k'}
+
g_{k,{\sf x}_l}\,
\frac{\partial}{\partial{\sf y}_k'}
\bigg)
\ \ \ \ \ \ \ \ \ \ \ \ \
{\scriptstyle{(l\,=\,1\,\cdots\,\NN)}},
\\
(f,g)_*
\bigg(
\frac{\partial}{\partial{\sf y}_l}
\bigg)
&
=
\sum_{k=1}^\NN\,
\bigg(
f_{k,{\sf y}_l}\,
\frac{\partial}{\partial{\sf x}_k'}
+
g_{k,{\sf y}_l}\,
\frac{\partial}{\partial{\sf y}_k'}
\bigg)
\ \ \ \ \ \ \ \ \ \ \ \ \
{\scriptstyle{(l\,=\,1\,\cdots\,\NN)}}.
\endaligned
\]
{\em Now by composing
with $h^{ -1}$, one can insure that in the right-hand side,
everything lives in the $( {\sf x}_\bullet', {\sf y}_\bullet')$-space:}
\[
\boxed{
\footnotesize
\aligned
(f,g)_*
\bigg(
\frac{\partial}{\partial{\sf x}_l}
\bigg)
&
\overset{\sf def}{\,:=\,}
\sum_{k=1}^\NN\,
\bigg(
f_{k,{\sf x}_l}
\circ
h^{-1}\big({\sf x}_\bullet',{\sf y}_\bullet'\big)\,
\frac{\partial}{\partial{\sf x}_k'}
+
g_{k,{\sf x}_l}
\circ
h^{-1}\big({\sf x}_\bullet',{\sf y}_\bullet'\big)\,
\frac{\partial}{\partial{\sf y}_k'}
\bigg)
\ \ \ \ \ \ \ \ \ \ \ \ \
{\scriptstyle{(l\,=\,1\,\cdots\,\NN)}},
\\
(f,g)_*
\bigg(
\frac{\partial}{\partial{\sf y}_l}
\bigg)
&
\overset{\sf def}{\,:=\,}
\sum_{k=1}^\NN\,
\bigg(
f_{k,{\sf y}_l}
\circ
h^{-1}\big({\sf x}_\bullet',{\sf y}_\bullet'\big)\,
\frac{\partial}{\partial{\sf x}_k'}
+
g_{k,{\sf y}_l}
\circ
h^{-1}\big({\sf x}_\bullet',{\sf y}_\bullet'\big)\,
\frac{\partial}{\partial{\sf y}_k'}
\bigg)
\ \ \ \ \ \ \ \ \ \ \ \ \
{\scriptstyle{(l\,=\,1\,\cdots\,\NN)}}.
\endaligned}
\]

\medskip\noindent{\bf Transfer of Lie brackets.} 
On an open subset:
\[
{\sf U}
\subset
\R^{2N},
\]
consider two $\mathcal{ C}^0$ vector field sections:
\[
\aligned
{\sf P}
&
=
\sum_{k=1}^{\NN}\,
\bigg(
a_k({\sf x}_\bullet,{\sf y}_\bullet)\,
\frac{\partial}{\partial{\sf x}_k}
+
b_k({\sf x}_\bullet,{\sf y}_\bullet)\,
\frac{\partial}{\partial{\sf y}_k}
\bigg),
\\
{\sf Q}
&
=
\sum_{k=1}^{\NN}\,
\bigg(
c_k({\sf x}_\bullet,{\sf y}_\bullet)\,
\frac{\partial}{\partial{\sf x}_k}
+
d_k({\sf x}_\bullet,{\sf y}_\bullet)\,
\frac{\partial}{\partial{\sf y}_k}
\bigg).
\endaligned
\]

\medskip\noindent{\bf Definition.}
The {\sl Lie bracket}:
\[
\big[{\sf P},\,{\sf Q}\big]
\]
between two such general vector fields is the {\em vector field}:
\[
\aligned
\big[{\sf P},\,{\sf Q}\big]
&
\,:=\,
\sum_{k=1}^{\NN}\,
\bigg(
\sum_{l=1}^\NN\,
\Big(
a_l\,c_{k,{\sf x}_l}
+
b_l\,c_{k,{\sf y}_l}
-
c_l\,a_{k,{\sf x}_l}
-
d_l\,a_{k,{\sf y}_l}
\Big)
\bigg)\,
\frac{\partial}{\partial{\sf x}_k}
+
\\
&
\ \ \ \ \
+
\sum_{k=1}^{\NN}\,
\bigg(
\sum_{l=1}^\NN\,
\Big(
a_l\,d_{k,{\sf x}_l}
+
b_l\,d_{k,{\sf y}_l}
-
c_l\,b_{k,{\sf x}_l}
-
d_l\,b_{k,{\sf y}_l}
\Big)
\bigg)\,
\frac{\partial}{\partial{\sf y}_k},
\endaligned
\]
where indices after commas abbreviate partial derivatives
of coefficient-functions.

\medskip\noindent{\bf Lemma.}
{\em Through a $\mathcal{ C}^\kappa$ ($\kappa \geqslant 1$), 
or $\mathcal{ C}^\infty$, 
or $\mathcal{ C}^\omega$
diffeomorphism:}
\[
h\colon\ \ \
{\sf U}
\overset{\sim}{\,\longrightarrow\,}
{\sf U}'
=
h({\sf U})
\subset
{\R'}^{2\NN},
\]
{\em Lie brackets between real vector fields transfer as}:
\[
h_*\big(\big[{\sf P},{\sf Q}\big]\big)
=
\big[h_*({\sf P}),h_*({\sf Q})\big].
\]

\proof
Considered to be known or can be reproved directly
by applying the formulas. This, in particular, 
assures that the Lie bracket is well defined, 
independently of coordinates.
\endproof

\medskip\noindent{\bf Biholomorphisms commute with complex
structures.}
Next, the complex structure $J'$ on ${\C'}^\NN$ acts as: 
\[
J'
\bigg(
\frac{\partial}{\partial{\sf x}_k'}
\bigg)
=
\frac{\partial}{\partial{\sf y}_k'},
\ \ \ \ \ \ \ \ \ \ \ \ \ \ \ \ \ \ \ \ \ \ \ \
J'
\bigg(
\frac{\partial}{\partial{\sf y}_k'}
\bigg)
=
-\,\frac{\partial}{\partial{\sf x}_k'},
\ \ \ \ \ \ \ \ \ \ \ \ \
{\scriptstyle{(k\,=\,1\,\cdots\,\NN)}}.
\]

\medskip\noindent{\bf Proposition.}
{\em
A 
$\mathcal{ C}^\kappa$ ($\kappa \geqslant 1$), or $\mathcal{ C}^\infty$, 
or $\mathcal{ C}^\omega$
diffeomorphism:}
\[
h
=
(f,g)
\colon\ \ \
\xymatrix{
{\sf U}
\rule[-3pt]{0pt}{15pt} 
\ar[r]^\sim 
\ar@{_{(}->}[d]
& 
{\sf U}'
\rule[-3pt]{0pt}{15pt} 
\ar@{_{(}->}[d]
\\
\R^{2\NN}=\C^\NN
&
{\C'}^\NN={\R'}^{2\NN}
}
\]
{\em between two open sets is a {\em biholomorphism} if
and only if:}
\[
\boxed{\,\,
h_*\circ J
=
J'\circ h_*,\,\,}
\]
{\em which means the diagram commutation:}
\[
\xymatrix{
T\R^{2\NN}
\ar[r]^{h_*} 
\ar[d]_J
& 
T{\R'}^{2\NN}
\rule[-3pt]{0pt}{15pt} 
\ar[d]^{J'}
\\
T\R^{2\NN}
\ar[r]_{h_*}
&
T{\R'}^{2\NN},
}
\]
{\em on restriction to the concerned open subsets.}

\proof
By linearity of $h_*$, $J$, $J'$, it suffices to show:
\[
\aligned
h_*\circ J\,
\bigg(
\frac{\partial}{\partial{\sf x}_l}
\bigg)
&
\overset{\text{\bf ?}}{\,=\,}
J'\circ h_*
\bigg(
\frac{\partial}{\partial{\sf x}_l}
\bigg),
\\
h_*\circ J\,
\bigg(
\frac{\partial}{\partial{\sf y}_l}
\bigg)
&
\overset{\text{\bf ?}}{\,=\,}
J'\circ h_*
\bigg(
\frac{\partial}{\partial{\sf y}_l}
\bigg),
\endaligned
\]
for $l = 1, \dots, \NNN$. Of course:
\[
h_*
=
(f,g)_*,
\]
in the real sense.

One computes:
\[
\footnotesize
\aligned
h_*
\bigg(
J
\bigg(
\frac{\partial}{\partial{\sf x}_l}
\bigg)
\bigg)
&
=
h_*
\bigg(
\frac{\partial}{\partial{\sf y}_l}
\bigg)
\\
&
=
\sum_{k=1}^\NN\,
\bigg(
f_{k,{\sf y}_l}\,
\frac{\partial}{\partial{\sf x}_k'}
+
g_{k,{\sf y}_l}\,
\frac{\partial}{\partial{\sf y}_k'}
\bigg)
\\
\explain{Insert $J'$}
\ \ \ \ \ \ \ \ \ \ \ \ \ \ \ \ \ \ \ \ \ \ \
&
=
\sum_{k=1}^\NN\,
\bigg(
-\,f_{k,{\sf y}_l}\,
J'\bigg(\frac{\partial}{\partial{\sf y}_k'}\bigg)
+
g_{k,{\sf y}_l}\,
J'\bigg(\frac{\partial}{\partial{\sf x}_k'}\bigg)
\bigg)
\\
\explain{Cauchy-Riemann equations}
\ \ \ \ \ \ \ \ \ \ \ \ \ \ \ \ \ \ \ \ \ \ \
&
=
\sum_{k=1}^\NN\,
\bigg(
g_{k,{\sf x}_l}\,
J'\bigg(\frac{\partial}{\partial{\sf y}_k'}\bigg)
+
f_{k,{\sf x}_l}\,
J'\bigg(\frac{\partial}{\partial{\sf x}_k'}\bigg)
\bigg)
\\
\explain{Extract $J'$}
\ \ \ \ \ \ \ \ \ \ \ \ \ \ \ \ \ \ \ \ \ \ \
&
=
J'
\bigg(
\sum_{k=1}^\NN\,
\bigg(
f_{k,{\sf x}_l}\,
\frac{\partial}{\partial{\sf x}_k'}
+
g_{k,{\sf x}_l}\,
\frac{\partial}{\partial{\sf y}_k'}
\bigg)
\bigg)
\\
\explain{Recognize}
\ \ \ \ \ \ \ \ \ \ \ \ \ \ \ \ \ \ \ \ \ \ \
&
=
J'
\bigg(
h_*
\bigg(
\frac{\partial}{\partial{\sf x}_l}
\bigg)
\bigg).
\endaligned
\]
The second family of identities is checked similarly.
The converse is just logical.
\endproof

\medskip\noindent{\bf Complexification.}
Define the {\em complexified} tangent vector bundle:
\[
\C\otimes_\R
T\R^{2\NN}
\]
by its fibers at points $p \in \R^{2\NN}$:
\[
\C\otimes T_p\R^{2\NN}
\overset{\sf def}{\,:=\,}
\C\,\frac{\partial}{\partial{\sf x}_1}
\bigg\vert_p
\oplus
\C\,\frac{\partial}{\partial{\sf y}_1}
\bigg\vert_p
\oplus\cdots\oplus
\C\,\frac{\partial}{\partial{\sf x}_\NN}
\bigg\vert_p
\oplus
\C\,\frac{\partial}{\partial{\sf y}_\NN}
\bigg\vert_p.
\]

Introducing the fields:
\[
\aligned
\frac{\partial}{\partial{\sf z}_1}
&
=
\frac{1}{2}\,\frac{\partial}{\partial{\sf x}_1}
-
\frac{\isqrt}{2}\,\frac{\partial}{\partial{\sf y}_1},\,\,
\dots\dots\dots,\,\,
\frac{\partial}{\partial{\sf z}_\NN}
=
\frac{1}{2}\,\frac{\partial}{\partial{\sf x}_\NN}
-
\frac{\isqrt}{2}\,\frac{\partial}{\partial{\sf y}_\NN},
\\
\frac{\partial}{\partial\overline{\sf z}_1}
&
=
\frac{1}{2}\,\frac{\partial}{\partial{\sf x}_1}
+
\frac{\isqrt}{2}\,\frac{\partial}{\partial{\sf y}_1},\,\,
\dots\dots\dots,\,\,
\frac{\partial}{\partial\overline{\sf z}_\NN}
=
\frac{1}{2}\,\frac{\partial}{\partial{\sf x}_\NN}
+
\frac{\isqrt}{2}\,\frac{\partial}{\partial{\sf y}_\NN},
\endaligned
\]
one realizes that equivalently:
\[
\C\otimes T_p\R^{2\NN}
\,=\,
\C\,\frac{\partial}{\partial{\sf z}_1}
\bigg\vert_p
\oplus\cdots\oplus
\C\,\frac{\partial}{\partial{\sf z}_\NN}
\bigg\vert_p
\,\oplus\,
\C\,\frac{\partial}{\partial\overline{\sf z}_1}
\bigg\vert_p
\oplus\cdots\oplus
\C\,\frac{\partial}{\partial\overline{\sf z}_\NN}
\bigg\vert_p.
\]

On an open set:
\[
{\sf U}
\subset
\C^\NN
=
\R^{2\NN},
\]
a $\mathcal{ C}^0$ vector field section of $\C\otimes_\R
T\R^{2\NN}$ writes:
\[
\sum_{k=1}^\NN\,
\bigg(
\alpha_k\big({\sf x}_\bullet,{\sf y}_\bullet\big)\,
\frac{\partial}{\partial{\sf x}_k}
+
\beta_k\big({\sf x}_\bullet,{\sf y}_\bullet\big)\,
\frac{\partial}{\partial{\sf y}_k}
\bigg),
\]
with $\mathcal{ C}^0$ complex-valued functions:
\[
\alpha_k,\ \ 
\beta_k\colon\ \ \ \ \
{\sf U}
\,\longrightarrow\,
\C
\ \ \ \ \ \ \ \ \ \ \ \ \
{\scriptstyle{(k\,=\,1\,\cdots\,\NN)}}.
\]
Replacing:
\[
\frac{\partial}{\partial{\sf x}_k}
=
\frac{\partial}{\partial{\sf z}_k}
+
\frac{\partial}{\partial\overline{\sf z}_k}
\ \ \ \ \ \ \ \ \ \ \ \ \
\text{\rm and}
\ \ \ \ \ \ \ \ \ \ \ \ \
\frac{\partial}{\partial{\sf y}_k}
=
\isqrt\,\bigg(
\frac{\partial}{\partial{\sf z}_k}
-
\frac{\partial}{\partial\overline{\sf z}_k}
\bigg),
\]
such a vector field section of $\C \otimes_\R T \R^{2\NN}$ may 
also be written:
\[
\sum_{k=1}^\NN\,
\bigg(
\widetilde{\alpha}_k\big({\sf x}_\bullet,{\sf y}_\bullet\big)\,
\frac{\partial}{\partial{\sf z}_k}
+
\widetilde{\beta}_k\big({\sf x}_\bullet,{\sf y}_\bullet\big)\,
\frac{\partial}{\partial\overline{\sf z}_k},
\bigg)
\]
with:
\[
\widetilde{\alpha}_k
=
\alpha_k
+
\isqrt\,\beta_k
\ \ \ \ \ \ \ \ \ \ \ \ \
\text{\rm and}
\ \ \ \ \ \ \ \ \ \ \ \ \
\widetilde{\beta}_k
=
\alpha_k
-
\isqrt\,\beta_k
\ \ \ \ \ \ \ \ \ \ \ \ \
{\scriptstyle{(k\,=\,1\,\cdots\,\NN)}}.
\]

\medskip\noindent{\bf Extension of $h_*$.}
Given two {\em real} vector field local sections:
\[
P
\ \ \ \ \ \ \ \ \ \ \ \ \
\text{\rm and}
\ \ \ \ \ \ \ \ \ \ \ \ \
Q
\]
of $T\R^{2\NN}$ over ${\sf U} \subset \R^{2\NN}$ open,
and given a $\mathcal{ C}^\kappa$ ($\kappa \geqslant 1$), 
or $\mathcal{ C}^\infty$, 
or $\mathcal{ C}^\omega$
diffeomorphism:
\[
h
\colon\ \ \
{\sf U}
\overset{\sim}{\,\longrightarrow\,}
{\sf U}',
\]
one extends:
\[
\boxed{\,\,
h_*\big(
P
+
\isqrt\,
Q
\big)
\overset{\sf def}{\,:=\,}
h_*(P)
+
\isqrt\,h_*(Q).\,\,}
\]

\medskip\noindent{\bf Lemma.}
{\em For any two {\em complex} 
vector field (local) sections of $\C\otimes_\R T\R^{2\NN}$:}
\[
\aligned
\mathcal{P}
&
=
P^r+\isqrt\,P^i,
\\
\mathcal{Q}
&
=
Q^r+\isqrt\,Q^i,
\endaligned
\]
{\em with real vector fields $P^r$, $P^i$, $Q^r$, $Q^i$, one has:}
\[
h_*\big([\mathcal{P},\,\mathcal{Q}]\big)
=
\big[h_*(\mathcal{P}),\,h_*(\mathcal{Q})\big].
\]

\proof
It suffices to compute:
\[
\footnotesize
\aligned
h_*\big([\mathcal{P},\,\mathcal{Q}]\big)
&
=
h_*\Big(
\big[P^r+\isqrt\,P^i,\,Q^r+\isqrt\,Q^i\big]
\Big)
\\
&
=
h_*
\Big(
\big[P^r,Q^r\big]
-
\big[P^i,Q^i\big]
+
\isqrt\,\big[P^r,Q^i\big]
+
\isqrt\,\big[P^i,Q^r\big]
\Big)
\\
&
=
h_*\big(\big[P^r,Q^r\big]\big)
-
h_*\big(\big[P^i,Q^i\big]\big)
+
\isqrt\,
h_*\big(\big[P^r,Q^i\big]\big)
+
\isqrt\,
h_*\big(\big[P^i,Q^r\big]\big)
\\
&
=
\big[h_*(P^r),h_*(Q^r)\big]
-
\big[h_*(P^i),h_*(Q^i)\big]
+
\isqrt\,
\big[h_*(P^r),h_*(Q^i)\big]
+
\isqrt\,
\big[h_*(P^i),h_*(Q^r)\big]
\\
&
=
\Big[
h_*(P^r)+\isqrt\,h_*(P^i),\,\,
h_*(Q^r)+\isqrt\,h_*(Q^i)
\Big]
\\
&
=
\big[
h_*(P^r+\isqrt\,P^i),\,
h_*(Q^r+\isqrt\,Q^i)
\big]
\\
&
=
\big[h_*(\mathcal{P}),\,h_*(\mathcal{Q})\big],
\endaligned
\]
as was easy to check.
\endproof

\medskip\noindent{\bf Lemma.}
{\em For any local section of $\C \otimes_\R T \C^\NN$:}
\[
\mathcal{P}
=
P^r+\isqrt\,P^i,
\]
{\em with real vector fields $P^r$, $P^i$, one has:}
\[
\overline{h_*(\mathcal{P})}
=
h_*
\big(
\overline{\mathcal{P}}\big).
\]

\proof
Again elementarily:
\[
\aligned
\overline{h_*(P^r+\isqrt\,P^i)}
&
=
\overline{
h_*(P^r)+\isqrt\,h_*(P^i)}
\\
&
=
h_*(P^r)
-
\isqrt\,h_*(P^i)
\\
&
=
h_*(P^r-\isqrt\,P^i)
\\
&
=
h_*(\mathcal{P}),
\endaligned
\]
which is so.
\endproof

Now, the quite central concept of $(1, 0)$ and
of $(0, 1)$ bundles enters the scene.

\noindent{\bf Definition.}
For $p \in \C^\NN = \R^{ 2\NN}$, set:
\[
\aligned
T_p^{1,0}\R^{2\NN}
&
\overset{\sf def}{\,:=\,}
\big\{
X_p-\isqrt\,J(X_p)
\colon
X_p\in T_p\R^{2\NN}
\big\}
\\
&\,\,\,
=
\C\,\frac{\partial}{\partial{\sf z}_1}
\bigg\vert_p
\oplus\cdots\oplus
\C\,\frac{\partial}{\partial{\sf z}_\NN}
\bigg\vert_p,
\endaligned
\]
and:
\[
\aligned
T_p^{0,1}\R^{2\NN}
&
\overset{\sf def}{\,:=\,}
\big\{
X_p+\isqrt\,J(X_p)
\colon
X_p\in T_p\R^{2\NN}
\big\}
\\
&\,\,\,
=
\C\,\frac{\partial}{\partial\overline{\sf z}_1}
\bigg\vert_p
\oplus\cdots\oplus
\C\,\frac{\partial}{\partial\overline{\sf z}_\NN}
\bigg\vert_p.
\endaligned
\]

\medskip
One checks:
\[
T^{0,1}\R^{2\NN}
=
\overline{T^{1,0}\R^{2\NN}},
\]
and:
\[
\C\otimes_\R \R^{2\NN}
=
T^{1,0}\R^{2\NN}
\oplus
T^{0,1}\R^{2\NN}.
\]

Clearly, a $\mathcal{ C}^0$ vector field section of 
$T^{1,0} \R^{2\NN}$ writes:
\[
\sum_{k=1}^\NN\,
\alpha_k\big({\sf x}_\bullet,{\sf y}_\bullet\big)\,
\frac{\partial}{\partial{\sf z}_k},
\]
with $\mathcal{ C}^0$ complex-valued functions:
\[
\alpha_k
\colon\ \ \ \ \
{\sf U}
\,\longrightarrow\,
\C
\ \ \ \ \ \ \ \ \ \ \ \ \
{\scriptstyle{(k\,=\,1\,\cdots\,\NN)}},
\]
while a $\mathcal{ C}^0$ vector field section of $T^{0,1} \R^{2\NN}$ writes:
\[
\sum_{k=1}^\NN\,
\beta_k\big({\sf x}_\bullet,{\sf y}_\bullet\big)\,
\frac{\partial}{\partial\overline{\sf z}_k},
\]
with $\mathcal{ C}^0$ complex-valued functions:
\[
\beta_k
\colon\ \ \ \ \
{\sf U}
\,\longrightarrow\,
\C
\ \ \ \ \ \ \ \ \ \ \ \ \
{\scriptstyle{(k\,=\,1\,\cdots\,\NN)}}.
\]

\medskip\noindent{\bf Lemma.}
{\em Through a biholomorphism:}
\[
h
\colon\ \ \
\xymatrix{
{\sf U}
\rule[-3pt]{0pt}{15pt} 
\ar[r]^\sim 
\ar@{_{(}->}[d]
& 
{\sf U}'
\rule[-3pt]{0pt}{15pt} 
\ar@{_{(}->}[d]
\\
\R^{2\NN}=\C^\NN
&
{\C'}^\NN={\R'}^{2\NN}
}
\]
{\em viewed as a real map $h = (f, g)$, one has:}
\[
\aligned
h_*\big(T_p^{1,0}\R^{2\NN}\big)
&
=
T_{h(p)}^{1,0}{\R'}^{2\NN},
\\
h_*\big(T_p^{0,1}\R^{2\NN}\big)
&
=
T_{h(p)}^{0,1}{\R'}^{2\NN},
\endaligned
\]
{\em at every $p \in {\sf U}$.}

\proof
Indeed:
\[
\aligned
h_*\big(X_p
-
\isqrt\,J(X_p)\big)
&
=
h_*(X_p)
-
\isqrt\,h_*\big(J(X_p)\big)
\\
\explain{Apply $h_*\circ J = J'\circ h_*$}
\ \ \ \ \ \ \ \ \ \ \ \ \ \ \ \ \ \ \ \ \ \ \ \ \ \ \ \ \ \ \ \ \ \ \ 
\ \ \ \ \ \ \ \ \ \ \ \ \ \ \ \
&
=
h_*(X_p)
-
\isqrt\,J'\big(h_*(X_p)\big)
\\
&
=:
X_{h(p)}'
-
\isqrt\,J'(X_{h(p)}'\big)
\\
&
\in
T^{1,0}{\R'}^{2\NN},
\endaligned
\]
whence:
\[
h_*
\big(T^{1,0}\R^{2\NN}\big)
\,\subset\,
T^{1,0}{\R'}^{2\NN}.
\]
For the reverse inclusion, proceed with $h^{-1}$. 
For $T^{0, 1}$, conjugate $T^{1, 0}$.
\endproof

\noindent{\bf Notation.}
Given a  
$\mathcal{ C}^\kappa$ ($\kappa \geqslant 1$), or $\mathcal{ C}^\infty$, 
or $\mathcal{ C}^\omega$
diffeomorphism:
\[
h
\colon\ \ \
\xymatrix{
{\sf U}
\rule[-3pt]{0pt}{15pt} 
\ar[r]^\sim 
\ar@{_{(}->}[d]
& 
{\sf U}'
\rule[-3pt]{0pt}{15pt} 
\ar@{_{(}->}[d]
\\
\R^{2\NN}=\C^\NN
&
{\C'}^\NN={\R'}^{2\NN}
}
\]
between two open subsets, use the symbol:
\[
\boxed{\,h_*\,}
\]
to denote both differentials:
\[
\boxed{\,\,
\aligned
h_*\colon
\ \ \ \ \ \ \ \ \ \ \ \ \ \ \,
T\C^\NN
&
\,\longrightarrow\,
T{\C'}^\NN,
\\
h_*\colon
\ \ \ \ \
\C\otimes_\R T\C^\NN
&
\,\longrightarrow\,
\C\otimes_\R T{\C'}^\NN,
\endaligned}
\]
on restriction to the concerned open sets.

\medskip\noindent{\bf Application to CR-generic submanifolds.}
Let:
\[
M^{2n+c}
\,\subset\,
\C^{n+c}
\]
be a connected $\mathcal{ C}^\kappa$ ($\kappa \geqslant 1$), or 
$\mathcal{ C}^\infty$, or $\mathcal{ C}^\omega$ submanifold
which is CR-generic:
\[
TM
+
J(TM)
=
T\C^{n+c}
\big\vert_M,
\]
with:
\[
\aligned
n
&
=
\CRdim\,M,
\\
c
&
=
\codim\,M,
\endaligned
\]
basic facts about CR manifolds having been already developed
in what precedes.

\medskip

At each point $q \in M$, one therefore views {\em extrinsically}:
\[
T_qM
\subset
T_q\C^{n+c}
=
T_q\R^{2n+2c}.
\]
Now, complexify:
\[
\aligned
\C\otimes_\R T_qM
&
\,\subset\,
\C\otimes_\R T_q\C^{n+c}
\\
&
\,\subset\,
T_q^{1,0}\C^{n+c}
\oplus
T_q^{0,1}\C^{n+c},
\endaligned
\]
using:
\[
\C\otimes_\R T_q\C^{n+c}
=
T_q^{1,0}\C^{n+c}
\oplus
T_q^{0,1}\C^{n+c}.
\]

\medskip\noindent{\bf Definition.}
At every $q \in M^{2n+c}$, set:
\[
\aligned
T_q^{1,0}M
&
\overset{\sf def}{\,:=\,}
T_q^{1,0}\C^{n+c}
\cap
\big(
\C\otimes_\R T_qM
\big),
\\
T_q^{0,1}M
&
\overset{\sf def}{\,:=\,}
T_q^{0,1}\C^{n+c}
\cap
\big[
\C\otimes_\R T_qM
\big],
\endaligned
\]
so that:
\[
T_q^{0,1}M
=
\overline{T_q^{1,0}M},
\]

\medskip

Now, an arbitrary vector:
\[
\mathcal{L}_q
\,\in\,
T_q^{1,0}\C^{n+c}
\]
writes:
\[
\mathcal{L}_q
=
L_q
-
\isqrt\,J(L_q),
\]
with a real vector:
\[
L_q
\,\in\,
\C\otimes_\R TM.
\]
If, moreover:
\[
\mathcal{L}_q
\,\in\,
\C\otimes_\R TM,
\]
then clearly:
\[
L_q
\in
T_qM
\ \ \ \ \ \ \ \ \ \ \ \ \
\text{\rm and}
\ \ \ \ \ \ \ \ \ \ \ \ \
J(L_q)
\in
T_qM,
\]
that is to say:
\[
L_q
\,\in\,
T_qM
\cap
J(T_qM)
=
T_q^cM.
\]

Conversely, for every $L_q \in T_q^cM$, one has:
\[
\aligned
L_q
-
\isqrt\,J(L_q)
&
\,\in\,
\big(\C\otimes_\R T_qM\big)
\cap
T_q^{1,0}\C^{n+c}
\\
&
\,=\,
T_q^{1,0}M.
\endaligned
\]

\medskip\noindent{\bf Summary.}
{\em On a CR-generic submanifold $M^{ 2n+c} \subset \C^{n+c}$ of CR
dimension $n$, one has:}
\[
\aligned
T_q^{1,0}M
&
\,=\,
\Big\{
L_q-\isqrt\,J(L_q)\colon\,
L_q\in T_qM\cap J(T_qM)
\Big\},
\\
T_q^{0,1}M
&
\,=\,
\Big\{
L_q+\isqrt\,J(L_q)\colon\,
L_q\in T_qM\cap J(T_qM)
\Big\},
\endaligned
\]
{\em and as $q \in M$ runs, these spaces gather coherently 
to constitute two $\C$-vector bundles of ranks:}
\[
\aligned
\rank_\C\big(T^{1,0}M\big)
=
\rank_\C\big(T^{0,1}M\big)
=
{\textstyle{\frac{1}{2}}}\,
\rank_\R
\big(TM\cap J(TM)\big).
\endaligned
\]

\medskip

As bundles, one may also write:
\[
\aligned
T^{1,0}M
&
\overset{\sf def}{\,:=\,}
\big(T^{1,0}\C^{n+c}
\big\vert_M
\big)
\cap
\big(
\C\otimes_\R TM
\big),
\\
T^{0,1}M
&
\overset{\sf def}{\,:=\,}
\big(T^{0,1}\C^{n+c}
\big\vert_M
\big)
\cap
\big(
\C\otimes_\R TM
\big).
\endaligned
\]

\medskip\noindent{\bf Fundamental Proposition.}
{\em On a CR-generic $M^{ 2n+c} \subset \C^{n+c}$ of CR dimension $n$,
both bundles $T^{ 1, 0} M$ and $T^{ 0, 1} M$ are
Frobenius-integrable:}
\[
\aligned
\big[T^{1,0}M,\,T^{1,0}M\big]
&
\,\subset\,
T^{1,0}M,
\\
\big[T^{0,1}M,\,T^{0,1}M\big]
&
\,\subset\,
T^{0,1}M,
\endaligned
\]
{\em which means that for any two (local) vector field sections 
$\mathcal{ M}$ and $\mathcal{ N}$ of
$T^{1, 0} M$ on some open subset of $M$, the two Lie brackets:}
\[
\aligned
&
\big[\mathcal{M},\,\mathcal{N}\big],
\\
&
\big[\overline{\mathcal{M}},\,\overline{\mathcal{N}}\big],
\endaligned
\]
{\em are again vector field sections of $T^{1, 0} M$ and
of $T^{0, 1}M$, respectively.}

\medskip

Of course as a difference, it is almost always the case in CR geometry that:
\[
\big[T^{1,0}M,\,T^{0,1}M\big]
\,\not\subset\,
T^{1,0}M
\oplus
T^{0,1}M.
\]

\proof
In brief (more explanations follow):

\medskip\noindent{\bf (i)}\,\,
Lie brackets of sections of $\C \otimes_\R TM$ remain sections of $\C
\otimes_\R TM$, because tangency to a submanifold is preserved after
taking brackets.

\medskip\noindent{\bf (ii)}\,\,
Lie brackets of sections of $T^{1, 0} \C^{ n+c} \big\vert_M$
remain sections of $T^{1, 0} \C^{ n+c} \big\vert_M$ because
linear combinations of $\frac{ \partial}{ \partial{\sf z}_1},
\dots, \frac{ \partial}{ \partial {\sf z}_{ n+c}}$ remain  
such after taking brackets.

\medskip
Consequently, Lie brackets of the {\em intersection}:
\[
T^{1,0}M
=
\big(T^{1,0}\C^{n+c}\big\vert_M\big)
\cap
\big(\C\otimes_\R TM\big)
\]
remain in {\em this} intersection, 
and similarly of course\,\,---\,\,alternatively,
use plain conjugation\,\,---\,\,for
$[ T^{ 0, 1} M, \, T^{ 0, 1} M ] \subset T^{ 0, 1} M$.

\medskip

One can produce an abstract proof of {\bf (i)} and {\bf (ii)}
following known differential-geometric lines, 
but granted the objectives of the
present memoir, it is preferable to introduce now local
coordinates in order to start making everything more explicit.

\medskip

As was already seen above, at every point $p \in M$, there
exist local coordinates:
\[
\big(z_1,\dots,z_n,w_1,\dots,w_c\big)
=
\Big(
x_1+\isqrt\,y_1,\dots,x_n+\isqrt\,y_n,\,
u_1+\isqrt\,v_1,\dots,u_c+\isqrt\,v_c
\Big)
\]
in which:
\[
T_pM
=
\big\{
0
=
v_1
=\cdots=
v_c
\big\},
\]
and there exist $c$ graphing functions:
\[
\varphi_1,\,\dots,\,\varphi_c
\]
that locally represent $M$ as:
\[
\left[
\aligned
v_1
&
=
\varphi_1
\big(x_1,\dots,x_n,y_1,\dots,y_n,u_1,\dots,u_c\big),
\\
\cdots
&
\cdots\cdots\cdots\cdots\cdots\cdots\cdots\cdots\cdots\cdots\cdots\cdots
\cdot\cdot
\\ 
v_c
&
=
\varphi_c
\big(x_1,\dots,x_n,y_1,\dots,y_n,u_1,\dots,u_c\big).
\endaligned\right.
\]
and which satisfy:
\[
\aligned
0
&
=
\varphi_1(0)
=
d\varphi_1(0),
\\
\cdots
&
\cdots\cdots\cdots\cdots\cdots\cdots
\\
0
&
=
\varphi_c(0)
=
d\varphi_c(0),
\endaligned
\]

In such coordinates, some $2n + c$ vector fields tangent to $M$
constituting a frame for $TM$ are:
\[
\aligned
X_1
&
=
\frac{\partial}{\partial x_1}
+
\varphi_{1,x_1}\,\frac{\partial}{\partial v_1}
+\cdots+
\varphi_{c,x_1}\,\frac{\partial}{\partial v_c},
\\
\cdots
&
\cdots\cdots\cdots\cdots\cdots\cdots\cdots\cdots\cdots\cdots\cdots\cdots
\\
X_n
&
=
\frac{\partial}{\partial x_n}
+
\varphi_{1,x_n}\,\frac{\partial}{\partial v_1}
+\cdots+
\varphi_{c,x_n}\,\frac{\partial}{\partial v_c},
\\
Y_1
&
=
\frac{\partial}{\partial y_1}
+
\varphi_{1,y_1}\,\frac{\partial}{\partial v_1}
+\cdots+
\varphi_{c,y_1}\,\frac{\partial}{\partial v_c},
\\
\cdots
&
\cdots\cdots\cdots\cdots\cdots\cdots\cdots\cdots\cdots\cdots\cdots\cdots
\\
Y_n
&
=
\frac{\partial}{\partial y_n}
+
\varphi_{1,y_n}\,\frac{\partial}{\partial v_1}
+\cdots+
\varphi_{c,y_n}\,\frac{\partial}{\partial v_c},
\\
U_1
&
=
\frac{\partial}{\partial u_1}
+
\varphi_{1,u_1}\,\frac{\partial}{\partial v_1}
+\cdots+
\varphi_{c,u_1}\,\frac{\partial}{\partial v_c},
\\
\cdots
&
\cdots\cdots\cdots\cdots\cdots\cdots\cdots\cdots\cdots\cdots\cdots\cdots
\\
U_c
&
=
\frac{\partial}{\partial u_1}
+
\varphi_{1,u_c}\,\frac{\partial}{\partial v_1}
+\cdots+
\varphi_{c,u_c}\,\frac{\partial}{\partial v_c}
\endaligned
\]
{\em Importantly}, one notices that these generators:
\[
\big\{
X_1,\dots,X_n,\,Y_1,\dots,Y_n,\,U_1,\dots,U_c
\big\}
\] 
are {\em extrinsic} in the sense that they all involve
partial derivatives:
\[
\frac{\partial}{\partial v_1},\,\,
\dots\dots,\,\,
\frac{\partial}{\partial v_c}
\]
that are {\em not} related to the {\em intrinsic} coordinates on $M$:
\[
\big(x_1,\dots,x_n,\,y_1,\dots,y_n,\,u_1,\dots,u_c\big)
=
\big(x_\bullet,y_\bullet,u_\bullet\big),
\]
so that the {\em geometric} vectors
that these $2n + c$ derivations represent truly live in $\C^{ n+c}$,
while their coefficient-functions:
\[
\varphi_{j,x_k}\big(x_\bullet,y_\bullet,u_\bullet\big),
\ \ \ \ \
\varphi_{j,y_k}\big(x_\bullet,y_\bullet,u_\bullet\big),
\ \ \ \ \
\varphi_{j,u_l}\big(x_\bullet,y_\bullet,u_\bullet\big),
\]
depend only upon the {\em horizontal, intrinsic} coordinates
of $M$.

\medskip

To think intrinsically, then, one introduces the projection:
\[
\aligned
\pi^{2n+c}
\colon
\ \ \ \ \ \ \ \ \ \ \ \ \ \ \ \
\R^{2n+2c}
&
\,\longrightarrow\,
\R^{2n+c}
\\
\big(x_\bullet,y_\bullet,u_\bullet,v_\bullet\big)
&
\,\longmapsto\,
\big(x_\bullet,y_\bullet,u_\bullet\big)
\endaligned
\]
which sends the $2n+c$ elements of the extrinsic frame:
\[
\aligned
\pi_*
(X_k)
&
=
\pi_*
\bigg(
\frac{\partial}{\partial x_k}
+
\sum_{j=1}^c\,\varphi_{j,x_k}\,
\frac{\partial}{\partial v_j}
\bigg)
\\
&
=
\frac{\partial}{\partial x_k},
\\
\pi_*
(Y_k)
&
=
\pi_*
\bigg(
\frac{\partial}{\partial y_k}
+
\sum_{j=1}^c\,\varphi_{j,y_k}\,
\frac{\partial}{\partial v_j}
\bigg)
\\
&
=
\frac{\partial}{\partial y_k},
\\
\pi_*
(U_l)
&
=
\pi_*
\bigg(
\frac{\partial}{\partial u_l}
+
\sum_{j=1}^c\,\varphi_{j,u_l}\,
\frac{\partial}{\partial v_j}
\bigg)
\\
&
=
\frac{\partial}{\partial u_l},
\endaligned
\]
to the straight intrinsic frame on $TM$ naturally associated to
$(x_\bullet, y_\bullet, u_\bullet)$.

\medskip

Now, seek $n$ local generators of $T^{1, 0}M$.
Because at the origin:
\[
T_0^{1,0}M
\,=\,
\frac{\partial}{\partial z_1}
\bigg\vert_0
\oplus\cdots\oplus
\frac{\partial}{\partial z_n}
\bigg\vert_0,
\]
it is natural to seek them under the form:
\[
\aligned
\mathcal{L}_1
&
=
\frac{\partial}{\partial z_1}
+
\sum_{j=1}^c\,
{\tt A}_1^j\,
\frac{\partial}{\partial w_j},
\\
\cdots
&
\cdots\cdots\cdots\cdots\cdots\cdots\cdots
\\
\mathcal{L}_n
&
=
\frac{\partial}{\partial z_n}
+
\sum_{j=1}^c\,
{\tt A}_n^j\,
\frac{\partial}{\partial w_j},
\endaligned
\]
with unknown functions:
\[
{\tt A}_k^j
=
{\tt A}_k^j
\big(
x_\bullet,y_\bullet,u_\bullet
\big).
\]
Again, one notices that such $(1, 0)$ fields live {\em in $\C^{ n+c}$},
while their coefficients depend only on intrinsic coordinates
of $M \cong \R^{ 2n + c}$.

By definition, such fields should be tangent to $M$ in order
to be sections of $T^{ 1, 0} M$.
Writing then the equations of $M$ under the adapted form:
\[
\aligned
0
&
=
-\,v_1
+
\varphi_1\big(x_\bullet,y_\bullet,u_\bullet\big),
\\
\cdot\,\cdot
&
\cdots\cdots\cdots\cdots\cdots\cdots\cdots\cdot\cdot
\\
0
&
=
-\,v_1
+
\varphi_1\big(x_\bullet,y_\bullet,u_\bullet\big),
\endaligned
\] 
and writing simultaneously the fields under the expanded form:
\[
\aligned
\mathcal{L}_i
&
=
\frac{\partial}{\partial z_i}
+
{\tt A}_i^1\,
\bigg(
\frac{1}{2}\,\frac{\partial}{\partial u_1}
-
\frac{\isqrt}{2}\,
\frac{\partial}{\partial v_1}
\bigg)
+\cdots\cdots+
{\tt A}_i^c\,
\bigg(
\frac{1}{2}\,\frac{\partial}{\partial u_c}
-
\frac{\isqrt}{2}\,
\frac{\partial}{\partial v_c}
\bigg)
\\
&
\ \ \ \ \ \ \ \ \ \ \ \ \ \ \ \ \ \ \ \ \ \ \ \ \ \ \ \ \ \ \ \
\ \ \ \ \ \ \ \ \ \ \ \ \ \ \ \ \ \ \ \ \ \ \ \ \ \ \ \ 
{\scriptstyle{(i\,=\,1\,\cdots\,n)}},
\endaligned
\]
one expresses the tangency of $\mathcal{ L}_i$ by applying it considered
as a derivation to the $c$ equations of $M$ and the result should be zero:
\[
\aligned
0
&
=
{\textstyle{\frac{\isqrt}{2}}}\,
{\tt A}_i^1
+
\varphi_{1,z_i}
+
{\tt A}_i^1\,
\big(
{\textstyle{\frac{1}{2}}}\,
\varphi_{1,u_1}
\big)
+\cdots\cdots+
{\tt A}_i^c\,
\big(
{\textstyle{\frac{1}{2}}}\,
\varphi_{1,u_c}
\big),
\\
\cdot\,\cdot
&
\cdots\cdots\cdots\cdots\cdots\cdots\cdots\cdots\cdots\cdots\cdots\cdots
\cdots\cdots\cdots\cdots\cdots\cdot
\\
0
&
=
{\textstyle{\frac{\isqrt}{2}}}\,
{\tt A}_i^c
+
\varphi_{c,z_i}
+
{\tt A}_i^1\,
\big(
{\textstyle{\frac{1}{2}}}\,
\varphi_{c,u_1}
\big)
+\cdots\cdots+
{\tt A}_i^c\,
\big(
{\textstyle{\frac{1}{2}}}\,
\varphi_{c,u_c}
\big).
\endaligned
\]
Reorganizing this linear system as:
\[
\aligned
-\,2\,\varphi_{1,z_i}
&
=
{\tt A}_i^1\,
\big(
\isqrt+\varphi_{1,u_1}
\big)
+
{\tt A}_i^2\,
\big(
\varphi_{1,u_2}
\big)
+\cdots\cdots+
{\tt A}_i^c\,
\big(
\varphi_{1,u_c}
\big),
\\
-\,2\,\varphi_{2,z_i}
&
=
{\tt A}_i^1\,
\big(
\varphi_{2,u_1}
\big)
+
{\tt A}_i^2\,
\big(
\isqrt+\varphi_{2,u_2}
\big)
+\cdots\cdots+
{\tt A}_i^c\,
\big(
\varphi_{2,u_c}
\big),
\\
\cdots\cdots\cdot\cdot
&
\cdots\cdots\cdots\cdots\cdots\cdots\cdots\cdots\cdots\cdots\cdots\cdots
\cdots\cdots\cdots\cdots\cdot\cdot
\\
-\,2\,\varphi_{c,z_i}
&
=
{\tt A}_i^1\,
\big(
\varphi_{c,u_1}
\big)
+
{\tt A}_i^2\,
\big(
\varphi_{c,u_2}
\big)
+\cdots\cdots+
{\tt A}_i^c\,
\big(
\isqrt+\varphi_{c,u_c}
\big),
\endaligned
\]
one sees that its determinant:
\[
\aligned
\left\vert\!
\begin{array}{cccc}
\isqrt+\varphi_{1,u_1} & \varphi_{1,u_2} & \cdots & \varphi_{1,u_c}
\\
\varphi_{2,u_1} & \isqrt+\varphi_{2,u_2} & \cdots & \varphi_{2,u_c}
\\
\vdots & \vdots & \ddots & \vdots
\\
\varphi_{c,u_1} & \varphi_{c,u_2} & \cdots & \isqrt+\varphi_{c,u_c}
\end{array}
\!\right\vert
\endaligned
\]
is locally nonvanishing, because at the origin:
\[
\aligned
\left\vert\!
\begin{array}{cccc}
\isqrt & 0 & \cdots & 0
\\
0 & \isqrt & \cdots & 0
\\
\vdots & \vdots & \ddots & \vdots
\\
0 & 0 & \cdots & \isqrt
\end{array}
\!\right\vert
=
\big(\isqrt\big)^n
\neq
0.
\endaligned
\]
An application of Cramer's rule then concludes a fundamental
explicit:

\medskip\noindent{\bf Proposition.}
{\em On a CR-generic submanifold $M^{ 2n + c} \subset \C^{ n+c}$
of smoothness $\mathcal{ C}^\kappa$ $(\kappa \geqslant 1)$, or 
$\mathcal{ C}^\infty$, or $\mathcal{ C}^\omega$ with:}
\[
\aligned
c
&
=
\codim\,M,
\\
n
&
=
\CRdim\,M
\endaligned
\]
{\em which is locally represented in coordinates:}
\[
\big(x_1,\dots,x_n,\,y_1,\dots,y_n,\,u_1,\dots,u_c\big)
\]
{\em as:}
\[
\aligned
v_1
&
=
\varphi_1
\big(x_1,\dots,x_n,y_1,\dots,y_n,u_1,\dots,u_c\big),
\\
\cdots
&
\cdots\cdots\cdots\cdots\cdots\cdots\cdots\cdots\cdots\cdots\cdots\cdots
\cdot\cdot
\\ 
v_c
&
=
\varphi_c
\big(x_1,\dots,x_n,y_1,\dots,y_n,u_1,\dots,u_c\big).
\endaligned
\]
{\em with graphing functions satisfying:}
\[
\aligned
0
&
=
\varphi_1(0)
=
d\varphi_1(0),
\\
\cdots
&
\cdots\cdots\cdots\cdots\cdots\cdots
\\
0
&
=
\varphi_c(0)
=
d\varphi_c(0),
\endaligned
\]
{\em a local frame for $T^{ 1, 0}M$:}
\[
\big\{
\mathcal{L}_1,\dots,\mathcal{L}_n
\big\}
\]
{\em is constituted of the $n$ vector fields:}
\[
\aligned
\mathcal{L}_1
&
=
\frac{\partial}{\partial z_1}
+
{\tt A}_1^1\big(x_\bullet,y_\bullet,u_\bullet\big)\,
\frac{\partial}{\partial w_1}
+\cdots\cdots+
{\tt A}_1^c\big(x_\bullet,y_\bullet,u_\bullet\big)\,
\frac{\partial}{\partial w_c}
\\
\cdots
&
\cdots\cdots\cdots\cdots\cdots\cdots\cdots\cdots\cdots\cdots\cdots\cdots
\cdots\cdots\cdots\cdots\cdots\cdots
\\
\mathcal{L}_n
&
=
\frac{\partial}{\partial z_n}
+
{\tt A}_n^1\big(x_\bullet,y_\bullet,u_\bullet\big)\,
\frac{\partial}{\partial w_1}
+\cdots\cdots+
{\tt A}_n^c\big(x_\bullet,y_\bullet,u_\bullet\big)\,
\frac{\partial}{\partial w_c},
\endaligned
\]
{\em whose coefficient functions are given, for $i = 1, \dots, n$, 
explicitly by:}
\[
\footnotesize
\aligned
\!\!\!\!\!\!\!\!\!\!\!\!\!\!\!\!\!\!\!\!\!
{\tt A}_i^1
\,=\,
\frac{
\left\vert\!\!
\begin{array}{cccc}
-\,2\,\varphi_{1,z_i} & \varphi_{1,u_2} & \cdots & \varphi_{1,u_c}
\\
-\,2\,\varphi_{2,z_i} & \isqrt+\varphi_{2,u_2} & \cdots & \varphi_{2,u_c}
\\
\vdots & \vdots & \ddots & \vdots
\\
-\,2\,\varphi_{c,z_i} & \varphi_{c,u_2} & \cdots & \isqrt+\varphi_{c,u_c}
\end{array}
\!\!\right\vert
}{
\left\vert\!\!
\begin{array}{ccc}
\isqrt+\varphi_{1,u_1} & \cdots & \varphi_{1,u_c}
\\
\varphi_{2,u_1} & \cdots & \varphi_{2,u_c}
\\
\vdots & \ddots & \vdots
\\
\varphi_{c,u_1} & \cdots & \isqrt+\varphi_{c,u_c}
\end{array}
\!\!\right\vert
},\,\,\,
\dots\dots\dots,\,\,\,
{\tt A}_i^c
\,=\,
\frac{
\left\vert\!\!
\begin{array}{cccc}
\isqrt+\varphi_{1,u_1} & \cdots & -\,2\,\varphi_{1,z_i}
\\
\varphi_{2,u_1} & \cdots & -\,2\,\varphi_{2,z_i}
\\
\vdots & \ddots & \vdots
\\
\varphi_{c,u_1} & \cdots & -\,2\,\varphi_{c,z_i}
\end{array}
\!\!\right\vert
}{
\left\vert\!\!
\begin{array}{ccc}
\isqrt+\varphi_{1,u_1} & \cdots & \varphi_{1,u_c}
\\
\varphi_{2,u_1} & \cdots & \varphi_{2,u_c}
\\
\vdots & \ddots & \vdots
\\
\varphi_{c,u_1} & \cdots & \isqrt+\varphi_{c,u_c}
\end{array}
\!\!\right\vert
}.
\qed
\endaligned
\]

\medskip

Come now back to the proof of the penultimate proposition.

\medskip

Allowing the notational coincidences:
\[
\aligned
\big({\sf z}_1,\dots,{\sf z}_n,{\sf z}_{n+1},\dots,{\sf z}_\NN\big)
&
\,\equiv\,
\big(z_1,\dots,z_n,u_1,\dots,u_c\big),
\\
\NNN
&
\,\equiv\,
n+c,
\endaligned
\]
take any two local sections of $T^{ 1, 0} M$:
\[
\mathcal{M}
\ \ \ \ \ \ \ \ \ \ \ \ \
\text{\rm and}
\ \ \ \ \ \ \ \ \ \ \ \ \
\mathcal{N}.
\]
Since they are of type $(1, 0)$, they both write under the form: 
\[
\aligned
\mathcal{L}
&
=
\sum_{k=1}^\NN\,
c_k(x_\bullet,y_\bullet,u_\bullet)\,
\frac{\partial}{\partial{\sf z}_k},
\\
\mathcal{M}
&
=
\sum_{k=1}^\NN\,
d_k(x_\bullet,y_\bullet,u_\bullet)\,
\frac{\partial}{\partial{\sf z}_k},
\endaligned
\]
whence their bracket:
\[
\aligned
\big[\mathcal{M},\mathcal{N}\big]
&
=
\bigg[
\sum_{k=1}^\NN\,c_k\,\frac{\partial}{\partial{\sf z}_k},\,\,
\sum_{k=1}^\NN\,d_k\,\frac{\partial}{\partial{\sf z}_k}
\bigg]
\\
&
=
\sum_{k=1}^\NN\,
\bigg(
\sum_{l=1}^\NN\,
\big(
c_l\,d_{k,z_l}
-
d_l\,c_{k,z_l}
\big)
\bigg)\,
\frac{\partial}{\partial{\sf z}_k}
\endaligned
\]
is visibly still again of type $(1, 0)$. This explains
with more precisions the claim {\bf (ii)} made above.

\medskip

Concerning {\bf (i)}, set:
\[
r_1(x_\bullet,y_\bullet,u_\bullet)
:=
v_1
-
\varphi_1(x_\bullet,y_\bullet,u_\bullet),\,\,
\dots\dots,\,\,
r_c(x_\bullet,y_\bullet,u_\bullet)
:=
v_c
-
\varphi_c(x_\bullet,y_\bullet,u_\bullet),\,\,
\]
so that reminding the notational coincidence:
\[
({\sf x}_\bullet,{\sf y}_\bullet)
\equiv
\big(x_\bullet,u_\bullet,y_\bullet,v_\bullet\big),
\]
$M$ is then represented as the common zero-set:
\[
M
=
\big\{
0
=
r_1({\sf x}_\bullet,{\sf y}_\bullet)
=\cdots=
r_c({\sf x}_\bullet,{\sf y}_\bullet)
\big\}.
\]

\medskip\noindent{\bf Definition, or Property.}
Next, recall that by a standard known conceptionalized
fact of elementary differential geometry, a (real or) complex
vector field:
\[
\mathcal{M}
=
\sum_{k=1}^\NN\,
\bigg(
\alpha_k({\sf x}_\bullet,{\sf y}_\bullet)\,
\frac{\partial}{\partial{\sf x}_k}
+
\beta_k({\sf x}_\bullet,{\sf y}_\bullet)\,
\frac{\partial}{\partial{\sf y}_k}
\bigg)
\]
defined in $\R^{2\NN}$ on some local 
open neighborhood ${\sf U}_p$ of some point $p \in M$
is {\sl tangent} to:
\[
M
=
\big\{
0
=
r_1({\sf x}_\bullet,{\sf y}_\bullet)
=\cdots=
r_c({\sf x}_\bullet,{\sf y}_\bullet)
\big\}
\]
if:
\[
\aligned
\mathcal{M}(r_1)
&
=
0
\ \ \ \ \
\text{\rm on restriction to}\,\,
\big\{0=r_1=\cdots=r_c\big\},
\\
\cdots\cdots
&
\cdots\cdots\cdots\cdots\cdots\cdots\cdots\cdots\cdots
\cdots\cdots\cdots\cdots\cdot
\\
\mathcal{M}(r_c)
&
=
0
\ \ \ \ \
\text{\rm on restriction to}\,\,
\big\{0=r_1=\cdots=r_c\big\}.
\endaligned
\]

\medskip

Classically, using the smoothness of $M$, namely the {\em independency}
of the $c$ differentials:
\[
dr_1,\dots,dr_c,
\]
a so-called {\sl Hadamard lemma} yields then that these
vanishings produce:
\[
\aligned
\mathcal{M}(r_1)
&
=
\function_1^1\,r_1
+\cdots+
\function_1^c\,r_c,
\\
\cdots\cdots
&
\cdots\cdots\cdots\cdots\cdots\cdots\cdots\cdots\cdots\cdots\cdot
\\
\mathcal{M}(r_c)
&
=
\function_c^1\,r_1
+\cdots+
\function_c^c\,r_c.
\endaligned
\]

\medskip

These reminders being done, here are precisions about claim 
{\bf (i)} left above.

\medskip

Take two tangent vector field sections of $\C \otimes_\R TM$:
\[
\mathcal{M}
\ \ \ \ \ \ \ \ \ \ \ \ \
\text{\rm and}
\ \ \ \ \ \ \ \ \ \ \ \ \
\mathcal{N},
\]
so that, simultaneously also:

\[
\aligned
\mathcal{N}(r_1)
&
=
\function_1^1\,r_1
+\cdots+
\function_1^c\,r_c,
\\
\cdots\cdots
&
\cdots\cdots\cdots\cdots\cdots\cdots\cdots\cdots\cdots\cdots\cdot
\\
\mathcal{N}(r_c)
&
=
\function_c^1\,r_1
+\cdots+
\function_c^c\,r_c.
\endaligned
\]

Then by the very definition of the Lie bracket acting as a derivation
on functions:
\[
\big[\mathcal{M},\mathcal{N}\big]
(r_j)
=
\mathcal{M}\big(\mathcal{N}(r_j)\big)
-
\mathcal{N}\big(\mathcal{M}(r_j)\big).
\]
The reader will then easily check that for $j = 1, \dots, c$, one yet has:
\[
\big[\mathcal{M},\mathcal{N}\big]
(r_j)
=
\function_j^1\,r_1
+\cdots+
\function_j^c\,r_c,
\]
which proves that tangency is preserved under taking Lie brackets,
and which concludes the proof of the fundamental proposition.
\endproof

\noindent{\bf Scholium.}
{\em The local frame for $T^{1, 0}M$ of the preceding proposition:}
\[
\aligned
\mathcal{L}_1
&
=
\frac{\partial}{\partial z_1}
+
{\tt A}_1^1\big(x_\bullet,y_\bullet,u_\bullet\big)\,
\frac{\partial}{\partial w_1}
+\cdots\cdots+
{\tt A}_1^c\big(x_\bullet,y_\bullet,u_\bullet\big)\,
\frac{\partial}{\partial w_c}
\\
\cdots
&
\cdots\cdots\cdots\cdots\cdots\cdots\cdots\cdots\cdots\cdots\cdots\cdots
\cdots\cdots\cdots\cdots\cdots\cdots
\\
\mathcal{L}_n
&
=
\frac{\partial}{\partial z_n}
+
{\tt A}_n^1\big(x_\bullet,y_\bullet,u_\bullet\big)\,
\frac{\partial}{\partial w_1}
+\cdots\cdots+
{\tt A}_n^c\big(x_\bullet,y_\bullet,u_\bullet\big)\,
\frac{\partial}{\partial w_c},
\endaligned
\]
{\em whose coefficient functions are given, for $i = 1, \dots, n$, 
explicitly by:}
\[
\footnotesize
\aligned
\!\!\!\!\!\!\!\!\!\!\!\!\!\!\!\!\!\!\!\!\!
{\tt A}_i^1
\,=\,
\frac{
\left\vert\!\!
\begin{array}{cccc}
-\,2\,\varphi_{1,z_i} & \varphi_{1,u_2} & \cdots & \varphi_{1,u_c}
\\
-\,2\,\varphi_{2,z_i} & \isqrt+\varphi_{2,u_2} & \cdots & \varphi_{2,u_c}
\\
\vdots & \vdots & \ddots & \vdots
\\
-\,2\,\varphi_{c,z_i} & \varphi_{c,u_2} & \cdots & \isqrt+\varphi_{c,u_c}
\end{array}
\!\!\right\vert
}{
\left\vert\!\!
\begin{array}{ccc}
\isqrt+\varphi_{1,u_1} & \cdots & \varphi_{1,u_c}
\\
\varphi_{2,u_1} & \cdots & \varphi_{2,u_c}
\\
\vdots & \ddots & \vdots
\\
\varphi_{c,u_1} & \cdots & \isqrt+\varphi_{c,u_c}
\end{array}
\!\!\right\vert
},\,\,\,
\dots\dots\dots,\,\,\,
{\tt A}_i^c
\,=\,
\frac{
\left\vert\!\!
\begin{array}{cccc}
\isqrt+\varphi_{1,u_1} & \cdots & -\,2\,\varphi_{1,z_i}
\\
\varphi_{2,u_1} & \cdots & -\,2\,\varphi_{2,z_i}
\\
\vdots & \ddots & \vdots
\\
\varphi_{c,u_1} & \cdots & -\,2\,\varphi_{c,z_i}
\end{array}
\!\!\right\vert
}{
\left\vert\!\!
\begin{array}{ccc}
\isqrt+\varphi_{1,u_1} & \cdots & \varphi_{1,u_c}
\\
\varphi_{2,u_1} & \cdots & \varphi_{2,u_c}
\\
\vdots & \ddots & \vdots
\\
\varphi_{c,u_1} & \cdots & \isqrt+\varphi_{c,u_c}
\end{array}
\!\!\right\vert
}.
\endaligned
\]
{\em which is closed under taking Lie brackets:}
\[
\big[\mathcal{L}_{i_1},\,\mathcal{L}_{i_2}\big]
\,\equiv\,
0
\ \ \ \ \
\mod\,\big(
\mathcal{L}_1,\dots,\mathcal{L}_n
\big)
\]
{\em in fact even satisfies the better property of being
commutative:}
\[
\big[\mathcal{L}_{i_1},\,\mathcal{L}_{i_2}\big]
=
0,
\]
{\em for $1 \leqslant i_1, \, i_2 \leqslant n$.}

\proof
When one looks at what such a Lie bracket can give:
\[
\big[\mathcal{L}_{i_1},\,\mathcal{L}_{i_2}\big]
=
\bigg[
\frac{\partial}{\partial z_{i_1}}
+
\sum_{j=1}^c\,{\tt A}_{i_1}^j\,
\frac{\partial}{\partial w_j},\,\,
\frac{\partial}{\partial z_{i_2}}
+
\sum_{j=1}^c\,{\tt A}_{i_2}^j\,
\frac{\partial}{\partial w_j}
\bigg],
\]
one realizes that because both coefficient-functions:
\[
1
\ \ 
\text{\rm of}
\ \
\frac{\partial}{\partial z_{i_1}}
\ \ \ \ \ \ \ \ \ \ \ \ \
\text{\rm and}
\ \ \ \ \ \ \ \ \ \ \ \ \
1
\ \ 
\text{\rm of}
\ \
\frac{\partial}{\partial z_{i_2}}
\]
are constant, they {\em disappear} after one derivation, so that:
\[
\big[\mathcal{L}_{i_1},\,\mathcal{L}_{i_2}\big]
=
\mathmotsf{absolutely no}\,\,
\frac{\partial}{\partial z}
+
\sum_{j=1}^c\,
\bigg(
\mathcal{L}_{i_1}\big({\tt A}_{i_2}^j\big)
-
\mathcal{L}_{i_2}\big({\tt A}_{i_1}^j\big)
\bigg)\,
\frac{\partial}{\partial w_j}
\]
without any need to expand more. Hence, because this result
must a linear combination of $\mathcal{ L}_1, \dots, 
\mathcal{ L}_n$ {\em which do truly contain the independent}
$\frac{ \partial}{ \partial z_1}, \dots, \frac{ \partial}{ \partial z_n}$, 
such a linear combination can only be plainly zero.

\medskip

Interestingly, one may also produce a proof of commutation
by direct computation.

\medskip

Restrict to $c = 1$ for simplicity: 
\[
v
=
\varphi(x_\bullet,y_\bullet,u),
\]
otherwise, one would have 
to spend time to set up an appropriate formalism with determinants.
In codimension $c = 1$:
\[
\mathcal{L}_i
=
\frac{\partial}{\partial z_i}
-
\frac{2\,\varphi_{z_i}}{\isqrt+\varphi_u}\,
\frac{\partial}{\partial w},
\]
the coefficient-functions:
\[
\frac{2\,\varphi_{z_i}}{\isqrt+\varphi_u}
=
\frac{2\,\varphi_{z_i}(x_\bullet,y_\bullet,u)}{
\isqrt+\varphi_u(x_\bullet,y_\bullet,u)}
\]
depend only on horizontal variables, although the
$\frac{ \partial}{ \partial w}$-component is {\em not}
intrinsic to $M$.

Then one computes frankly:
\[
\footnotesize
\aligned
\big[\mathcal{L}_{i_1},\,\mathcal{L}_{i_2}\big]
&
=
\bigg[
\frac{\partial}{\partial z_{i_1}}
-
\frac{2\,\varphi_{z_{i_1}}}{\isqrt+\varphi_u}\,
\frac{\partial}{\partial w},\,\,
\frac{\partial}{\partial z_{i_2}}
-
\frac{2\,\varphi_{z_{i_2}}}{\isqrt+\varphi_u}\,
\frac{\partial}{\partial w}
\bigg]
\\
&
=
\bigg(
\frac{-\,2\,\varphi_{z_{i_1}z_{i_2}}}{\isqrt+\varphi_u}
+
\frac{2\,\varphi_{z_{i_2}}\,\varphi_{uz_{i_1}}}{
(\isqrt+\varphi_u)^2}
+
\frac{2\,\varphi_{z_{i_1}}}{\isqrt+\varphi_u}
\bigg[
\frac{2\varphi_{z_{i_2}w}}{\isqrt+\varphi_u}
-
\frac{2\,\varphi_{z_{i_2}}\,\varphi_{uw}}{
(\isqrt+\varphi_u)^2}
\bigg]
\bigg)\,
\frac{\partial}{\partial w}
-
\\
&
\ \ 
-
\bigg(
\mathmotsf{precisely the same expression but after the permutation}\,\,
i_1\,\longleftrightarrow\,i_2
\bigg)\,
\frac{\partial}{\partial w},
\endaligned
\]
which indeed simplifies to $0$, thanks to a mental exercise.
\endproof

\noindent{\bf Intrinsic generators
for the $T^{1, 0}M$ and $T^{0, 1} M$ bundles.}
Introduce the projection:
\[
\pi\colon
\ \ \ 
M
\,\longrightarrow\,
T_0M,
\]
namely:
\[
\aligned
\R^{2n+2c}
&
\,\longrightarrow\,
\R^{2n+c}
\\
\big(x_\bullet,y_\bullet,u_\bullet,v_\bullet\big)
&
\,\longmapsto\,
\big(x_\bullet,y_\bullet,u_\bullet\big).
\endaligned
\]
Restrict it to $M$:
\[
\pi\big\vert_M
\colon\ \ \
M
\,\longrightarrow\,
\R^{2n+c}.
\]

The {\sl extrinsic} generators for $T^{ 1, 0} M$:
\[
\aligned
\mathcal{L}_1
&
=
\frac{\partial}{\partial z_1}
+
{\tt A}_1^1\,
\bigg(
\frac{1}{2}\,\frac{\partial}{\partial u_1}
\zero{-
\frac{\isqrt}{2}\,\frac{\partial}{\partial v_1}}
\bigg)
+\cdots\cdots+
{\tt A}_1^c\,
\bigg(
\frac{1}{2}\,\frac{\partial}{\partial u_c}
\zero{-
\frac{\isqrt}{2}\,\frac{\partial}{\partial v_c}}
\bigg)
\\
\cdots\cdot
&
\cdots\cdots\cdots\cdots\cdots\cdots\cdots\cdots\cdots\cdots\cdots\cdots
\cdots\cdots\cdots\cdots\cdots\cdots\cdots\cdots\cdots
\\
\mathcal{L}_n
&
=
\frac{\partial}{\partial z_n}
+
{\tt A}_n^1\,
\bigg(
\frac{1}{2}\,\frac{\partial}{\partial u_1}
\zero{-
\frac{\isqrt}{2}\,\frac{\partial}{\partial v_1}}
\bigg)
+\cdots\cdots+
{\tt A}_n^c\,
\bigg(
\frac{1}{2}\,\frac{\partial}{\partial u_c}
\zero{-
\frac{\isqrt}{2}\,\frac{\partial}{\partial v_c}}
\bigg)
\endaligned
\]
have then, through $h_*$, image vector fields
in which the $\frac{ \partial}{ \partial v_\bullet}$-components
are suppressed:

\[
\aligned
\pi_*\big(\mathcal{L}_1\big)
&
=
\frac{\partial}{\partial z_1}
+
\frac{{\tt A}_1^1}{2}\,
\frac{\partial}{\partial u_1}
+\cdots\cdots+
\frac{{\tt A}_1^c}{2}\,
\frac{\partial}{\partial u_c},
\\
\cdots\cdots\cdot
&
\cdots\cdots\cdots\cdots\cdots\cdots\cdots\cdots\cdots\cdots\cdots\cdot
\\
\pi_*\big(\mathcal{L}_n\big)
&
=
\frac{\partial}{\partial z_n}
+
\frac{{\tt A}_n^1}{2}\,
\frac{\partial}{\partial u_1}
+\cdots\cdots+
\frac{{\tt A}_n^c}{2}\,
\frac{\partial}{\partial u_c}.
\endaligned
\]

Accordingly\,\,---\,\,mind font differences\,\,---, set:
\[
\boxed{\,\,
\aligned
A_1^1
&
:=
\frac{{\tt A}_1^1}{2},
\,\dots\dots,\,
A_1^c
:=
\frac{{\tt A}_1^c}{2},\,\,
\\
\cdots
&
\cdots\cdots\cdots\cdots\cdots\cdots\cdots\cdots
\\
A_n^1
&
:=
\frac{{\tt A}_n^1}{2},
\,\dots\dots,\,
A_n^c
:=
\frac{{\tt A}_n^c}{2}.\,\,
\endaligned}
\]

The {\sl intrinsic} generators for $T^{1, 0}M$ will be written
plainly:
\[
\aligned
\mathcal{L}_1
&
=
\frac{\partial}{\partial z_1}
+
A_1^1\,\frac{\partial}{\partial u_1}
+\cdots\cdots+
A_1^c\,\frac{\partial}{\partial u_c},
\\
\cdots
&
\cdots\cdots\cdots\cdots\cdots\cdots\cdots\cdots\cdots\cdots\cdots
\cdot\cdot
\\
\mathcal{L}_n
&
=
\frac{\partial}{\partial z_n}
+
A_n^1\,\frac{\partial}{\partial u_1}
+\cdots\cdots+
A_n^c\,\frac{\partial}{\partial u_c},
\endaligned
\]
while those for $T^{0, 1} M$ are their conjugates:
\[
\aligned
\overline{\mathcal{L}}_1
&
=
\frac{\partial}{\partial\overline{z}_1}
+
\overline{A}_1^1\,\frac{\partial}{\partial u_1}
+\cdots\cdots+
\overline{A}_1^c\,\frac{\partial}{\partial u_c},
\\
\cdots
&
\cdots\cdots\cdots\cdots\cdots\cdots\cdots\cdots\cdots\cdots\cdots
\cdot\cdot
\\
\overline{\mathcal{L}}_n
&
=
\frac{\partial}{\partial\overline{z}_n}
+
\overline{A}_n^1\,\frac{\partial}{\partial u_1}
+\cdots\cdots+
\overline{A}_n^c\,\frac{\partial}{\partial u_c}.
\endaligned
\]

\noindent{\bf Real analytic CR functions.}
Consider as above a CR-generic:
\[
\aligned
&
M^{2n+c}
\,\subset\,
\C^{n+c}
\\
&
{\scriptstyle{(c\,=\,{\sf codim}\,M,\,\,\,
n\,=\,{\sf CRdim}\,M)}},
\endaligned
\]

\medskip\noindent{\bf Definition.}
A $\mathcal{ C}^\kappa$ ($\kappa \geqslant 1$), or
$\mathcal{ C}^\infty$, or $\mathcal{ C}^\omega$ function:
\[
f
\colon\ \ \
M
\longrightarrow
\C
\]
is called {\sl Cauchy-Riemann} (CR for short) when:
\[
0
\equiv
\overline{\mathcal{L}}(f),
\]
for every (local) section:
\[
\overline{\mathcal{L}}
\]
of $T^{0,1}M$.

\medskip\noindent{\bf Theorem.}
{\em 
On a $\mathcal{ C}^\omega$ CR-generic $M^{2n+c} \subset \C^{ n+c}$, 
a $\mathcal{ C}^\omega$ function is CR if and only if 
it is the restriction to $M$: 
\[
f
=
F\big\vert_M
\]
of a function
$F$ holomorphic in some neighborhood of $M$.}

\proof
The statement being local, pick $p \in M$, 
choose $p$-centered affine holomorphic coordinates:
\[
\aligned
\big(z_1,\dots,z_n,w_1,\dots,w_c\big)
&
=
\big(z_\bullet,w_\bullet\big)
\\
&
=
\big(z_\bullet,u_\bullet+\isqrt\,v_\bullet\big),
\endaligned
\]
in which $c$ graphing equations for $M$ are:
\[
\aligned
v_1
&
=
\varphi_1\big(z_\bullet,\overline{z}_\bullet,u_\bullet\big),
\\
\cdots
&
\cdots\cdots\cdots\cdots\cdots\cdot
\\
v_c
&
=
\varphi_c\big(z_\bullet,\overline{z}_\bullet,u_\bullet\big),
\endaligned
\]
with:
\[
\aligned
0
&
=
\varphi_1(0)
=
d\varphi_1(0),
\\
\cdots
&
\cdots\cdots\cdots\cdots\cdots\cdots
\\
0
&
=
\varphi_c(0)
=
d\varphi_c(0);
\endaligned
\]
here, since the graphing functions $\varphi_j$ are assumed to be
{\em real analytic}, they are converging power series 
in $(x_\bullet, y_\bullet, u_\bullet)$, but it
will be more appropriate to consider them as converging
power series in $(z_\bullet, \overline{ z}_\bullet, 
u_\bullet)$:
\[
\varphi_j
\big(z_\bullet,\overline{z}_\bullet,u_\bullet\big)
=
\sum_{\alpha_\bullet\in\N^n}\,
\sum_{\beta_\bullet\in\N^n}\,
\sum_{\gamma_\bullet\in\N^c}\,
\underbrace{\varphi_{j,\alpha_\bullet,\beta_\bullet,\gamma_\bullet}
}_{\in\,\C}
\big(z_\bullet\big)^{\alpha_\bullet}\,
\big(\overline{z}_\bullet\big)^{\beta_\bullet}\,
\big(u_\bullet\big)^{\gamma_\bullet},
\]
whose coefficients satisfy a Cauchy-type estimate:
\[
\big\vert
\varphi_{j,\alpha_\bullet,\beta_\bullet,\gamma\bullet}
\big\vert
\,\leqslant\,
\constant\,
\bigg(
\frac{1}{\radius}
\bigg)^{\vert\alpha_\bullet\vert
+
\vert\beta_\bullet\vert
+
\vert\gamma_\bullet\vert},
\]
the two constants:
\[
\constant>0,
\ \ \ \ \ \ \ \ \ \ \ \ \ \ \ \ \ \ \ \ \ \
\radius>0
\]
being positive, the second one not necessarily assumed
to be close to the true radius of convergence.

By what precedes, a frame for $T^{0, 1}M$ is of the form:
\[
\aligned
\overline{\mathcal{L}}_1
&
=
\frac{\partial}{\partial\overline{z}_1}
+
\overline{{\tt A}_1^1}\,
\frac{\partial}{\partial\overline{w}_1}
+\cdots+
\overline{{\tt A}_1^c}\,
\frac{\partial}{\partial\overline{w}_c},
\\
\cdots
&
\cdots\cdots\cdots\cdots\cdots\cdots\cdots\cdots\cdots
\cdots\cdots
\\
\overline{\mathcal{L}}_n
&
=
\frac{\partial}{\partial\overline{z}_n}
+
\overline{{\tt A}_n^1}\,
\frac{\partial}{\partial\overline{w}_1}
+\cdots+
\overline{{\tt A}_n^c}\,
\frac{\partial}{\partial\overline{w}_c}.
\endaligned
\]

Now, take a function:
\[
F\big(z_1,\dots,z_n,w_1,\dots,w_c\big)
\]
which is holomorphic in some neighborhood of the origin,
namely:
\[
F(z_\bullet,w_\bullet)
=
\sum_{\alpha_\bullet\in\N^n}\,
\sum_{\beta_\bullet\in\N^c}\,
F_{\alpha_\bullet,\beta_\bullet}\,
\big(z_\bullet\big)^{\alpha_\bullet}\,
\big(w_\bullet\big)^{\beta_\bullet}
\]
with Cauchy-type control:
\[
\big\vert
F_{\alpha_\bullet,\beta_\bullet}
\big\vert
\,\leqslant\,
\constant\,
\bigg(
\frac{1}{\radius}
\bigg)^{\vert\alpha_\bullet\vert+\vert\beta_\bullet\vert}.
\]

Then because the $\overline{ \mathcal{ L}}_i$ are
antiholomorphic derivations, from:
\[
0
\equiv
\frac{\partial}{\partial\overline{\sf z}_{i_1}}
\big({\sf z}_{i_2}\big)
\ \ \ \ \ \ \ \ \ \ \ \ \
{\scriptstyle{(1\,\leqslant\,i_1,\,\,i_2\,\leqslant\,n)}},
\]
one easily gets:
\[
0
\equiv
\overline{\mathcal{L}}_1(F)
\equiv\cdots\equiv
\overline{\mathcal{L}}_n(F),
\]
which shows that the restriction:
\[
\aligned
f
:=
&\,
F\big\vert_M
\\
=
&\,
F\big(
z_\bullet,\,u_\bullet+
\isqrt\,\varphi_\bullet(z_\bullet,\overline{z}_\bullet,
u_\bullet)
\big)
\endaligned
\]
is CR, and of course $\mathcal{ C}^\omega$ too.

\medskip

{\em Conversely}, start with a function:
\[
f\,\in\,
\mathcal{C}_{CR}^\omega(M).
\]
Again, localize the study at some point $p \in M$,
and take $(z_\bullet, w_\bullet)$-coordinates as above.

The natural (horizontal) coordinates on $M$ being:
\[
\big(z_\bullet,\overline{z}_\bullet,u_\bullet\big),
\]
one expresses the real analyticity of $f$ as:
\[
f
=
\sum_{\alpha_\bullet\in\N^n}\,
\sum_{\beta_\bullet\in\N^n}\,
\sum_{\gamma_\bullet\in\N^c}\,
\underbrace{f_{\alpha_\bullet,\beta_\bullet,\gamma_\bullet}}_{
\in\,\C}\,
\big(z_\bullet\big)^{\alpha_\bullet}\,
\big(\overline{z}_\bullet\big)^{\beta_\bullet}\,
\big(u_\bullet\big)^{\gamma_\bullet},
\]
with:
\[
\big\vert
f_{\alpha_\bullet,\beta_\bullet,\gamma_\bullet}
\big\vert
\,\leqslant\,
\constant\,
\bigg(
\frac{1}{\radius}
\bigg)^{\vert\alpha_\bullet\vert+\vert\beta_\bullet\vert
+\vert\gamma_\bullet\vert}.
\]

To analyze the problem, assume temporarily that a holomorphic
$F$ exists:
\[
f\big(z_\bullet,\overline{z}_\bullet,u_\bullet\big)
\equiv
F\big(z_\bullet,\,
u_\bullet+\isqrt\,\varphi_\bullet
(z_\bullet,\overline{z}_\bullet,u_\bullet)\big).
\]
Complexify the real variables $u_1, \dots, u_c$, namely
introduce new complex variables:
\[
\big(\nu_1,\dots,\nu_c\big)
\]
with:
\[
\big(u_1,\dots,u_c\big)
=
\big(\Re\,\nu_1,\dots,\Re\,\nu_c\big).
\] 
Because it just concerns power series, 
the above identity transfers to complexified
variables:
\[
f\big(z_\bullet,\overline{z}_\bullet,\nu_\bullet\big)
\equiv
F\big(z_\bullet,\,
\underbrace{\nu_\bullet+\isqrt\,\varphi_\bullet
(z_\bullet,\overline{z}_\bullet,\nu_\bullet)}_{=:\,w_\bullet}\big).
\]

Dropping now the heuristic assumption, using:
\[
0
=
\varphi_\bullet(0)
=
d\varphi_\bullet(0),
\]
the analytic implicit function theorem solves the $c$ equations:
\[
w_\bullet
=
\nu_\bullet
+
\isqrt\,\varphi_\bullet
(z_\bullet,\overline{z}_\bullet,\nu_\bullet),
\]
yielding a local analytic solution:
\[
\nu_\bullet
=
\Lambda_\bullet
\big(z_\bullet,\overline{z}_\bullet,w_\bullet\big),
\]
which by definition satisfies identically:
\[
v_\bullet
\equiv
\Lambda_\bullet
\big(
z_\bullet,\overline{z}_\bullet,\,
\nu_\bullet+\isqrt\,\varphi_\bullet
(z_\bullet,\overline{z}_\bullet,\nu_\bullet)\big),
\]
and also:
\[
\boxed{\,\,
\aligned
w_j
&
\equiv
\Lambda_j\big(z_\bullet,\overline{z}_\bullet,
w_\bullet\big)
+
\isqrt\,\varphi_j
\big(
z_\bullet,\overline{z}_\bullet,
\Lambda_\bullet
(z_\bullet,\overline{z}_\bullet,w_\bullet)
\big)\,\,
\\
&
\ \ \ \ \ \ \ \ \ \ \ \ \ \ \ \ \ \ \ \ \ \ \ \ \ \ \ \
{\scriptstyle{(j\,=\,1\,\cdots\,c)}.}
\endaligned
}
\]

Now, define:
\[
\boxed{\,\,
F(z_\bullet,\overline{z}_\bullet,w_\bullet)
:=
f\big(z_\bullet,\overline{z}_\bullet,
\Lambda_\bullet
(z_\bullet,\overline{z}_\bullet,w_\bullet)\big).\,\,}
\]

This function satisfies:
\[
\aligned
F\big(z_\bullet,\overline{z}_\bullet,
u_\bullet+\isqrt\,\varphi_\bullet(z_\bullet,
\overline{z}_\bullet,u_\bullet)\big)
&
=
f\big(z_\bullet,\overline{z}_\bullet,\,
\underbrace{
\Lambda_\bullet\big(z_\bullet,\overline{z}_\bullet,\,
u_\bullet
+
\isqrt\,\varphi_\bullet(z_\bullet,\overline{z}_\bullet,
u_\bullet)\big)}_{
\equiv\,u_\bullet}
\big)
\\
&
=
f\big(z_\bullet,\overline{z}_\bullet,u_\bullet),
\endaligned
\]
that is to say:
\[
F\big\vert_M
=
f.
\]

Moreover, $F$ is visibly holomorphic with respect to:
\[
w_\bullet
=
(w_1,\dots,w_c).
\]

\medskip\noindent{\bf Claim.}
{\em Because $f$ was assumed to 
be CR, this function $F$ is in fact also holomorphic with respect
to {\em all} the remaining variables:}
\[
z_\bullet
=
(z_1,\dots,z_n).
\]

\proof
For fixed $i = 1, \dots, n$, one would like to see that:
\[
\aligned
0
&
\overset{?}{\equiv}
\frac{\partial F}{\partial\overline{z}_i}
\\
&
=
\frac{\partial f}{\partial\overline{z}_i}
+
\sum_{j=1}^c\,
\frac{\partial\Lambda_j}{\partial\overline{z}_i}\,
\frac{\partial f}{\partial u_j}\,\,?
\endaligned
\]
But classically, fixing $i$, the partial derivatives:
\[
\frac{\partial\Lambda_j}{\partial\overline{z}_i}
\]
can be computed by coming back to the implicit equation
(boxed above) that $\Lambda_\bullet$ solves, and
by differentiating it with respect to 
$\overline{ z}_i$: 
\[
\aligned
0
&
\,\equiv\,
\frac{\partial\Lambda_j}{\partial\overline{z}_i}
+
\isqrt\,\varphi_{j,\overline{z}_i}
+
\isqrt\,
\sum_{l=1}^c\,
\varphi_{j,u_l}\,
\frac{\partial\Lambda_l}{\partial\overline{z}_i}
\\
&
\ \ \ \ \ \ \ \ \ \ \ \ \ \ \ \ \ \ \ \ \
{\scriptstyle{(j\,=\,1\,\cdots\,c)}}.
\endaligned
\]

But here, one recognizes, up to an overall 
multiplication by $\isqrt$,
the linear system seen previously:
\[
0
=
-\,\isqrt\,
\underline{{\textstyle{\frac{1}{2}}}\,
\overline{\tt A}_i^j}
+
\varphi_{j,\overline{z}_i}
+
\sum_{l=1}^c\,
\varphi_{j,u_l}\,
\underline{{\textstyle{\frac{1}{2}}}\,
\overline{\tt A}_i^l}
\]
that the coefficients $\overline{\tt A}_i^1, \dots, \overline{\tt A}_i^c$
of the:
\[
\overline{\mathcal{L}}_i
=
\frac{\partial}{\partial\overline{z}_i}
+
\sum_{j=1}^c\,
\overline{\tt A}_i^j\,
\frac{\partial}{\partial\overline{w}_j}
\]
satisfied, whence 
because the solutions to this linear system were unique,
it must be that:
\[
\boxed{\,
\frac{\partial\Lambda_j}{\partial\overline{z}_i}
=
{\textstyle{\frac{1}{2}}}\,
\overline{\tt A}_i^j.
\,}
\]
Lastly:
\[
\aligned
\frac{\partial F}{\partial\overline{z}_i}
&
=
\frac{\partial f}{\partial\overline{z}_i}
+
\sum_{j=1}^c\,
\overline{\tt A}_i^j\,
\frac{1}{2}\,
\frac{\partial f}{\partial u_j}
\\
&
=
\overline{\mathcal{L}}_i(f)
\\
&
\equiv
0,
\endaligned
\]
since $f$ was assumed to be CR.
\endproof

Thus, $F$ is an extension of $f$ that is holomorphic, which
concludes.
\endproof

\medskip\noindent{\bf Action of local biholomorphisms.}
Again, consider:
\[
M^{2n+c}
\subset
\C^{n+c}
\]
a connected $\mathcal{ C}^\kappa$ ($\kappa \geqslant 1$), or
$\mathcal{ C}^\infty$, or $\mathcal{ C}^\omega$ submanifold
which is CR generic: 
\[
TM
+
J(TM)
=
T\C^{n+c}
\big\vert_M,
\]
with:
\[
\aligned
n
&
=
\CRdim\,M,
\\
c
&
=
\codim\,M.
\endaligned
\]

Pick a point $p \in M$ and take a small open polydisc neighborhood:
\[
p\in{\sf U}_p
\subset
\C^{n+c}.
\]
Consider a local biholomorphism:
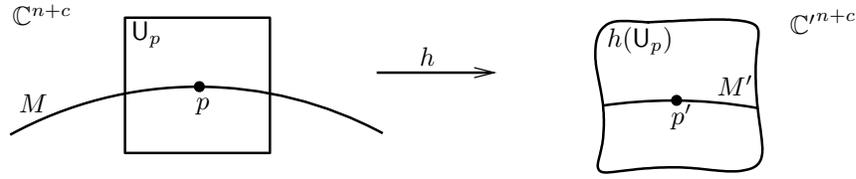
\[
h\colon\ \ \
{\sf U}_p
\overset{\sim}{\,\longrightarrow\,}
h({\sf U}_p)
\subset
{\C'}^{n+c}.
\]

\begin{center}
\begin{picture}(0,0)%
\includegraphics{M-h-polydiscs.pstex}%
\end{picture}%
\setlength{\unitlength}{4144sp}%
\begingroup\makeatletter\ifx\SetFigFont\undefined%
\gdef\SetFigFont#1#2#3#4#5{%
  \reset@font\fontsize{#1}{#2pt}%
  \fontfamily{#3}\fontseries{#4}\fontshape{#5}%
  \selectfont}%
\fi\endgroup%
\begin{picture}(4693,1052)(874,-2319)
\put(2002,-1922){\makebox(0,0)[lb]{\smash{{\SetFigFont{10}{12.0}{\familydefault}{\mddefault}{\updefault}{\color[rgb]{0,0,0}$p$}%
}}}}
\put(4844,-2009){\makebox(0,0)[lb]{\smash{{\SetFigFont{10}{12.0}{\familydefault}{\mddefault}{\updefault}{\color[rgb]{0,0,0}$p'$}%
}}}}
\put(1612,-1507){\makebox(0,0)[lb]{\smash{{\SetFigFont{10}{12.0}{\familydefault}{\mddefault}{\updefault}{\color[rgb]{0,0,0}${\sf U}_p$}%
}}}}
\put(943,-1939){\makebox(0,0)[lb]{\smash{{\SetFigFont{10}{12.0}{\familydefault}{\mddefault}{\updefault}{\color[rgb]{0,0,0}$M$}%
}}}}
\put(4459,-1546){\makebox(0,0)[lb]{\smash{{\SetFigFont{10}{12.0}{\familydefault}{\mddefault}{\updefault}{\color[rgb]{0,0,0}$h({\sf U}_p)$}%
}}}}
\put(3339,-1656){\makebox(0,0)[lb]{\smash{{\SetFigFont{10}{12.0}{\familydefault}{\mddefault}{\updefault}{\color[rgb]{0,0,0}$h$}%
}}}}
\put(5117,-1833){\makebox(0,0)[lb]{\smash{{\SetFigFont{10}{12.0}{\familydefault}{\mddefault}{\updefault}{\color[rgb]{0,0,0}$M'$}%
}}}}
\put(902,-1414){\makebox(0,0)[lb]{\smash{{\SetFigFont{10}{12.0}{\familydefault}{\mddefault}{\updefault}{\color[rgb]{0,0,0}$\C^{n+c}$}%
}}}}
\put(5552,-1470){\makebox(0,0)[lb]{\smash{{\SetFigFont{10}{12.0}{\familydefault}{\mddefault}{\updefault}{\color[rgb]{0,0,0}${\C'}^{n+c}$}%
}}}}
\end{picture}%

\end{center}

Set:
\[
p'
:=
h(p)
\]
and:
\[
M'
:=
h\big(M\cap{\sf U}_p\big).
\]

\medskip\noindent{\bf Lemma.}
{\em Then $M'$ is a $\mathcal{ C}^\kappa$
$(\kappa \geqslant 1)$, or $\mathcal{ C}^\infty$,
or $\mathcal{ C}^\omega$ CR-generic submanifold of 
the target ${\C'}^{n+c}$ having the same CR dimension and
the same codimension:}
\[
\aligned
n
&
=
\CRdim\,M'
=
\CRdim\,M
\\
c
&
=
\codim\,M'
=
\codim\,M,
\endaligned
\]
{\em and moreover, for every $q \in M \cap {\sf U}_p$:}
\[
h_*\big(T_q^cM\big)
=
T_{h(q)}^cM'.
\]

\proof
Genericity of $M \cap {\sf U}_p$ is:
\[
T_qM+J(T_qM)
=
T_q\C^{n+c}
\ \ \ \ \ \ \ \ \ \ \ \ \
{\scriptstyle{(q\,\in\,{\sf U}_p)}}.
\]
Since the biholomorphism $h$ is in particular a diffeomorphism:
\[
\aligned
h_*(T_q\C^{n+c})
&
=
T_{h(q)}{\C'}^{n+c},
\\
h_*(T_qM)
&
=
T_{h(q)}M',
\endaligned
\]
again for all $q \in {\sf U}_p$.

Also:
\[
\aligned
h_*\big(J(T_qM)\big)
&
=
J'\big(h_*(T_qM)\big)
\\
&
=
J'\big(T_{h(q)}M'\big),
\endaligned
\]
whence:
\[
\aligned
T_{h(q)}{\C'}^{n+c}
&
=
h_*\big(T_q\C^{n+c}\big)
\\
&
=
h_*\big(T_qM+J(T_qM)\big)
\\
&
=
h_*(T_qM)
+
h_*\big(J(T_qM)\big)
\\
\explain{Apply $h_*\circ J = J'\circ h_*$}
\ \ \ \ \ \ \ \ \ \ \ \ \ \ \ \ \ \ \ \ \ \ \ \ \ \ \ \ \ \ \ \ \ \ \ 
\ \ \ \
&
=
T_{h(q)}M'
+
J'\big(h_*(T_qM)\big)
\\
&
=
T_{h(q)}M'
+
J'\big(T_{h(q)}M'\big),
\endaligned
\]
namely the image $M' = h (M)$ is also
CR-generic at every $q \in {\sf U}_p$.

Lastly:
\[
\aligned
h_*\big(T_q^cM\big)
&
=
h_*\big(T_qM\cap J(T_qM)\big)
\\
&
=
h_*(T_qM)
\cap
J'\big(h_*(T_qM)\big)
\\
&
=
T_{h(q)}M
\cap
J'\big(T_{h(q)}M'\big)
\\
&
=
T_{h(q)}^cM',
\endaligned
\]
which concludes.
\endproof

\medskip\noindent{\bf Lemma.}
{\em Under the same assumptions, one has in addition:}
\[
\aligned
h_*\big(T_q^{1,0}M\big)
&
=
T_{h(q)}^{1,0}M',
\\
h_*\big(T_q^{0,1}M\big)
&
=
T_{h(q)}^{0,1}M'.
\endaligned
\]

\proof
Compute:
\[
\aligned
h_*\big(T_q^{1,0}M\big)
&
=
h_*
\Big(
T_q^{1,0}\C^{n+c}
\cap
\big[\C\otimes_\R T_qM\big]
\Big)
\\
&
=
h_*\big(T_q^{1,0}\C^{n+c}\big)
\cap
h_*\big(\C\otimes T_qM\big)
\\
&
=
T_{h(q)}^{1,0}{\C'}^{n+c}
\cap
\big[\C\otimes T_{h(q)}M\big]
\\
&
=
T_{h(q)}^{1,0}M'.
\endaligned
\]
For $T^{0, 1}$, proceed similarly or else, conjugate this.
\endproof

\noindent{\bf Exercise.}
Using the preceding lemma, get the same conclusion from:
\[
T_q^{1,0}M
=
\big\{
X_q
-
\isqrt\,J(X_q)
\colon\,
X_q
\in
T_q^cM
\big\}.
\qed
\]

\medskip\noindent{\bf Lemma.}
{\em Given two local vector field sections:}
\[
\mathcal{P}
\ \ \ \ \ \ \ \ \ \ \ \ \
\text{\rm and}
\ \ \ \ \ \ \ \ \ \ \ \ \
\mathcal{Q}
\]
{\em of the complexified tangent bundle:}
\[
\C\otimes_\R TM,
\]
{\em and given as above a local biholomorphism:}
\[
h\colon\ \ \
M
\,\longrightarrow\,
M',
\]
{\em one has:}
\[
h_*\big([\mathcal{P},\mathcal{Q}]\big)
=
\big[h_*(\mathcal{P}),\,h_*(\mathcal{Q})\big].
\]

\proof
On restriction to $M$, this follows from the general
fact that Lie brackets pass through any push forward 
by any 
$\mathcal{C}^\kappa$ ($\kappa \geqslant 1$), or $\mathcal{ C}^\infty$, 
or $\mathcal{ C}^\omega$
map between manifolds.
\endproof

\noindent{\bf Lemma.}
{\em Given a $\mathcal{ C}^\kappa$
$(\kappa \geqslant 1)$, or $\mathcal{ C}^\infty$,
or $\mathcal{ C}^\omega$ CR-generic submanifold:}
\[
\aligned
&
M^{2n+c}
\,\subset\,
\C^{n+c}
\\
&
{\scriptstyle{(c\,=\,{\sf codim}\,M,\,\,\,
n\,=\,{\sf CRdim}\,M)}},
\endaligned
\]
{\em given a local biholomorphism:}
\[
h\colon\ \ \
{\sf U}_p
\overset{\sim}{\,\longrightarrow\,}
h({\sf U}_p)
=
{\sf U}_{p'}'
\subset
{\C'}^{n+c},
\]
{\em with $p \in M$, $p' = h (p)$, setting:}
\[
M'
:=
h(M)
\,\subset\,
{\C'}^{n+c}
\ \ \ \ \ \ \ \ \ \ \ \ \
{\scriptstyle{(c\,=\,{\sf codim}\,M',\,\,\,
n\,=\,{\sf CRdim}\,M')}},
\]
{\em then for any two local frames:}
\[
\aligned
&
\big\{\mathcal{L}_1,\dots,\mathcal{L}_n\big\},
\\
&
\big\{\mathcal{L}_1',\dots,\mathcal{L}_n'\big\},
\endaligned
\]
{\em for $T^{1, 0}M$ and for $T^{ 1, 0} M'$, there exist
uniquely defined $\mathcal{ C}^{ \kappa - 1}$, 
or $\mathcal{ C}^\infty$, or $\mathcal{ C}^\omega$ local
coefficient-functions:}
\[
a_{i_1i_2}'\colon\ \ \
M'
\,\longrightarrow\,
\C
\ \ \ \ \ \ \ \ \ \ \ \ \
{\scriptstyle{(1\,\leqslant\,i_1,\,i_2\,\leqslant\,n)}},
\]
satisfying:
\[
\aligned
h_*\big(\mathcal{L}_1\big)
&
=
a_{11}'\,\mathcal{L}_1'
+\cdots+
a_{n1}'\,\mathcal{L}_n',
\\
\cdots\cdots\cdot
&
\cdots\cdots\cdots\cdots\cdots\cdots\cdots\cdots
\\
h_*\big(\mathcal{L}_n\big)
&
=
a_{1n}'\,\mathcal{L}_1'
+\cdots+
a_{nn}'\,\mathcal{L}_n'.
\endaligned
\]

\proof
This follows from the fact that both $T^{1, 0} M$ are
complex vector bundles of the same rank $n$, 
and from the already seen fact that 
$h_* ( T^{1, 0} M ) = T^{0, 1} M'$.
\endproof

\noindent{\bf Differential forms on $\C^\NN = \R^{ 2\NN}$.}
Now, on some open subset:
\[
{\sf U}
\subset
\C^\NN
=
\R^{2\NN},
\]
in the standard coordinates:
\[
\big({\sf x}_1+\isqrt\,{\sf y}_1,\,\dots\dots,\,
{\sf x}_\NN+\isqrt\,{\sf y}_\NN\big),
\]
a natural coframe for the cotangent bundle: 
\[
T^*\R^{2\NN}
\]
consists of the $2\NNN$ differential $1$-forms:
\[
d{\sf x}_1,\,\,d{\sf y}_1,
\,\,\dots\dots,\,\,
d{\sf x}_\NN,\,\,d{\sf y}_\NN,
\]
in the sense that every (local) differential $1$-form writes:
\[
\sum_{k=1}^\NN\,\big(
a_k({\sf x}_\bullet,{\sf y}_\bullet)\,
d{\sf x}_k
+
b_k({\sf x}_\bullet,{\sf y}_\bullet)\,
d{\sf y}_k
\big),
\]
with $2\NNN$ real-valued coefficient-functions:
\[
a_k,\,\,b_k\colon\ \ \
{\sf U}
\,\longrightarrow\,\R.
\]

Similarly, a (local) complex-valued differential $1$-form,
namely a section of:
\[
\C\otimes_\R T^*\R^{2\NN}
\]
over ${\sf U}$ writes:
\[
\sum_{k=1}^\NN\,
\big(
\alpha_k({\sf x}_\bullet,{\sf y}_\bullet)\,d{\sf x}_k
+
\beta_k({\sf x}_\bullet,{\sf y}_\bullet)\,d{\sf y}_k
\big),
\]
with $2\NNN$ complex-valued coefficient-functions:
\[
\alpha_k,\,\,\beta_k\colon\ \ \
{\sf U}
\,\longrightarrow\,\C.
\]

The differentials of:
\[
\aligned
{\sf z}_k
&
=
{\sf x}_k+\isqrt\,{\sf y}_k,
\\
\overline{\sf z}_k
&
=
{\sf x}_k-\isqrt\,{\sf y}_k
\endaligned
\]
being:
\[
\aligned
d{\sf z}_k
&
=
d{\sf x}_k+\isqrt\,d{\sf y}_k,
\\
d\overline{\sf z}_k
&
=
d{\sf x}_k-\isqrt\,d{\sf y}_k,
\endaligned
\]
if one solves:
\[
\aligned
d{\sf x}_k
&
=
{\textstyle{\frac{1}{2}}}
\big(
d{\sf z}_k+d\overline{\sf z}_k
\big),
\\
d{\sf y}_k
&
=
{\textstyle{\frac{\isqrt}{2}}}
\big(
-\,d{\sf z}_k+d\overline{\sf z}_k
\big),
\endaligned
\]
one obtains equivalently that sections of $\C \otimes_\R T^* \R^{2\NN}$
also write:
\[
\sum_{k=1}^\NN\,
\big(
\widetilde{\alpha}_k({\sf x}_\bullet,{\sf y}_\bullet)\,
d{\sf z}_k
+
\widetilde{\beta}_k({\sf x}_\bullet,{\sf y}_\bullet)\,
d\overline{\sf z}_k
\big),
\]
with $2\NNN$ complex-valued coefficient-functions:
\[
\widetilde{\alpha}_k,\,\,\widetilde{\beta}_k\colon\ \ \
{\sf U}
\,\longrightarrow\,\C
\]
related to the previous ones by:
\[
\aligned
\widetilde{\alpha}_k
&
=
{\textstyle{\frac{1}{2}}}\,
\alpha_k
-
{\textstyle{\frac{\isqrt}{2}}}\,
\beta_k,
\\
\widetilde{\beta}_k
&
=
{\textstyle{\frac{1}{2}}}\,
\alpha_k
+
{\textstyle{\frac{\isqrt}{2}}}\,
\beta_k.
\endaligned
\]

\medskip

Now, consider a $\mathcal{ C}^\kappa$
($\kappa \geqslant 1$), or $\mathcal{ C}^\infty$, or
$\mathcal{ C}^\omega$ diffeomorphism:
\[
(f,g)\colon\ \ \ 
{\sf U}
\overset{\sim}{\,\longrightarrow\,}
{\sf U}'
=
(f,g)({\sf U}),
\]
written as:
\[
\aligned
\big({\sf x}_1,{\sf y}_1,\dots,{\sf x}_\NN,{\sf y}_\NN\big)
&
\,\longmapsto\,
\big(
f_1({\sf x}_\bullet,{\sf y}_\bullet),
g_1({\sf x}_\bullet,{\sf y}_\bullet),
\,\dots\dots,\,
f_\NN({\sf x}_\bullet,{\sf y}_\bullet),
g_\NN({\sf x}_\bullet,{\sf y}_\bullet)
\big)
\\
&\ \ \ \,
=:
\big({\sf x}_1',{\sf y}_1',\dots,{\sf x}_\NN',{\sf y}_\NN'\big),
\endaligned
\]
the target coordinates having a prime, with target coframe:
\[
d{\sf x}_1',d{\sf y}_1',
\,\dots\dots,\,
d{\sf x}_\NN',d{\sf y}_\NN'.
\]

\medskip\noindent{\bf Pulback of differential $1$-forms.}
{\em Under the diffeomorphism $(f, g)$:}
\[
\boxed{
\aligned
(f,g)^*\big(d{\sf x}_l'\big)
&
=
\sum_{k=1}^\NN\,
\Big(
f_{l,{\sf x}_k}({\sf x}_\bullet,{\sf y}_\bullet)\,
d{\sf x}_k
+
f_{l,{\sf y}_k}({\sf x}_\bullet,{\sf y}_\bullet)\,
d{\sf y}_k
\Big),
\\
(f,g)^*\big(d{\sf y}_l'\big)
&
=
\sum_{k=1}^\NN\,
\Big(
g_{l,{\sf x}_k}({\sf x}_\bullet,{\sf y}_\bullet)\,
d{\sf x}_k
+
g_{l,{\sf y}_k}({\sf x}_\bullet,{\sf y}_\bullet)\,
d{\sf y}_k
\Big).
\endaligned
}
\]

\medskip

When dealing with differential $1$-forms having
complex coefficients, one still uses the symbol:
\[
(f,g)^*(\cdot).
\]

When such a diffeomorphism comes from a {\em biholomorphism:}
\[
\aligned
h\colon\ \ \
{\sf U}
&
\overset{\sim}{\,\longrightarrow\,}
{\sf U}'=h({\sf U})
\\
({\sf z}_1,\dots,{\sf z}_\NN)
&
\,\longmapsto\,
\big(h_1({\sf z}_\bullet),\dots,h_\NN({\sf z}_\bullet)\big)
\\
&\ \ \ \,
=:
({\sf z}_1',\dots,{\sf z}_\NN'),
\endaligned
\]
one has:
\[
\boxed{\,
\aligned
h^*\big(d{\sf z}_l'\big)
&
=
\sum_{k=1}^\NN\,
h_{l,{\sf z}_k}({\sf z}_\bullet)\,
d{\sf z}_k,
\\
h^*\big(d\overline{\sf z}_l'\big)
&
=
\sum_{k=1}^\NN\,
\overline{h_{l,{\sf z}_k}}(\overline{\sf z}_\bullet)\,
d\overline{\sf z}_k,
\endaligned
}
\]
and this motivates the introduction of two specific
subbundles of $\C \otimes_\R T^*\R^{ 2\NN}$.

\medskip\noindent{\bf Definition.}
The subbundle:
\[
T^{*(1,0)}\R^{2\NN}
\,\subset\,
\C\otimes_\R T^*\R^{2\NN}
\]
is defined in terms of its local sections of the form:
\[
\sum_{k=1}^\NN\,
\alpha_k({\sf x}_\bullet,{\sf y}_\bullet)\,
d{\sf z}_k,
\]
with complex-valued coefficient-functions $\alpha_k$.
Similarly:
\[
T^{*(0,1)}\R^{2\NN}
\,\subset\,
\C\otimes_\R T^*\R^{2\NN}
\]
has local sections of the form:
\[
\sum_{k=1}^\NN\,
\beta_k({\sf x}_\bullet,{\sf y}_\bullet)\,
d\overline{\sf z}_k.
\]

\medskip

One easily checks:
\[
\C\otimes_\R T^*\R^{2\NN}
=
T^{*(1,0)}\R^{2\NN}
\oplus
T^{*(0,1)}\R^{2\NN}.
\]

\medskip
Also, by what precedes, these bundles are invariant
through local biholomorphisms:
\[
\aligned
h^*\big(T^{*(1,0)}{\C'}^\NN\big)
&
=
T^{*(1,0)}\C^\NN,
\\
h^*\big(T^{*(0,1)}{\C'}^\NN\big)
&
=
T^{*(0,1)}\C^\NN,
\endaligned
\]
since for instance, a local section of $T^{*(1,0)} {\C'}^\NN$:
\[
\sum_{k=1}^\NN\,\alpha_k'\,d{\sf z}_k'
\]
is transformed to:
\[
\aligned
h^*\bigg(
\sum_{k=1}^\NN\,\alpha_k'\,d{\sf z}_k'
\bigg)
&
=
\sum_{k=1}^\NN\,
\alpha_k'\circ h^{-1}
\cdot
h^*\big(d{\sf z}_k'\big)
\\
&
=
\sum_{k=1}^\NN\,\sum_{l=1}^\NN\,
\alpha_k'\circ h^{-1}\,
h_{k,{\sf z}_l}
\cdot
d{\sf z}_l,
\endaligned
\]
which is still under the form of a general section of
$T^{*(1, 0)}\C^\NN$.

\medskip

Next, given any $\mathcal{ C}^\kappa$
($\kappa \geqslant 1$), or $\mathcal{ C}^\infty$, 
or $\mathcal{ C}^\omega$ (local) function:
\[
f
=
f({\sf x}_\bullet,{\sf y}_\bullet),
\]
one defines its {\sl differential:}
\[
df
:=
\sum_{k=1}^\NN\,
\big(
f_{{\sf x}_k}\,d{\sf x}_k
+
f_{{\sf y}_k}\,d{\sf y}_k
\big),
\]
together with:
\[
\aligned
\partial f
&
:=
\sum_{k=1}^\NN\,
f_{{\sf z}_k}\,d{\sf z}_k,
\\
\overline{\partial}f
&
:=
\sum_{k=1}^\NN\,
f_{\overline{\sf z}_k}\,d\overline{\sf z}_k,
\endaligned
\]
so that one checks:
\[
df
=
\partial f
+
\overline{\partial}f.
\]

\medskip\noindent{\bf Pairings.} The action of
basic $1$-forms on basic vector fields on $\R^{ 2 \NN}$ is:
\[
\aligned
d{\sf x}_{k_1}
\bigg(
\frac{\partial}{\partial{\sf x}_{k_2}}
\bigg)
&
=
\delta_{k_1,k_2},
\\
d{\sf x}_k
\bigg(
\frac{\partial}{\partial{\sf y}_l}
\bigg)
&
=
0,
\\
d{\sf y}_l
\bigg(
\frac{\partial}{\partial{\sf x_k}}
\bigg)
&
=
0,
\\
d{\sf y}_{l_1}
\bigg(
\frac{\partial}{\partial{\sf y}_{l_2}}
\bigg)
&
=
\delta_{l_1,l_2}.
\endaligned
\]
Also:
\[
\aligned
d{\sf z}_{k_1}
\bigg(
\frac{\partial}{\partial{\sf z}_{k_2}}
\bigg)
&
=
\delta_{k_1,k_2},
\\
d{\sf z}_k
\bigg(
\frac{\partial}{\partial\overline{\sf z}_l}
\bigg)
&
=
0,
\\
d\overline{\sf z}_l
\bigg(
\frac{\partial}{\partial{\sf z}_k}
\bigg)
&
=
0,
\\
d\overline{\sf z}_{l_1}
\bigg(
\frac{\partial}{\partial\overline{\sf z}_{l_2}}
\bigg)
&
=
\delta_{l_1,l_2},
\endaligned
\]
everything being coherent, for instance because:
\[
\aligned
d{\sf z}_{k_1}
\bigg(
\frac{\partial}{\partial{\sf z}_{k_2}}
\bigg)
&
=
\big(
d{\sf x}_{k_1}
+
\isqrt\,
d{\sf y}_{k_1}
\big)
\bigg(
\frac{1}{2}\,\frac{\partial}{\partial{\sf x}_{k_2}}
-
\frac{\isqrt}{2}\,\frac{\partial}{\partial{\sf y}_{k_2}}
\bigg)
\\
&
=
1\cdot\frac{1}{2}\,
\delta_{k_1,k_2}
-
\isqrt\,\frac{\isqrt}{2}\,
\delta_{k_1,k_2}
\\
&
=
\delta_{k_1,k_2}.
\endaligned
\]

Generally, a differential $1$-form:
\[
\omega
=
\sum_{k=1}^\NN\,
\Big(
\alpha_k({\sf x}_\bullet,{\sf y}_\bullet)\,
d{\sf x}_k
+
\beta_k({\sf x}_\bullet,{\sf y}_\bullet)\,
d{\sf y}_k
\Big)
\]
acts on a general vector field:
\[
L
=
\sum_{k=1}^\NN\,
\Big(
\gamma_k({\sf x}_\bullet,{\sf y}_\bullet)\,
\frac{\partial}{\partial{\sf x}_k}
+
\delta_k({\sf x}_\bullet,{\sf y}_\bullet)\,
\frac{\partial}{\partial{\sf y}_k}
\Big)
\]
(both having either real or complex coefficient-functions) 
to provide the function:
\[
\omega(L)
:=
\sum_{k=1}^\NN\,
\Big(
\alpha_k({\sf x}_\bullet,{\sf y}_\bullet)\,
\gamma_k({\sf x}_\bullet,{\sf y}_\bullet)
+
\beta_k({\sf x}_\bullet,{\sf y}_\bullet)\,
\delta_k({\sf x}_\bullet,{\sf y}_\bullet)
\Big).
\]

Similarly and equivalently:
\[
\bigg(
\sum_{k=1}^\NN
\Big(
\widetilde{\alpha}_k\,d{\sf z}_k
+
\widetilde{\beta}_k\,d\overline{\sf z}_k
\Big)
\bigg)
\bigg(
\sum_{k=1}^\NN\,
\Big(
\widetilde{\gamma}_k\,
\frac{\partial}{\partial{\sf z}_k}
+
\widetilde{\delta}_k\,
\frac{\partial}{\partial\overline{\sf z}_k}
\Big)
\bigg)
:=
\sum_{k=1}^\NN\,
\Big(
\widetilde{\alpha}_k\,\widetilde{\gamma}_k
+
\widetilde{\beta}_k\,\widetilde{\delta}_k
\Big).
\]

\medskip

Now, consider a CR-generic:
\[
\aligned
&
M^{2n+c}
\,\subset\,
\C^{n+c}
\\
&
{\scriptstyle{(c\,=\,{\sf codim}\,M,\,\,\,
n\,=\,{\sf CRdim}\,M)}},
\endaligned
\]
having equations:
\[
\aligned
v_1
&
=
\varphi_1(x_\bullet,y_\bullet,u_\bullet),
\\
\cdots
&
\cdots\cdots\cdots\cdots\cdots\cdot
\\
v_c
&
=
\varphi_c(x_\bullet,y_\bullet,u_\bullet),
\endaligned
\]
and intrinsic $T^{1, 0}M$-generators:
\[
\aligned
\mathcal{L}_1
&
=
\frac{\partial}{\partial z_1}
+
A_1^1\,\frac{\partial}{\partial u_1}
+\cdots\cdots+
A_1^c\,\frac{\partial}{\partial u_c},
\\
\cdots
&
\cdots\cdots\cdots\cdots\cdots\cdots\cdots\cdots\cdots\cdots\cdots
\cdot\cdot
\\
\mathcal{L}_n
&
=
\frac{\partial}{\partial z_n}
+
A_n^1\,\frac{\partial}{\partial u_1}
+\cdots\cdots+
A_n^c\,\frac{\partial}{\partial u_c}.
\endaligned
\]

Introduce the $n$
differential $1$-forms together with their $n$ conjugates:
\[
\aligned
\zeta_{01}
&
=
dz_1,
\ \ \ \ \ \ \ \ \ \ \ \ \ \ \ \ \ \ \ \ \ \ \ \ \ \ \ \
\overline{\zeta}_{01}
=
d\overline{z}_1,
\\
\cdots
&
\cdots\cdots\cdot
\\
\zeta_{0n}
&
=
dz_n,
\ \ \ \ \ \ \ \ \ \ \ \ \ \ \ \ \ \ \ \ \ \ \ \ \ \ \ \
\overline{\zeta}_{0n}
=
d\overline{z}_n,
\endaligned
\]
which satisfy:
\[
\aligned
\zeta_{0i_1}
\big(
\mathcal{L}_{i_2}
\big)
&
=
\delta_{i_1,i_2},
\\
\zeta_{0i}
\big(
\overline{\mathcal{L}}_l
\big)
&
=
0,
\\
\overline{\zeta}_{0l}
\big(
\mathcal{L}_i
\big)
&
=
0,
\\
\overline{\zeta}_{0l_1}
\big(
\overline{\mathcal{L}}_{l_2}
\big)
&
=
\delta_{l_1,l_2}.
\endaligned
\]

Also, introduce the {\sl real-valued} differential $1$-forms:
\[
\aligned
\rho_{01}
&
=
du_1
-
A_1^1\,dz_1
-\cdots-
A_n^1\,dz_n
-
\overline{A}_1^1\,
d\overline{z}_1
-\cdots-
\overline{A}_n^1\,d\overline{z}_n,
\\
\cdots\cdot
&
\cdots\cdots\cdots\cdots\cdots\cdots\cdots\cdots\cdots\cdots
\cdots\cdots\cdots\cdots\cdots\cdots\cdots\cdot\cdot
\\
\rho_{0c}
&
=
du_1
-
A_1^c\,dz_1
-\cdots-
A_n^c\,dz_n
-
\overline{A}_1^c\,
d\overline{z}_1
-\cdots-
\overline{A}_n^c\,d\overline{z}_n,
\endaligned
\]
which visibly satisfy:
\[
\aligned
0
&
=
\rho_{01}\big(\mathcal{L}_1\big)
=\cdots=
\rho_{01}\big(\mathcal{L}_n\big)
=
\rho_{01}\big(\overline{\mathcal{L}}_1\big)
=\cdots=
\rho_{01}\big(\overline{\mathcal{L}}_n\big),
\\
\cdot\cdot
&
\cdots\cdots\cdots\cdots\cdots\cdots\cdots\cdots\cdots\cdots
\cdots\cdots\cdots\cdots\cdots\cdots\cdots\cdots\cdot\cdot
\\
0
&
=
\rho_{0c}\big(\mathcal{L}_1\big)
=\cdots=
\rho_{0c}\big(\mathcal{L}_n\big)
=
\rho_{0c}\big(\overline{\mathcal{L}}_1\big)
=\cdots=
\rho_{0c}\big(\overline{\mathcal{L}}_n\big).
\endaligned
\]

These vanishings will be regularly abbreviated as:
\[
\boxed{\,\,
\big\{
0
=
\rho_{01}
=\cdots=
\rho_{0c}
\big\}
\,=\,
T^{1,0}M\oplus T^{0,1}M,\,\,
}
\]
within $\C \otimes_\R TM$, of course.

\medskip

Also, one realizes that:
\[
\boxed{\,\,
\aligned
\big\{
0
=
\rho_{01}
=\cdots=
\rho_{0c}
=
\overline{\zeta}_{01}
=\cdots=
\overline{\zeta}_{0n}
\big\}
&
\,=\,
T^{1,0}M,\,\,
\\
\big\{
0
=
\rho_{01}
=\cdots=
\rho_{0c}
=
\zeta_{01}
=\cdots=
\zeta_{0n}
\big\}
&
\,=\,
T^{0,1}M.\,\,
\endaligned
}
\]

\medskip

Now, since any local biholomorphism:
\[
h\colon\ \ \
M\,\longrightarrow\,M'
\]
satisfies:
\[
\aligned
h_*\big(T^{1,0}M\big)
&
=
T^{1,0}M',
\\
h_*\big(T^{0,1}M\big)
&
=
T^{0,1}M',
\\
h_*\big(
\underbrace{T^{1,0}M\oplus T^{0,1}M}_{
\{0=\rho_{01}=\cdots=\rho_{0c}\}}
\big)
&
=
\underbrace{T^{1,0}M'\oplus T^{0,1}M'}_{
\{0=\rho_{01}'=\cdots=\rho_{0c}'\}};
\endaligned
\] 
if one allows to denote by the same symbol:
\[
h_*
\equiv
\big(h^{-1}\big)^*
\]
the natural
pullback action of the {\sl inverse} 
$h^{-1}$
of $h$ on differential $1$-forms,
it follows that:
\[
\aligned
h_*(\rho_{01})
&
=
b_{11}'\,\rho_{01}'
+\cdots+
b_{1c}'\,\rho_{0c}',
\\
\cdots\cdots\cdot
&
\cdots\cdots\cdots\cdots\cdots\cdots\cdots\cdot\cdot
\\
h_*(\rho_{0c})
&
=
b_{c1}'\,\rho_{01}'
+\cdots+
b_{cc}'\,\rho_{0c}',
\endaligned
\]
for some local functions:
\[
b_{j_1j_2}'\colon\ \ \
M'
\,\longrightarrow\,
\C
\]
whose $c \times c$ matrix is invertible.

Similarly:
\[
\aligned
h_*(\zeta_{01})
&
=
d_{11}'\,\rho_{01}'
+\cdots+
d_{1c}'\,\rho_{0c}'
+
e_{11}'\,\zeta_{01}'
+\cdots+
e_{1n}'\,\zeta_{0n}',
\\
\cdots\cdots\cdot
&
\cdots\cdots\cdots\cdots\cdots\cdots\cdots\cdots\cdots\cdots\cdots
\cdots\cdots\cdots\cdots\cdot
\\
h_*(\zeta_{0n})
&
=
d_{n1}'\,\rho_{01}'
+\cdots+
d_{nc}'\,\rho_{0c}'
+
e_{n1}'\,\zeta_{01}'
+\cdots+
e_{nn}'\,\zeta_{0n}',
\endaligned
\]
for some local functions:
\[
d_{ij}'\colon\ \ \ 
M'
\,\longrightarrow\,
\C,
\ \ \ \ \ \ \ \ \ \ \ \ \ \ \ \ \ \ \ \ \ \ \ \
e_{i_1i_2}'\colon\ \ \ 
M'
\,\longrightarrow\,
\C,
\]
the $n \times n$ matrix of the latter being invertible.

Conjugating and reminding that the $\rho_{0\bullet}$
and the $\rho_{ 0\bullet}'$ are real, one has:
\[
\aligned
h_*(\overline{\zeta}_{01})
&
=
\overline{d}_{11}'\,\rho_{01}'
+\cdots+
\overline{d}_{1c}'\,\rho_{0c}'
+
\overline{e}_{11}'\,\overline{\zeta}_{01}'
+\cdots+
\overline{e}_{1n}'\,\overline{\zeta}_{0n}',
\\
\cdots\cdots\cdot
&
\cdots\cdots\cdots\cdots\cdots\cdots\cdots\cdots\cdots\cdots\cdots
\cdots\cdots\cdots\cdots\cdot
\\
h_*(\overline{\zeta}_{0n})
&
=
\overline{d}_{n1}'\,\rho_{01}'
+\cdots+
\overline{d}_{nc}'\,\rho_{0c}'
+
\overline{e}_{n1}'\,\overline{\zeta}_{01}'
+\cdots+
\overline{e}_{nn}'\,\overline{\zeta}_{0n}'.
\endaligned
\]

%%%%%%%%%%%%%%%%%%%%%%%%%%%%%%%%%%%%%%%%%%%%%%%%%%%%%%%%%%%%%%%%%%%%%

\bigskip

\section{\sf Application: 
\\
invariance of an archetypical nondegeneracy condition}
\label{invariance-archetypical}
\HEAD{\ref{invariance-archetypical}.~Application: 
invariance of an archetypical nondegeneracy condition}{
Jo\"el {\sc Merker} (Paris-Sud), 
Samuel {\sc Pocchiola} (Paris-Sud), 
Masoud {\sc Sabzevari} (Shahrekord)}

\medskip

Now, in CR dimension:
\[
n
=
{\bf 1},
\]
introduce local vector field generators:
\[
\aligned
\mathcal{L}
\ \ 
&
\text{\rm for}\ \
T^{1,0}M
\ \ \ \ \ \ \ \ \ \ \ \ \
{\scriptstyle{({\bf 1}\,=\,{\sf rank}_\C(T^{1,0}M))}},
\\
\mathcal{L}'
\ \ 
&
\text{\rm for}\ \
T^{1,0}M'
\ \ \ \ \ \ \ \ \ \ \ \ \
{\scriptstyle{({\bf 1}\,=\,{\sf rank}_\C(T^{1,0}M’))}}.
\endaligned
\]

\noindent{\bf Lemma.}
{\em At every point $q \in M$ near $p$, one has the equivalence:}
\[
\aligned
{\bf 3}
&
=
\rank_\C
\Big(
\mathcal{L}\big\vert_q,\,\,
\overline{\mathcal{L}}\big\vert_q,\,\,
\big[\mathcal{L},\overline{\mathcal{L}}\big]\big\vert_q
\Big)
\\
&
\ \ \ \ \
\Updownarrow
\\
{\bf 3}
&
=
\rank_\C
\Big(
\mathcal{L}'\big\vert_{h(q)},\,\,
\overline{\mathcal{L}}'\big\vert_{h(q)},\,\,
\big[\mathcal{L}',\overline{\mathcal{L}}'\big]\big\vert_{h(q)}
\Big).
\endaligned
\]

\proof
By what precedes:
\[
h_*\big(
T^{1,0}M
\big)
=
T^{1,0}M',
\]
hence there must exist a nowhere vanishing function:
\[
a'
\colon\ \ \
M'\longrightarrow\C
\]
defined near $p$ with:
\[
\aligned
h_*(\mathcal{L})
&
=
a'\,\mathcal{L}',
\\
h_*\big(\overline{\mathcal{L}}\big)
&
=
\overline{a}'\,\overline{\mathcal{L}}'.
\endaligned
\]
Consequently:
\[
\aligned
h_*\big(\big[\mathcal{L},\overline{\mathcal{L}}\big]\big)
&
=
\big[h_*(\mathcal{L}),\,h_*\big(
\overline{\mathcal{L}}\big)\big]
\\
&
=
\big[a'\mathcal{L}',\,\overline{a}'\overline{\mathcal{L}}'\big]
\\
&
=
a'\overline{a}'\,
\big[\mathcal{L}',\,\overline{\mathcal{L}}'\big]
+
a'\,\mathcal{L}'(\overline{a}')\cdot\overline{\mathcal{L}}'
-
\overline{a}'\,\overline{\mathcal{L}}'(a')\cdot\mathcal{L}'.
\endaligned
\] 
Dropping the mention of $h_*$, since the {\sl change of frame} matrix:
\[
\left(\!
\begin{array}{c}
\mathcal{L}
\\
\overline{\mathcal{L}}
\\
\big[\mathcal{L},\overline{\mathcal{L}}\big]
\end{array}
\!\right)
\,=\,
\left(\!
\begin{array}{ccc}
a' & 0 & 0
\\
0 & \overline{a}' & 0
\\
\ast & \ast & a'\overline{a}'
\end{array}
\!\right)
\left(\!
\begin{array}{c}
\mathcal{L}'
\\
\overline{\mathcal{L}}'
\\
\big[\mathcal{L}',\overline{\mathcal{L}}'\big]
\end{array}
\!\right)
\]
is visibly of rank $3$, the result follows.
\endproof

%%%%%%%%%%%%%%%%%%%%%%%%%%%%%%%%%%%%%%%%%%%%%%%%%%%%%%%%%%%%%%%%%%%%%

\medskip

\section{\sf Concept of Levi form}
\label{concept-levi-form}
\HEAD{\ref{concept-levi-form}.~Concept of Levi form}{
Jo\"el {\sc Merker} (Paris-Sud), 
Samuel {\sc Pocchiola} (Paris-Sud), 
Masoud {\sc Sabzevari} (Shahrekord)}

\medskip

Consider a CR-generic:
\[
M^{2n+c}
\,\subset\,
\C^{n+c}.
\]
Locally, $T^{1, 0} M$ has $n$ vector field generators:
\[
\mathcal{L}_1,\dots,\mathcal{L}_n,
\]
(more will be soon said about these), with conjugate:
\[
\overline{\mathcal{L}}_1,\dots,\overline{\mathcal{L}}_n,
\]
making a frame for $T^{0,1}M = \overline{ T^{1,0}M}$.

\smallskip

Given local $\mathcal{ C}^\omega$ functions:
\[
\aligned
\mu_1,\dots,\mu_n
\colon\ \ \
M
&
\,\longrightarrow\,
\C,
\\
\nu_1,\dots,\nu_n
\colon\ \ \
M
&
\,\longrightarrow\,
\C,
\endaligned
\]
and consider the Lie bracket:
\[
\Big[
\mu_1\mathcal{L}_1
+\cdots+
\mu_n\mathcal{L}_n,\ \
\overline{\nu}_1\overline{\mathcal{L}}_1
+\cdots+
\overline{\nu}_n\overline{\mathcal{L}}_n.
\Big]
\]
A direct expansion gives:
\[
\footnotesize
\aligned
\bigg[
\sum_{j=1}^n\,\mu_j\mathcal{L}_j,\,\,
\sum_{k=1}^n\,\overline{\nu}_k\overline{\mathcal{L}}_k
\bigg]
&
=
\sum_{j=1}^n\,\sum_{k=1}^n\,
\mu_j\,\overline{\nu}_k\,
\big[\mathcal{L}_j,\,\overline{\mathcal{L}}_k\big]
+
\\
&
\ \ \ \ \ \
+
\zero{\sum_{k=1}^n\,
\bigg(
\sum_{j=1}^n\,\mu_j\,\mathcal{L}_j
\big(\overline{\nu}_k\big)
\bigg)
\cdot
\overline{\mathcal{L}}_k}
-
\zero{\sum_{j=1}^n\,\bigg(
\sum_{k=1}^n\,\overline{\nu}_k\,\overline{\mathcal{L}}_k
\big(\mu_j\big)
\bigg)
\cdot
\mathcal{L}_j}
\\
&
\equiv
\sum_{j=1}^n\,\sum_{k=1}^n\,
\mu_j\overline{\nu}_k\,
\big[\mathcal{L}_j,\overline{\mathcal{L}}_k\big]
\ \ \ \ \ 
\mod\big(
\mathcal{L}_\bullet,\,
\overline{\mathcal{L}}_\bullet
\big).
\endaligned
\]
Implicitly here, only the cases $n = 1$ or $n = 2$ are in mind.

\medskip\noindent{\bf Definition of the complex Levi form.}
At any point $p \in M$, given two vectors:
\[
\aligned
\mathcal{M}_p
&
\in
T_p^{1,0}M,
\\
\mathcal{N}_p
&
\in
T_p^{1,0}M,
\endaligned
\]
which can then both be decomposed along the $(1, 0)$ frame:
\[
\aligned
\mathcal{M}_p
&
=
\mu_{1p}\,\mathcal{L}_1\big\vert_p
+\cdots+
\mu_{np}\,\mathcal{L}_n\big\vert_p
\\
\mathcal{N}_p
&
=
\nu_{1p}\,\mathcal{L}_1\big\vert_p
+\cdots+
\nu_{np}\,\mathcal{L}_n\big\vert_p,
\endaligned
\]
by means of certain {\em constants}:
\[
\aligned
\mu_{1p},\dots,\mu_{np}
&
\,\in\,\C,
\\
\nu_{1p},\dots,\nu_{np}
&
\,\in\,\C,
\endaligned
\]
the {\sl complex Levi form} is the Hermitian skew-bilinear form:
\[
\aligned
T_p^{1,0}M\times T_p^{1,0}M
&
\,\longrightarrow\,
\C\otimes_\R T_pM
\ \ \ \ \
\mod\,\,\big(T_p^{1,0}M\oplus T_p^{0,1}M\big)
\\
\big(\mathcal{M}_p,\mathcal{N}_p\big)
&
\,\longmapsto\,
\isqrt\,\big[\mathcal{M},\,\overline{\mathcal{N}}\big]
\big\vert_p
\,\,\,\ \ \
\mod\,\,
T_p^{1,0}M\oplus T_p^{0,1}M,
\endaligned
\]
for any two local $(1, 0)$ vector {\em fields}: 
\[
\aligned
&
\mathcal{M},
\\
&
\mathcal{N},
\endaligned
\]
which {\em extend} the vectors:
\[
\aligned
\mathcal{M}\big\vert_p
&
=
\mathcal{M}_p,
\\
\mathcal{N}\big\vert_p
&
=
\mathcal{N}_p.
\endaligned
\]

\medskip\noindent{\bf Assertion.}
{\em The result:}
\[
\isqrt\,
\big[\mathcal{M},\,\overline{\mathcal{N}}\big]\big\vert_p
\]
{\em does not depend upon choice of vector field extensions $\mathcal{ M}$, $\mathcal{ N}$.}

\proof
Consider two pairs of vector fields:
\[
\mathcal{M}^1,
\ \
\mathcal{M}^2
\ \ \ \ \ \ \ \ \ \ \ \ \
\text{\rm and}
\ \ \ \ \ \ \ \ \ \ \ \ \
\mathcal{N}^1,
\ \
\mathcal{N}^2
\]
that are sections of $T^{ 1, 0} M$ with:
\[
\mathcal{M}^1\big\vert_p
=
\mathcal{M}_p
=
\mathcal{M}^2\big\vert_p
\ \ \ \ \ \ \ \ \ \ \ \ \
\text{\rm and}
\ \ \ \ \ \ \ \ \ \ \ \ \
\mathcal{N}^1\big\vert_p
=
\mathcal{N}_p
=
\mathcal{N}^2\big\vert_p.
\]
The goal is to check:
\[
\big[\mathcal{M}^1,\,\overline{\mathcal{N}}^1\big]
\big\vert_p
\,=\,
\big[\mathcal{M}^2,\,\overline{\mathcal{N}}^2\big]
\big\vert_p.
\]

Decompose everybody along the frame:
\[
\aligned
\mathcal{M}^1
&
=
\mu_1^1\,\mathcal{L}_1
+\cdots+
\mu_n^1\,\mathcal{L}_n,
\\
\mathcal{M}^2
&
=
\mu_1^2\,\mathcal{L}_1
+\cdots+
\mu_n^2\,\mathcal{L}_n,
\\
\mathcal{N}^1
&
=
\nu_1^1\,\mathcal{L}_1
+\cdots+
\nu_n^1\,\mathcal{L}_n,
\\
\mathcal{N}^2
&
=
\nu_1^2\,\mathcal{L}_1
+\cdots+
\nu_n^2\,\mathcal{L}_n,
\endaligned
\]
with coefficient-functions having value at $p$:
\[
\aligned
\mu_1^1(p)
=
\mu_{1p}
=
\mu_1^2(p),
&
\,\,\dots\dots\dots,\,\, 
\mu_n^1(p)
=
\mu_{np}
=
\mu_n^2(p),
\\
\nu_1^1(p)
=
\nu_{1p}
=
\nu_1^2(p),
&
\,\,\dots\dots\dots,\,\, 
\nu_n^1(p)
=
\nu_{np}
=
\nu_n^2(p).
\endaligned
\]

Then the same Lie bracket expansion as above taken at $p$:
\[
\footnotesize
\aligned
\big[\mathcal{M}^1,\,\overline{\mathcal{N}}^1\big]
\big\vert_p
&
=
\bigg[
\sum_{j=1}^n\,\mu_j^1\,\mathcal{L}_j,\,\,
\sum_{k=1}^n\,\overline{\nu}_k^1\,\overline{\mathcal{L}}_k
\bigg]
\bigg\vert_p
\\
&
=
\sum_{j=1}^n\,\sum_{k=1}^n\,
\mu_j^1(p)\,\overline{\nu}_k^1(p)\,
\big[\mathcal{L}_j,\,\overline{\mathcal{L}}_k\big]
\big\vert_p
+
\\
&
\ \ \ \ \ \
+
\zero{\sum_{k=1}^n\,
\bigg(
\sum_{j=1}^n\,\mu_j^1(p)\,\mathcal{L}_j\big(
\overline{\nu}_k^1\big)(p)
\bigg)
\cdot
\overline{\mathcal{L}}_k\big\vert_p}
-
\zero{\sum_{j=1}^n\,\bigg(
\sum_{k=1}^n\,\overline{\nu}_k^1(p)\,
\overline{\mathcal{L}}_k
\big(\mu_j^1\big)\bigg)(p)
\cdot
\mathcal{L}_j\big\vert_p},
\endaligned
\]
continued as:
\[
\aligned
\big[\mathcal{M}^1,\,\overline{\mathcal{N}}^1\big]
\big\vert_p
&
\equiv
\sum_{j=1}^n\,\sum_{k=1}^n\,
\mu_j^1(p)\,\overline{\nu_k^1(p)}\,
\big[\mathcal{L}_j,\,\overline{\mathcal{L}}_k\big]
\big\vert_p
\ \ \ \ \
\mod\,\big(\mathcal{L}_\bullet,\overline{\mathcal{L}}_\bullet\big)
\\
&
\equiv
\sum_{j=1}^n\,\sum_{k=1}^n\,
\mu_{1p}\,\overline{\nu_{1p}}\,
\big[\mathcal{L}_j,\,\overline{\mathcal{L}}_k\big]
\big\vert_p
\ \ \ \ \
\mod\,\big(\mathcal{L}_\bullet,\overline{\mathcal{L}}_\bullet\big)
\\
&
\equiv
\sum_{j=1}^n\,\sum_{k=1}^n\,
\mu_j^2(p)\,\overline{\nu_k^2(p)}\,
\big[\mathcal{L}_j,\,\overline{\mathcal{L}}_k\big]
\big\vert_p
\ \ \ \ \
\mod\,\big(\mathcal{L}_\bullet,\overline{\mathcal{L}}_\bullet\big)
\\
&
\equiv
\big[\mathcal{M}^2,\,\overline{\mathcal{N}}^2\big]
\big\vert_p,
\endaligned
\]
provides the desired equality. 
\endproof

\smallskip

Assume from now that:
\[
c
=
\codim_\R\,M
=
{\bf 1},
\]
namely that:
\[
M^{2n+1}
\subset
\C^{n+1}
\]
is a hypersurface, since only the last case of:
\[
M^5
\subset
\C^3
\]
will require Levi form considerations.

\smallskip

Now, present another equivalent view of the Levi form.

\smallskip 

The quotient bundle:
\[
TM
\big/
T^cM
\]
being then of rank ${\bf 1}$, choose any $1$-form:
\[
\rho_0
:=
\text{\rm local section of}\,\,
T^*M
\]
satisfying:
\[
\big\{\rho_0=0\big\}
=
T^cM
\]
which is {\em real-valued}:
\[
\rho_0
\colon\ \ \
TM
\,\longrightarrow\,
\R.
\]
Extend it to $\C \otimes_\R TM$ by:
\[
\rho_0\big(P+\isqrt\,Q\big)
\,:=\,
\rho_0\big(P\big)
+
\isqrt\,\rho_0\big(Q\big),
\]
where $P$ and $Q$ are any two
(local) real-valued vector fields on $M$, so that:
\[
\boxed{\,
\rho_0
\Big(
T^{1,0}M
\oplus
T^{0,1}M
\Big)
=
0,\,}
\]
too.

\medskip\noindent{\bf Lemma.}
{\em One has:}
\[
\overline{\rho_0(\mathcal{M})}
=
\rho_0\big(\overline{\mathcal{M}}\big)
\]
{\em for any local section $\mathcal{ M}$ of $\C \otimes_\R TM$.}

\proof
Indeed, one decomposes:
\[
\mathcal{M}
=
P
+
\isqrt\,Q
\]
whence:
\[
\aligned
\overline{\rho_0(\mathcal{M})}
&
=
\overline{\rho_0(P)+\isqrt\,\rho_0(Q)}
\\
&
=
\rho_0(P)
-
\isqrt\,\rho_0(Q)
\\
&
=
\rho_0\big(P-\isqrt\,Q\big)
\\
&
=
\rho_0\big(\overline{\mathcal{M}}\big),
\endaligned
\]
which is so.
\endproof

For instance:
\[
\aligned
\overline{\rho_0\big[\mathcal{L}_j,\,
\overline{\mathcal{L}}_k\big]}
&
=
\rho\big(
\big[\overline{\mathcal{L}}_j,\,\mathcal{L}_k\big]
\big)
\\
&
=
-\,\rho_0\big(
\big[\mathcal{L}_k,\,\overline{\mathcal{L}}_k\big]
\big),
\endaligned
\]
hence to counterbalance this change of sign,
it will be natural to put below an $\isqrt$ factor.

\medskip\noindent{\bf Definition.}
On a hypersurface:
\[
M
\subset
\C^{n+1},
\]
in terms of any (local) frame for $T^{1, 0} M$:
\[
\big\{
\mathcal{L}_1,\dots,\mathcal{L}_n
\big\},
\]
and in terms of any $1$-form (local) section of $T^*M$:
\[
\rho_0
\colon\ \ \
TM
\,\longrightarrow\,
\R
\]
satisfying:
\[
TM\cap J(T M)
=
\big\{\rho_0=0\big\},
\]
the {\sl Levi form} of $M$ at various points $q \in M$ is: 
\[
\big((\mu_{1q},\dots,\mu_{nq}),\,\,
(\nu_{1q},\dots,\nu_{nq})\big)
\,\longmapsto\,
\isqrt\,
\sum_{j=1}^n\,\sum_{k=1}^n\,
\mu_{jq}\,\overline{\nu}_{kq}\,
\rho_0\Big(
\isqrt\,
\big[\mathcal{L}_j,\,\overline{\mathcal{L}}_k\big]
\Big)
(q).
\]

\medskip

More intuitively:
\[
\left(\!
\begin{array}{ccc}
\rho_0\big(\isqrt\,\big[\mathcal{L}_1,\overline{\mathcal{L}}_1\big]\big)
& \cdots &
\rho_0\big(\isqrt\,\big[\mathcal{L}_1,\overline{\mathcal{L}}_n\big]\big)
\\
\vdots & \ddots & \vdots
\\
\rho_0\big(\isqrt\,\big[\mathcal{L}_n,\overline{\mathcal{L}}_1\big]\big)
& \cdots &
\rho_0\big(\isqrt\,\big[\mathcal{L}_n,\overline{\mathcal{L}}_n\big]\big)
\end{array}
\!\right)
=
\text{\footnotesize\sf Hermitian matrix of the Levi form},
\]
the extra factor $\isqrt$ being present in order
to counterbalance the change of sign:
\[
\overline{\big[\mathcal{L}_j,\,\overline{\mathcal{L}}_k\big]}
=
-\,\big[\mathcal{L}_k,\,\overline{\mathcal{L}}_j\big].
\] 
Later on,
an explicit treatment of the biholomorphic 
invariance of the Levi form will be provided. 

\medskip\noindent{\bf Elementary CR-Frobenius theorem.}
{\em The Levi form of a $\mathcal{ C}^\omega$ connected CR-generic
submanifold $M^{ 2n+c} \subset \C^{n+c}$ is identically zero:}
\[
\big[T^{1,0}M,\,\overline{T^{1,0}M}\big]
\,\subset\,
\Span_\C\big(T^{1,0}M\oplus T^{0,1}M\big),
\]
{\em if and only if:}
\[
M^{2n+c}
\cong
\C^n\times\R^c.\qed
\]

\medskip

One therefore {\em excludes} such a very degenerate circumstance, 
where $M$ is usually called {\sl Levi-flat}.

\medskip

Concretely, at any point $p$ of a Levi-flat $M$, there
exist coordinates:
\[
\big(z_1,\dots,z_n,\,w_1,\dots,w_c\big)
\]
in which the graphed equations are the simplest possible:
\[
\left[
\aligned
\Im\,w_1
&
=
0,
\\
\cdots\cdot\cdot
&
\cdots\cdot\cdot
\\
\Im\,w_c
&
=
0.
\endaligned\right.
\]

So from now on, one will assume that the Levi form
is not identically zero, and because the: 
\[
\text{\sl Lie-Cartan Principle of Relocalization},
\]
is admitted, one will in fact assume that the Levi form 
is {\em nowhere} zero,
{\em i.e.}:
\[
{\bf 3}
\,\leqslant\,
\dim_\C
\Big(
T_q^{1,0}M
+
T_q^{0,1}M
+
\big[T^{1,0}M,\,T^{0,1}M\big](q),
\Big)
\]
at {\em every} point $q \in M$.

%%%%%%%%%%%%%%%%%%%%%%%%%%%%%%%%%%%%%%%%%%%%%%%%%%%%%%%%%%%%%%%%%%%%%

\bigskip

\section{\sf $\big\{\mathcal{L},\, 
\overline{ \mathcal{ L}} \big\}$-nondegeneracies 
\\
in CR dimension $n = {\bf 1}$}
\label{nondegeneracies-levi-n-1}
\HEAD{\ref{nondegeneracies-levi-n-1}.~$\big\{\mathcal{L},\, 
\overline{ \mathcal{ L}} \big\}$-nondegeneracies 
in CR dimension $n = {\bf 1}$}{
Jo\"el {\sc Merker} (Paris-Sud), 
Samuel {\sc Pocchiola} (Paris-Sud), 
Masoud {\sc Sabzevari} (Shahrekord)}

\medskip\noindent{\bf Practical consequence.}
In CR dimensions:
\[
\aligned
n
&
=
1,
\\
n
&
=
2,
\endaligned
\]
given a (local) frame for $T^{1,0} M$:
\[
\aligned
&
\big\{
\mathcal{L}_1
\big\},
\\
&
\big\{
\mathcal{L}_1,\mathcal{L}_2
\big\},
\endaligned
\]
introducing the field:
\[
\mathcal{T}
:=
\isqrt\,\big[\mathcal{L}_1,\,\overline{\mathcal{L}}_1\big],
\]
which is {\em real}:
\[
\aligned
\overline{\mathcal{T}}
=
&
\,\,\overline{\isqrt\,\big[\mathcal{L}_1,\,\overline{\mathcal{L}}_1\big]}
\\
&
=
-\,\isqrt\,
\big[\overline{\mathcal{L}}_1,\,\mathcal{L}_1\big]
\\
&
=
\mathcal{T},
\endaligned
\]
one may and will assume throughout the paper that
the Levi form is not identically zero, whence:
\[
\boxed{\,\,
\aligned
{\bf 3}
&
=
\rank_\C
\big(
\mathcal{L},\,
\overline{\mathcal{L}},\,
\mathcal{T}
\big)
\\
&
=
\rank_\C
\big(
\mathcal{L},\,
\overline{\mathcal{L}},\,
\big[\mathcal{L},\overline{\mathcal{L}}\big]
\big)
\endaligned}
\]
{\em at every point of $M$} (after allowed relocalization), this being
justified both in CR dimension $1$ and $2$ because of an:

\medskip\noindent{\bf Exercise.}
In CR dimension $2$, starting from:
\[
\left(\!
\begin{array}{cc}
\rho_0\big(\isqrt\,\big[\mathcal{L}_1,
\overline{\mathcal{L}}_1\big]\big)
&
\rho_0\big(\isqrt\,\big[\mathcal{L}_2,
\overline{\mathcal{L}}_1\big]\big)
\\
\rho_0\big(\isqrt\,\big[\mathcal{L}_1,
\overline{\mathcal{L}}_2\big]\big)
&
\rho_0\big(\isqrt\,\big[\mathcal{L}_2,
\overline{\mathcal{L}}_2\big]\big)
\end{array}
\!\right)
\]
being nonzero at a point
$p \in M$, there is a change of frame for
$T^{1, 0}M$:
\[
\left(\!
\begin{array}{c}
\mathcal{L}_1^\sharp
\\
\mathcal{L}_2^\sharp
\end{array}
\!\right)
:=
\left(\!
\begin{array}{cc}
a_{11} & a_{12}
\\
a_{21} & a_{22}
\end{array}
\!\right) 
\left(\!
\begin{array}{c}
\mathcal{L}_1
\\
\mathcal{L}_2
\end{array}
\!\right)
\]
with constant coefficients which makes:
\[
{\bf 3}
=
\rank_\C
\big(
\mathcal{L}_1^\sharp,\,
\overline{\mathcal{L}}_1^\sharp,\,
\big[\mathcal{L}_1^\sharp,\,
\overline{\mathcal{L}}_1^\sharp\big]
\big).
\qed
\]

\smallskip

In CR dimension $n = 1$, one drops index mention:
\[
\mathcal{L}
\equiv
\mathcal{L}_1.
\]

Among the four cases:
\[
\aligned
&
M^3
\,\subset\,\C^2
\colon\ \ \ 
n
=
{\bf 1},
\\
&
M^4
\,\subset\,\C^3
\colon\ \ \ 
n
=
{\bf 1},
\\
&
M^5
\,\subset\,
\aligned
&
\C^4
\colon\ \ \ 
n
=
{\bf 1},
\\
&
\C^3
\colon\ \ \ 
n
=
{\bf 2},
\endaligned
\endaligned
\]

\medskip\noindent$\bullet$\,\,
The first ends up:
\[
\text{\rm For}\ \
M^3\subset\C^2,
\ \
\text{\rm already}\ \
\big\{
\mathcal{L},\,\overline{\mathcal{L}},\,\mathcal{T}\big\}\ \
\text{\rm makes up a frame for}\ \
\C\otimes_\R TM.
\]

\medskip\noindent$\bullet$\,\,
The second is awaiting:
\[
\aligned
\text{\rm For}\ \
M^4\subset\C^3,
\ \
&
{\bf 1}\,\,
\text{\rm field is still missing in}
\big\{\mathcal{L},\,\overline{\mathcal{L}},\,\mathcal{T}\big\}\ \
\text{\rm to complete a frame}
\\
&
\text{\rm since}\ \
{\bf 4}=\rank_\C\big(\C\otimes_\R TM\big).
\endaligned
\]

\medskip\noindent$\bullet$\,\,
The third is also awaiting:
\[
\aligned
\text{\rm For}\ \
M^5\subset\C^4,
\ \
&
{\bf 2}\,\,
\text{\rm fields are still missing in}
\big\{\mathcal{L},\,\overline{\mathcal{L}},\,\mathcal{T}\big\}\ \
\text{\rm to complete a frame}
\\
&
\text{\rm since}\ \
{\bf 5}=\rank_\C\big(\C\otimes_\R TM\big).
\endaligned
\]

\medskip\noindent$\bullet$\,\,
And the fourth and last will be dealt with subsequently:
\[
\aligned
\text{\rm For}\ \
M^5\subset\C^3,\ \ 
&
\text{\rm already}\ \
\big\{\mathcal{L}_1,\,\mathcal{L}_2,\,\overline{\mathcal{L}}_1,\,
\overline{\mathcal{L}}_2,\,\mathcal{T}\big\}\ \
\text{\rm makes up a frame for}\ \
\C\otimes_\R TM,
\\
&
\text{\rm but the generic rank of the Levi matrix}\ 
\rho_0\left(\!\!
\begin{array}{cc}
\isqrt\,\big[\mathcal{L}_1,\overline{\mathcal{L}}_1\big]
&
\isqrt\,\big[\mathcal{L}_2,\overline{\mathcal{L}}_1\big]
\\
\isqrt\,\big[\mathcal{L}_1,\overline{\mathcal{L}}_2\big]
&
\isqrt\,\big[\mathcal{L}_2,\overline{\mathcal{L}}_2\big]
\end{array}
\!\!\right)
\\
&
\text{\rm may be equal to}\ \
{\bf 1}\ \
\text{\rm or}\ \
{\bf 2},\ \
\text{\rm two subcases to study independently};
\endaligned
\]

\medskip

Discuss here only CR dimension $n = {\bf 1}$, 
postponing $n = {\bf 2}$.

\medskip

At least, one arrives 
at a first {\sl general class} of CR-generic manifolds:
\[
\aligned
&
\boxed{\text{\sf General Class $\text{\sf I}$:}}
\\
&
\boxed{\,\,
\aligned
M^3\subset\C^2
\ \
&
\text{\rm with}\ \
\Big\{\mathcal{L},\,\overline{\mathcal{L}},\,\,
\big[\mathcal{L},\overline{\mathcal{L}}\big]\Big\}\ \
\\
&
\text{\rm constituting a frame for}\ \
\C\otimes_\R TM.\,\,
\endaligned}
\endaligned
\]

\medskip\noindent{\bf Zariski-Generic degeneracies of $M^4 \subset \C^3$.}
Consider therefore next a $\mathcal{ C}^\omega$ CR-generic:
\[
M^4
\subset
\C^3.
\]
Pick a local generator:
\[
\mathcal{L}
\]
of $T^{1, 0} M$ and consider:
\[
\big[\mathcal{L},\overline{\mathcal{L}}\big].
\]
By what precedes:
\[
{\bf 3}
=
\rank_\C
\Big(
\mathcal{L},\overline{\mathcal{L}},\,
\big[\mathcal{L},\overline{\mathcal{L}}\big]
\Big),
\]
at every point. But this is still less than:
\[
{\bf 4}
=
\rank_\C
\big(
\C\otimes_\R TM
\big).
\]

Accordingly, introduce next the two possible further Lie brackets 
of length $3$:
\[
\aligned
&
\big[\mathcal{L},\,[\mathcal{L},\overline{\mathcal{L}}]\big],
\\
&
\big[\overline{\mathcal{L}},\,[\mathcal{L},\overline{\mathcal{L}}]\big],
\endaligned
\]
satisfying:
\[
\overline{\big[\mathcal{L},\,
\big[\mathcal{L},\overline{\mathcal{L}}\big]\big]}
=
-\,\big[\overline{\mathcal{L}},\,
\big[\mathcal{L},\overline{\mathcal{L}}\big]\big].
\]

\medskip\noindent{\bf Degeneracy assumption.} 
{\em Inspect the exceptional supposition:}
\[
{\bf 3}
=
\rank_\C\Big(
\mathcal{L},
\overline{\mathcal{L}},\,
\big[\mathcal{L},\overline{\mathcal{L}}\big],\,\,
\big[\mathcal{L},
\big[\mathcal{L},\overline{\mathcal{L}}\big]\big],\,\,
\big[\overline{\mathcal{L}},\,
\big[\mathcal{L},\overline{\mathcal{L}}\big]\big] 
\Big),
\]
{\em assumed to hold at every point.}

\medskip

Then locally, there exist $\mathcal{ C}^\omega$ functions so that:
\[
\big[\mathcal{L},\,
\big[\mathcal{L},\overline{\mathcal{L}}\big]
\big]
=
a\cdot\mathcal{L}
+
b\cdot\overline{\mathcal{L}}
+
c\cdot\big[\mathcal{L},\overline{\mathcal{L}}\big],
\]
whence conjugating:
\[
\big[\overline{\mathcal{L}},\,
\big[\mathcal{L},\overline{\mathcal{L}}\big]\big]
=
\overline{b}\cdot\mathcal{L}
+
\overline{a}\cdot\overline{\mathcal{L}}
-
\overline{c}\cdot\big[\mathcal{L},\overline{\mathcal{L}}\big].
\]

Compute then:
\[
\aligned
\big[\mathcal{L},\big[\mathcal{L},\,
\big[\mathcal{L},\overline{\mathcal{L}}\big]\big]\big]
&
=
\mathcal{L}(a)\cdot\mathcal{L}
+
\mathcal{L}(b)\cdot\overline{\mathcal{L}}
+
\mathcal{L}(c)\cdot\big[\mathcal{L},\overline{\mathcal{L}}\big]
+
\\
&
\ \ \ \ \
+
a\cdot\zero{[\mathcal{L},\mathcal{L}]}
+
b\cdot
\big[\mathcal{L},\overline{\mathcal{L}}\big]
+\,
c\cdot
\!\!\!
\underbrace{\big[\mathcal{L},
\big[\mathcal{L},\overline{\mathcal{L}}\big]\big]}_{
a\mathcal{L}+b\overline{\mathcal{L}}
+c[\mathcal{L},\overline{\mathcal{L}}]}
\\
&
\equiv
0
\ \ \ \ \
\mod\big(\mathcal{L},\overline{\mathcal{L}},
\big[\mathcal{L},\overline{\mathcal{L}}\big]\big).
\endaligned
\]

Similarly (exercise):
\[
\aligned
\big[\overline{\mathcal{L}},[\mathcal{L},
\big[\mathcal{L},\overline{\mathcal{L}}\big]]\big]
&
\equiv
0
\ \ \ \ \
\mod\big(\mathcal{L},\overline{\mathcal{L}},
\big[\mathcal{L},\overline{\mathcal{L}}\big]\big),
\\
\big[\mathcal{L},[\overline{\mathcal{L}},
\big[\mathcal{L},\overline{\mathcal{L}}\big]]\big]
&
\equiv
0
\ \ \ \ \
\mod\big(\mathcal{L},\overline{\mathcal{L}},
\big[\mathcal{L},\overline{\mathcal{L}}\big]\big),
\\
\big[\overline{\mathcal{L}},[\overline{\mathcal{L}},
\big[\mathcal{L},\overline{\mathcal{L}}\big]]\big]
&
\equiv
0
\ \ \ \ \
\mod\big(\mathcal{L},\overline{\mathcal{L}},
\big[\mathcal{L},\overline{\mathcal{L}}\big]\big).
\endaligned
\]

All further iterated Lie brackets, {\em e.g.}:
\[
\footnotesize
\aligned
\big[\mathcal{L},\big[\mathcal{L},[\mathcal{L},
\big[\mathcal{L},\overline{\mathcal{L}}\big]]\big]\big]
&
=
\big[
\mathcal{L},\,\mod\,(\same)\big]
\\
&
=
\Big[\mathcal{L},\,
\function\cdot\mathcal{L}
+
\function\cdot\overline{\mathcal{L}}
+
\function\cdot
\big[\mathcal{L},\overline{\mathcal{L}}\big]
\Big]
\\
&
=
\mathcal{L}(\function)\cdot\mathcal{L}
+
\mathcal{L}(\function)\cdot\overline{\mathcal{L}}
+
\mathcal{L}(\function)\cdot
\big[\mathcal{L},\overline{\mathcal{L}}\big]
+
\\
&
\ \ \ \ \
+
\function\cdot
\zero{\big[\mathcal{L},\mathcal{L}\big]}
+
\function\cdot
\big[\mathcal{L},\overline{\mathcal{L}}\big]
+\,
\function\cdot
\!\!\!
\underbrace{\big[\mathcal{L},
\big[\mathcal{L},\overline{\mathcal{L}}\big]\big]}_{
a\mathcal{L}+b\overline{\mathcal{L}}
+c[\mathcal{L},\overline{\mathcal{L}}]}
\\
&
\equiv
0
\ \ \ \ \
\mod\big(\mathcal{L},\overline{\mathcal{L}},\mathcal{T}\big)
\endaligned
\]
also mod out to zero, and similarly also:
\[
\big[\overline{\mathcal{L}},\,
\mod\,(\same)\big]
\equiv
0
\ \ \ \ \
\mod\big(\mathcal{L},\overline{\mathcal{L}},
\big[\mathcal{L},\overline{\mathcal{L}}\big]\big).
\]

\medskip\noindent{\bf Definition.}
For a $\mathcal{ C}^\omega$ connected CR-generic submanifold:
\[
M^{2n+c}
\subset
\C^{n+c}
\] 
with $T^{1, 0}M$ having local generators:
\[
\mathcal{L}_1,\dots,\mathcal{L}_n,
\]
on various open subsets, set:
\[
\aligned
\LL_{\mathcal{L},\overline{\mathcal{L}}}^1
&
:=
\mathcal{C}^\omega\text{\rm -linear combinations of}\,\,
\mathcal{L}_1,\dots,\mathcal{L}_n,\,\,
\overline{\mathcal{L}}_1,\dots,\overline{\mathcal{L}}_n,
\\
\LL_{\mathcal{L},\overline{\mathcal{L}}}^2
&
:=
\mathcal{C}^\omega\text{\rm -linear combinations of vector fields}\,\,
\mathcal{M}^1\in\LL_{\mathcal{L},\overline{\mathcal{L}}}^1
\\
&
\ \ \ \ \ \ \ \ \ \ \ \ \ 
\text{\rm and of brackets}\,\,
\big[\mathcal{L}_k,\mathcal{M}^1\big],\,\,
\big[\overline{\mathcal{L}}_k,\mathcal{M}^1\big],
\\
\LL_{\mathcal{L},\overline{\mathcal{L}}}^3
&
:=
\mathcal{C}^\omega\text{\rm -linear combinations of vector fields}\,\,
\mathcal{M}^2\in\LL_{\mathcal{L},\overline{\mathcal{L}}}^2
\\
&
\ \ \ \ \ \ \ \ \ \ \ \ \ 
\text{\rm and of brackets}\,\,
\big[\mathcal{L}_k,\mathcal{M}^2\big],\,\,
\big[\overline{\mathcal{L}}_k,\mathcal{M}^2\big],
\\
\cdots\cdots
&
\cdots\cdots\cdots\cdots\cdots\cdots\cdots\cdots\cdots\cdots\cdots\cdots
\cdots\cdots\cdots\cdots\cdots\cdots
\\
\LL_{\mathcal{L},\overline{\mathcal{L}}}^{\nu+1}
&
:=
\mathcal{C}^\omega\text{\rm -linear combinations of vector fields}\,\,
\mathcal{M}^\nu\in\LL_{\mathcal{L},\overline{\mathcal{L}}}^\nu
\\
&
\ \ \ \ \ \ \ \ \ \ \ \ \ 
\text{\rm and of brackets}\,\,
\big[\mathcal{L}_k,\mathcal{M}^\nu\big],\,\,
\big[\overline{\mathcal{L}}_k,\mathcal{M}^\nu\big],
\endaligned
\]
these being defined on local open subsets of $M$ (sheaf language is skipped).

\medskip\noindent{\bf Definition.}
Set:
\[
\LL_{\mathcal{L},\overline{\mathcal{L}}}^{\sf Lie}
:=
\bigcup_{\nu\geqslant 1}\,
\LL_{\mathcal{L},\overline{\mathcal{L}}}^\nu
=
\LL_{\mathcal{L},\overline{\mathcal{L}}}^\infty.
\]

\medskip
Alternatively, working with {\em real vector fields}, one introduces:
\[
\LL_{{\sf Re}\,\mathcal{L},{\sf Im}\,\mathcal{L}}^{\sf Lie}.
\]

\medskip\noindent{\bf Known real analytic fact.}
{\em There exists an integer:}
\[
\boxed{\,\,c_M\,\,}
\]
{\em with:}
\[
0\leqslant c_M\leqslant c
\]
{\em and a proper real analytic subset:}
\[
\Sigma\,
\subsetneqq\,
M
\]
{\em such that at every point:}
\[
q\,\in\,
M\backslash\Sigma,
\]
{\em one has:}
\[
\aligned
\dim_\C\,
\Big(
\LL_{\mathcal{L},\overline{\mathcal{L}}}^{\sf Lie}(q)
\Big)
&
=
2n+c_M
\\
&
=
\constant,
\endaligned
\]
{\em or equivalently (mental exercise):}
\[
\dim_\R
\Big(
\LL_{{\sf Re}\,\mathcal{L},{\sf Im}\,\mathcal{L}}^{\sf Lie}(q)
\Big)
=
2n+c_M.\qed
\]

\medskip\noindent{\bf Known generalized CR Frobenius theorem
(\cite{ Merker-Porten-2006}).}
{\em Every point:}
\[
q\,\in\,M\backslash\Sigma
\]
{\em has a small open neighbordhood:}
\[
{\sf U}_q
\subset
\C^{n+c}
\]
{\em in which:}
\[
M^{2n+c}
\cong\,
\underline{M}^{2n+c}
\]
{\em biholomorphically, with:}
\[
\underline{M}^{2n+c}
\subset
\C^{n+c_M}
\times
\R^{c-c_M}
\]
{\em $\mathcal{ C}^\omega$ and CR-generic. Furthermore, in appropriate coordinates:}
\[
\big(z_1,\dots,z_n,w_1,\dots,w_{c_M},w_{c_M+1},\dots,w_c\big),
\]
{\em the equations of $\underline{ M}$ are:}
\[
\left[
\aligned
v_1
&
=
\varphi_1\big(x_1,\dots,x_n,y_1,\dots,y_n,\,
u_1,\dots,u_{c_M},u_{c_M+1},\dots,u_c\big),
\\
\cdots
&
\cdots\cdots\cdots\cdots\cdots\cdots\cdots\cdots\cdots
\cdots\cdots\cdots\cdots\cdots\cdots\cdots\cdot\cdot
\cdots
\\ 
v_{c_M}
&
=
\varphi_{c_M}\big(x_1,\dots,x_n,y_1,\dots,y_n,\,
u_1,\dots,u_{c_M},u_{c_M+1},\dots,u_c\big),
\\
v_{c_M+1}
&
=
0,
\\
\cdot\cdots
&
\cdots\cdots\cdots\cdots\cdots\cdots\cdots\cdots\cdots
\cdots\cdots\cdots\cdots\cdots\cdots\cdots\cdot\cdot
\\
v_c
&
=
0.
\qed
\endaligned\right.
\]

\medskip

For the case $M^4 \subset \C^3$ at hand, the above exceptional supposition:
\[
{\bf 3}
=
\rank_\C\Big(
\mathcal{L},\,
\overline{\mathcal{L}},\,
\big[\mathcal{L},\overline{\mathcal{L}}\big],\,\,
\big[\mathcal{L},
\big[\mathcal{L},\overline{\mathcal{L}}\big]\big],\,\,
\big[\overline{\mathcal{L}},
\big[\mathcal{L},\overline{\mathcal{L}}\big]\big]
\Big),
\]
therefore yields:
\[
\aligned
{\bf 3}
&
=
\rank_\C
\Big(
\LL_{\mathcal{L},\overline{\mathcal{L}}}^3
\Big)
\\
&
=
\rank_\C
\Big(
\LL_{\mathcal{L},\overline{\mathcal{L}}}^4
\Big)
=\cdots=
\\
&
=
\rank_\C
\Big(
\LL_{\mathcal{L},\overline{\mathcal{L}}}^{\sf Lie}
\Big)
\\
&
=
\rank_\R
\Big(
\LL_{{\sf Re}\,\mathcal{L},{\sf Im}\,\mathcal{L}}^{\sf Lie}
\Big)
\,<\,
\,{\bf 4}\,
=
\dim_\R\,M,
\endaligned
\]
whence:
\[
\aligned
M^4
\,\cong\,
\underline{M}^4,
\endaligned
\]
with:
\[
\underline{M}^4
\,\subset\,
\C^2\times\R
\]
having, in appropriate coordinates, equations:
\[
\aligned
v_1
&
=
\varphi_1\big(x,y,u_1,u_2\big),
\\
v_2
&
=
0.
\endaligned
\]

\medskip\noindent{\bf Interpretation.} 
{\em One sets aside such an exceptional supposition, because the
equivalence problem reduces to that of an:}
\[
M^3
\subset
\C^2
\]
{\em in smaller dimension, plus $1$ real parameter coming
from $(\cdot) \times \R$.\qed}

\medskip

Consequently, neglecting such a degenerate subclass, 
one arrives at an interesting class of CR manifolds:
\[
M^4
\subset
\C^3
\]
satisfying:
\[
{\bf 4}
=
\rank_\C
\Big(\mathcal{L},\,\overline{\mathcal{L}},\,
[\mathcal{L},\overline{\mathcal{L}}\big],\,\,
\big[\mathcal{L},\,
[\mathcal{L},\overline{\mathcal{L}}\big]\big],\,\,
\big[\overline{\mathcal{L}},\,[\mathcal{L},\overline{\mathcal{L}}\big]\big]
\Big),
\]
which was discovered by Beloshapka (\cite{ Beloshapka-1998}).

\medskip\noindent{\bf Observational lemma.}
{\em In fact, then simultaneously:}
\[
\aligned
{\bf 4}
&
=
\rank_\C
\Big(\mathcal{L},\,\overline{\mathcal{L}},\,
[\mathcal{L},\overline{\mathcal{L}}\big],\,\,
\big[\mathcal{L},\,[\mathcal{L},\overline{\mathcal{L}}\big]\big]
\Big)
\\
&
=
\rank_\C
\Big(\mathcal{L},\,\overline{\mathcal{L}},\,
[\mathcal{L},\overline{\mathcal{L}}\big],\,\,
\big[\overline{\mathcal{L}},\,[\mathcal{L},\overline{\mathcal{L}}\big]\big]
\Big).\,\,
\endaligned
\]

\proof
Indeed, near points $q \in M$ where:
\[
\aligned
{\bf 4}
&
=
\rank_\C
\Big(\mathcal{L}\big\vert_q,\,\overline{\mathcal{L}}\big\vert_q,\,
[\mathcal{L},\overline{\mathcal{L}}\big]\big\vert_q,\,\,
\big[\mathcal{L},\,
[\mathcal{L},\overline{\mathcal{L}}\big]\big]\big\vert_q
\Big)
\\
&
=
\dim_\R\,M,
\endaligned
\]
the corresponding frame for $\C \otimes_\R TM$:
\[
\Big\{
\mathcal{L},\,\overline{\mathcal{L}},\,
\big[\mathcal{L},\overline{\mathcal{L}}\big],\,\,
\big[\mathcal{L},\,[\mathcal{L},\overline{\mathcal{L}}\big]\big]
\Big\}
\]
enables one to express:
\[
\big[\overline{\mathcal{L}},\,
\big[\mathcal{L},\overline{\mathcal{L}}\big]\big]
=
a\cdot\mathcal{L}
+
b\cdot\overline{\mathcal{L}}
+
c\cdot\big[\mathcal{L},\overline{\mathcal{L}}\big]
+
d\cdot\big[\mathcal{L},
\big[\mathcal{L},\overline{\mathcal{L}}\big]\big],
\]
whence by conjugation:
\[
\aligned
-\,\big[\mathcal{L},
\big[\mathcal{L},\overline{\mathcal{L}}\big]\big]
&
=
\overline{b}\cdot\mathcal{L}
+
\overline{a}\cdot\overline{\mathcal{L}}
-
\overline{c}\cdot
\big[\mathcal{L},\overline{\mathcal{L}}\big]
-
\overline{d}\cdot
\underbrace{\big[\overline{\mathcal{L}},
\big[\mathcal{L},\overline{\mathcal{L}}\big]\big]}_{
{\sf replace}}
\\
&
\equiv
-\,d\overline{d}\cdot
\big[\mathcal{L},
\big[\mathcal{L},\overline{\mathcal{L}}\big]\big]
\ \ \ \ \
\mod\big(\mathcal{L},\overline{\mathcal{L}},
\big[\mathcal{L},\overline{\mathcal{L}}\big]\big),
\endaligned
\]
so that:
\[
d\overline{d}
\equiv
1,
\]
identically as $\mathcal{ C}^\omega$ functions defined on $M$ near $q$.
\endproof

Hence one arrives at a second {\sl general class} of CR-generic manifolds:
\[
\aligned
&
\boxed{\text{\sf General Class $\text{\sf II}$:}}
\\
&
\boxed{\,\,
\aligned
M^4\subset\C^3
\ \
&
\text{\rm with}\ \
\Big\{\mathcal{L},\,\overline{\mathcal{L}},\,\,
\big[\mathcal{L},\overline{\mathcal{L}}\big],\,\,
\big[\mathcal{L},\,\big[\mathcal{L},\overline{\mathcal{L}}\big]\big]
\Big\}\ \
\\
&
\text{\rm constituting a frame for}\ \
\C\otimes_\R TM.\,\,
\endaligned}
\endaligned
\]

\noindent{\bf Zariski-Generic degeneracies of $M^5 \subset \C^4$.}
Now, consider an:
\[
M^5
\subset 
\C^4,
\]
having:
\[
\aligned
{\bf 1}
&
=
\CRdim\,M,
\\
{\bf 3}
&
=
\codim\,M.
\endaligned
\]

Similarly to the case of $M^4 \subset \C^3$, 
inspect the exceptional supposition:
\[
{\bf 3}
=
\rank_\C\Big(
\mathcal{L},\,
\overline{\mathcal{L}},\,
\big[\mathcal{L},\overline{\mathcal{L}}\big],\,\,
\big[\mathcal{L},
\big[\mathcal{L},\overline{\mathcal{L}}\big]\big],\,\,
\big[\overline{\mathcal{L}},
\big[\mathcal{L},\overline{\mathcal{L}}\big]\big]
\Big),
\]
assumed to hold at every point.

\medskip

Again (mental exercise), this entails:
\[
M^5
\,\cong\,
\underline{M}^5,
\]
with:
\[
\underline{M}^5
\,\subset\,
\C^2\times\R^2
\]
being represented in local coordinates as:
\[
\aligned
v_1
&
=
\varphi_1\big(x,y,u_1,u_2,u_3\big),
\\
v_2
&
=
0,
\\
v_3
&
=
0.
\endaligned
\]

\medskip\noindent{\bf Same interpretation.} {\em One sets aside
such an exceptional supposition, for the
equivalence problem reduces to that of an:}
\[
M^3
\subset
\C^2
\]
{\em in smaller dimension, plus $2$ real parameters coming
from $(\cdot) \times \R^2$.\qed}

\medskip

Consequently, neglecting such a degenerate subclass, 
one may assume:
\[
{\bf 4}
\,\leqslant\,
\rank_\C
\Big(\mathcal{L},\,\overline{\mathcal{L}},\,
[\mathcal{L},\overline{\mathcal{L}}\big],\,\,
\big[\mathcal{L},\,[\mathcal{L},\overline{\mathcal{L}}\big]\big],\,\,
\big[\overline{\mathcal{L}},\,
[\mathcal{L},\overline{\mathcal{L}}\big]\big]
\Big),
\]
and since ${\bf 5}$ fields are present,
one arrives at an interesting class of CR-generic submanifolds:
\[
M^5
\subset
\C^4
\]
satisfying:
\[
{\bf 5}
=
\rank_\C
\Big(\mathcal{L},\,\overline{\mathcal{L}},\,
[\mathcal{L},\overline{\mathcal{L}}\big],\,\,
\big[\mathcal{L},\,[\mathcal{L},\overline{\mathcal{L}}\big]\big],\,\,
\big[\overline{\mathcal{L}},\,[\mathcal{L},\overline{\mathcal{L}}\big]\big]
\Big),
\]
which was also discovered by Beloshapka.

Hence one arrives at a third {\sl general class} of CR-generic manifolds:
\[
\aligned
&
\boxed{\text{\sf General Class $\text{\sf III}_{\text{\sf 1}}$:}}
\\
&
\boxed{\,\,
\aligned
M^5\subset\C^4
\ \
&
\text{\rm with}\ \
\Big\{\mathcal{L},\,\overline{\mathcal{L}},\,\,
\big[\mathcal{L},\overline{\mathcal{L}}\big],\,\,
\big[\mathcal{L},\,
\big[\mathcal{L},\overline{\mathcal{L}}\big]\big],\,\,
\big[\overline{\mathcal{L}},\,
\big[\mathcal{L},\overline{\mathcal{L}}\big]\big]
\Big\}\,\,
\\
&
\text{\rm constituting a frame for}\ \
\C\otimes_\R TM.
\endaligned}
\endaligned
\]

%%%%%%%%%%%%%%%%%%%%%%%%%%%%%%%%%%%%%%%%%%%%%%%%%%%%%%%%%%%%%%%%%%%%%

\bigskip

\section{\sf Yet a last new general class 
\\
of {\bf 5}-dimensional CR manifolds $M^5 \subset \C^4$}
\label{yet-new-class}
\HEAD{\ref{yet-new-class}.~Yet a last new general class of 
{\bf 5}-dimensional CR manifolds $M^5 \subset \C^4$}{
Jo\"el {\sc Merker} (Paris-Sud), 
Samuel {\sc Pocchiola} (Paris-Sud), 
Masoud {\sc Sabzevari} (Shahrekord)}

\medskip

But for $M^5 \subset \C^4$, from:
\[
{\bf 4}
\,\leqslant\,
\rank_\C
\Big(\mathcal{L},\,\overline{\mathcal{L}},\,
[\mathcal{L},\overline{\mathcal{L}}\big],\,\,
\big[\mathcal{L},\,[\mathcal{L},\overline{\mathcal{L}}\big]\big],\,\,
\big[\overline{\mathcal{L}},\,[\mathcal{L},\overline{\mathcal{L}}\big]\big]
\Big),
\]
it can yet very well happen that:
\[
{\bf 4}
\,=\,
\rank_\C
\Big(\mathcal{L},\,\overline{\mathcal{L}},\,
[\mathcal{L},\overline{\mathcal{L}}\big],\,\,
\big[\mathcal{L},\,[\mathcal{L},\overline{\mathcal{L}}\big]\big],\,\,
\big[\overline{\mathcal{L}},\,[\mathcal{L},\overline{\mathcal{L}}\big]\big]
\Big),
\]
of course at every point after relocalization.

\smallskip

As above, from:
\[
\big[\overline{\mathcal{L}},[\mathcal{L},\overline{\mathcal{L}}]\big]
=
a\cdot\mathcal{L}
+
b\cdot\overline{\mathcal{L}}
+
c\cdot[\mathcal{L},\overline{\mathcal{L}}]
+
d\cdot\big[\mathcal{L},[\mathcal{L},\overline{\mathcal{L}}]\big],
\]
one gets:
\[
d\overline{d}
\equiv
1,
\]
so that it is legitimate to assume at the same time:
\[
\aligned
{\bf 4}
&
\,=\,
\rank_\C
\Big(\mathcal{L},\,\overline{\mathcal{L}},\,
[\mathcal{L},\overline{\mathcal{L}}\big],\,\,
\big[\mathcal{L},\,[\mathcal{L},\overline{\mathcal{L}}\big]\big]
\Big)
\\
&
\,=\,
\rank_\C
\Big(\mathcal{L},\,\overline{\mathcal{L}},\,
[\mathcal{L},\overline{\mathcal{L}}\big],\,\,
\big[\overline{\mathcal{L}},\,
[\mathcal{L},\overline{\mathcal{L}}\big]\big]
\Big).
\endaligned
\]

Hence {\bf 1} iterated Lie bracket is still missing to generate:
\[
{\bf 5}
=
\rank_\C
\big(
\C\otimes_\R TM
\big).
\]
All next candidates, namely the length $4$ Lie brackets in:
\[
\LL_{\mathcal{L},\overline{\mathcal{L}}}^4,
\]
are four in sum:
\[
\aligned
&
\big[\mathcal{L},\,\big[\mathcal{L},\,
[\mathcal{L},\overline{\mathcal{L}}]\big]\big],
\\
&
\big[\mathcal{L},\,\big[\overline{\mathcal{L}},\,
[\mathcal{L},\overline{\mathcal{L}}]\big]\big],
\\
&
\big[\overline{\mathcal{L}},\,\big[\mathcal{L},\,
[\overline{\mathcal{L}},\overline{\mathcal{L}}]\big]\big],
\\
&
\big[\overline{\mathcal{L}},\,\big[\overline{\mathcal{L}},\,
[\mathcal{L},\overline{\mathcal{L}}]\big]\big],
\endaligned
\]
but the Jacobi identity yields (exercise, or see below) that the
second equals the third.

Focusing attention on the first, suppose it does {\em not} complete
a frame.

\medskip\noindent{\bf Proposition.}
{\em When simultaneously, the following two degeneracy conditions hold:} 
\[
\aligned
{\bf 4}
&
\,=\,
\rank_\C
\Big(\mathcal{L},\,\overline{\mathcal{L}},\,
[\mathcal{L},\overline{\mathcal{L}}\big],\,\,
\big[\mathcal{L},\,[\mathcal{L},\overline{\mathcal{L}}\big]\big],\,\,
\big[\overline{\mathcal{L}},\,[\mathcal{L},\overline{\mathcal{L}}\big]\big]
\Big)
\\
{\bf 4}
&
=
\rank_\C
\Big(\mathcal{L},\overline{\mathcal{L}},\,
\big[\mathcal{L},\overline{\mathcal{L}}\big],\,\,
\big[\mathcal{L},\,\big[\mathcal{L},\overline{\mathcal{L}}\big]\big],\,\,
\big[\overline{\mathcal{L}},\,
\big[\mathcal{L},\overline{\mathcal{L}}\big]\big],\,\,
\\
&
\ \ \ \ \ \ \ \ \ \ \ \ \ \ \ \ \ \ \ \ \ \ \ \ \ \ \ \ \ \ \ \ \ \ \ \ \ 
\ \ \ \ \ \ \ \ \ \ \ \ \ \ \ \ \ \ \ \ \ \ \ \ \ \ \ \ \ \ \ \ \ \ \ \ \ 
\big[\mathcal{L},\,
\big[\mathcal{L},\,
\big[\mathcal{L},\overline{\mathcal{L}}\big]\big]\big]
\Big),
\endaligned
\]
{\em then:}
\[
\aligned
{\bf 4}
&
=
\rank_\C
\Big(
\LL_{\mathcal{L},\overline{\mathcal{L}}}^4
\Big)
\\
&
=
\rank_\C
\Big(
\LL_{\mathcal{L},\overline{\mathcal{L}}}^5
\Big)
=\cdots=
\\
&
=
\rank_\C
\Big(
\LL_{\mathcal{L},\overline{\mathcal{L}}}^{\sf Lie}
\Big)
\\
&
=
\rank_\R
\Big(
\LL_{{\sf Re}\,\mathcal{L},{\sf Im}\,\mathcal{L}}^{\sf Lie}
\Big)
\,<\,
\,{\bf 5}\,
=
\dim_\R\,M.
\endaligned
\]

\proof
Equivalently, one assumes simultaneously:
\[
\aligned
\big[\overline{\mathcal{L}},[\mathcal{L},\overline{\mathcal{L}}]\big]
&
=
a\cdot\mathcal{L}
+
b\cdot\overline{\mathcal{L}}
+
c\cdot[\mathcal{L},\overline{\mathcal{L}}]
+
d\cdot\big[\mathcal{L},[\mathcal{L},\overline{\mathcal{L}}]\big],
\\
\big[\mathcal{L},\,\big[\mathcal{L},\,
[\mathcal{L},\overline{\mathcal{L}}]\big]\big]
&
=
e\cdot\mathcal{L}
+
f\cdot\overline{\mathcal{L}}
+
g\cdot[\mathcal{L},\overline{\mathcal{L}}]
+
h\cdot\big[\mathcal{L},[\mathcal{L},\overline{\mathcal{L}}]\big].
\endaligned
\]
The claim (proof below) is that this entails:
\[
\aligned
\big[\mathcal{L},\,\big[\overline{\mathcal{L}},\,
[\mathcal{L},\overline{\mathcal{L}}]\big]\big]
&
\equiv
0
\ \ \ \ \
\mod\,\Big(
\mathcal{L},\,
\overline{\mathcal{L}},\,
\big[\mathcal{L},\overline{\mathcal{L}}\big],\,
\big[\mathcal{L},[\mathcal{L},\overline{\mathcal{L}}]\big]
\Big)
\\
\big[\overline{\mathcal{L}},\,\big[\mathcal{L},\,
[\mathcal{L},\overline{\mathcal{L}}]\big]\big]
&
\equiv
0
\ \ \ \ \
\mod\,\Big(
\mathcal{L},\,
\overline{\mathcal{L}},\,
\big[\mathcal{L},\overline{\mathcal{L}}\big],\,
\big[\mathcal{L},[\mathcal{L},\overline{\mathcal{L}}]\big]
\Big)
\\
\big[\overline{\mathcal{L}},\,\big[\overline{\mathcal{L}},\,
[\mathcal{L},\overline{\mathcal{L}}]\big]\big]
&
\equiv
0
\ \ \ \ \
\mod\,\Big(
\mathcal{L},\,
\overline{\mathcal{L}},\,
\big[\mathcal{L},\overline{\mathcal{L}}\big],\,
\big[\mathcal{L},[\mathcal{L},\overline{\mathcal{L}}]\big]
\Big)
\endaligned
\] 
and beyond up to infinity:
\[
\aligned
{\bf 4}
&
=
\rank_\C
\Big(
\LL_{\mathcal{L},\overline{\mathcal{L}}}^4
\Big)
\\
&
=
\rank_\C
\Big(
\LL_{\mathcal{L},\overline{\mathcal{L}}}^5
\Big)
=\cdots=
\\
&
=
\rank_\C
\Big(
\LL_{\mathcal{L},\overline{\mathcal{L}}}^\infty
\Big)
\endaligned
\]
because of the bracketing stability:
\[
\footnotesize
\!\!\!\!\!\!\!\!\!\!
\aligned
&
\Big[\mathcal{L},\,\mod\,\big(\same\big)\Big]
=
\Big[\mathcal{L},\,
\function\cdot\mathcal{L}
+
\function\cdot\overline{\mathcal{L}}
+
\function\cdot\big[\mathcal{L},\overline{\mathcal{L}}\big]
+
\function\cdot\big[\mathcal{L},[\mathcal{L},\overline{\mathcal{L}}]\big]
\Big]
\\
&
=
\mathcal{L}(\function)\cdot\mathcal{L}
+
\mathcal{L}(\function)\cdot\overline{\mathcal{L}}
+
\mathcal{L}(\function)\cdot\big[\mathcal{L},\overline{\mathcal{L}}\big]
+
\mathcal{L}(\function)\cdot
\big[\mathcal{L},[\mathcal{L},\overline{\mathcal{L}}]\big]
+
\\
&
\ \ \ \ \
+
\function\cdot
\zero{[\mathcal{L},\mathcal{L}]}
+
\function\cdot
\big[\mathcal{L},\overline{\mathcal{L}}\big]
+
\function\cdot
\big[\mathcal{L},[\mathcal{L},\overline{\mathcal{L}}]\big]
+
\function\cdot
\!\!\!\!\!\!\!\!
\underbrace{
\big[\mathcal{L},[\mathcal{L},[\mathcal{L},\overline{\mathcal{L}}]]\big]}_{
e\mathcal{L}+f\overline{\mathcal{L}}+
g[\mathcal{L},\overline{\mathcal{L}}]+
h[\mathcal{L},[\mathcal{L},\overline{\mathcal{L}}]]}
\\
&
\equiv
0
\ \ \ \ \
\mod\,\Big(
\mathcal{L},\,
\overline{\mathcal{L}},\,
\big[\mathcal{L},\overline{\mathcal{L}}\big],\,
\big[\mathcal{L},[\mathcal{L},\overline{\mathcal{L}}]\big]
\Big),
\endaligned
\]
and similarly also:
\[
\big[\overline{\mathcal{L}},\,\mod(\same)\big]
\equiv
0
\ \ \ \ \
\mod\,\Big(
\mathcal{L},\,
\overline{\mathcal{L}},\,
\big[\mathcal{L},\overline{\mathcal{L}}\big],\,
\big[\mathcal{L},[\mathcal{L},\overline{\mathcal{L}}]\big]
\Big),
\]
so that induction is clear.

\smallskip

For the claim, observe at first that the Jacobi identity:
\[
0
=
\big[\mathcal{L},\,\big[\overline{\mathcal{L}},\mathcal{T}\big]\big]
+
\!\!\!
\underbrace{
\zero{\big[\mathcal{T},\,\big[\mathcal{L},\overline{\mathcal{L}}\big]\big]}}_{
{\sf remind}\,\,
\mathcal{T}
=
\isqrt[\mathcal{L},\overline{\mathcal{L}}]}
\!\!\!\!
+
\big[\overline{\mathcal{L}},\,\big[\mathcal{T},\mathcal{L}\big]\big]
\]
indeed yields (erase $\isqrt$):
\[
\big[\mathcal{L},\,\big[\overline{\mathcal{L}},\,
[\mathcal{L},\overline{\mathcal{L}}]\big]\big]
=
\big[\overline{\mathcal{L}},\,\big[\mathcal{L},\,
[\overline{\mathcal{L}},\overline{\mathcal{L}}]\big]\big],
\]
so that it suffices to prove the claim only for two of three lines.

\smallskip

Treat the first line of the claim:
\[
\footnotesize
\aligned
\big[\mathcal{L},\,\big[\overline{\mathcal{L}},
[\mathcal{L},\overline{\mathcal{L}}\big]\big]\big]
&
=
\Big[
\mathcal{L},\,\,
a\,\mathcal{L}
+
b\,\overline{\mathcal{L}}
+
c\,\big[\mathcal{L},\overline{\mathcal{L}}\big]
+
d\,\big[\mathcal{L},\,[\mathcal{L},\overline{\mathcal{L}}]\big]
\Big]
\\
&
=
\mathcal{L}(a)\cdot\mathcal{L}
+
\mathcal{L}(b)\cdot\overline{\mathcal{L}}
+
\mathcal{L}(c)\cdot
\big[\mathcal{L},\overline{\mathcal{L}}\big]
+
\mathcal{L}(d)\cdot
\big[\mathcal{L},[\mathcal{L},\overline{\mathcal{L}}]\big]
+
\\
&
\ \ \ \ \
+
\zero{a\,[\mathcal{L},\mathcal{L}]}
+
b\,\big[\mathcal{L},\overline{\mathcal{L}}\big]
+
c\big[\mathcal{L},\,[\mathcal{L},\overline{\mathcal{L}}]\big]
+
d
\!\!\!\!\!
\underbrace{
\big[\mathcal{L},\,\big[\mathcal{L},\,
[\mathcal{L},\overline{\mathcal{L}}]\big]\big]
}_{
e\mathcal{L}+f\overline{\mathcal{L}}+
g[\mathcal{L},\overline{\mathcal{L}}]
+h[\mathcal{L},[\mathcal{L},\overline{\mathcal{L}}]]}
\\
&
\equiv
0
\ \ \ \ \
\mod\,\Big(
\mathcal{L},\,
\overline{\mathcal{L}},\,
\big[\mathcal{L},\overline{\mathcal{L}}\big],\,
\big[\mathcal{L},[\mathcal{L},\overline{\mathcal{L}}]\big]
\Big).
\endaligned
\]
(The same could also be done directly with the second line.)

Treat the third and last line of the claim:
\[
\aligned
\big[\overline{\mathcal{L}},\,
\big[\overline{\mathcal{L}},\,
[\mathcal{L},\overline{\mathcal{L}}]\big]\big]
&
=
-\,
\overline{
\big[\mathcal{L},\,\big[\mathcal{L},\,
[\mathcal{L},\overline{\mathcal{L}}]\big]\big]}
\\
&
=
-\,
\overline{
e\cdot\mathcal{L}
-
f\cdot\overline{\mathcal{L}}
-
g\cdot\big[\mathcal{L},\overline{\mathcal{L}}\big]
-
h\cdot
\big[\mathcal{L},\,[\mathcal{L},\overline{\mathcal{L}}]\big]}
\\
&
=
-\,\overline{e}\cdot\overline{\mathcal{L}}
-
\overline{f}\cdot\mathcal{L}
-
\overline{g}\cdot
\big[\overline{\mathcal{L}},\mathcal{L}\big]
-
\overline{h}\cdot
\!\!\!\!\!\!\!\!\!\!\!\!\!\!
\underbrace{
\big[\overline{\mathcal{L}},[\overline{\mathcal{L}},\mathcal{L}]\big]
}_{
-a\mathcal{L}-b\overline{\mathcal{L}}
-c[\mathcal{L},\overline{\mathcal{L}}]
-d[\mathcal{L},[\mathcal{L},\overline{\mathcal{L}}]]}
\\
&
\equiv
0
\ \ \ \ \
\mod\,\Big(
\mathcal{L},\,
\overline{\mathcal{L}},\,
\big[\mathcal{L},\overline{\mathcal{L}}\big],\,
\big[\mathcal{L},[\mathcal{L},\overline{\mathcal{L}}]\big]
\Big),
\endaligned
\]
which finishes.
\endproof

The generalized CR-Frobenius theorem concludes then that:
\[
M^5
\,\cong\,
\underline{M}^5,
\]
{\em with:}
\[
\underline{M}^5
\,\subset\,
\C^4\times\R
\]
represented in local coordinates as:
\[
\aligned
v_1
&
=
\varphi_1\big(x,y,u_1,u_2,u_3\big),
\\
v_2
&
=
\varphi_2\big(x,y,u_1,u_2,u_3\big),
\\
v_3
&
=
0.
\endaligned
\]
Naturally, one again excludes such a degenerate circumstance.

\medskip

Consequently, it must be that:
\[
\aligned
{\bf 5}
=
\rank_\C
\Big(\mathcal{L},\overline{\mathcal{L}},\,
\big[\mathcal{L},\overline{\mathcal{L}}\big],\,\,
&
\big[\mathcal{L},\,\big[\mathcal{L},\overline{\mathcal{L}}\big]\big],\,\,
\big[\overline{\mathcal{L}},\,
\big[\mathcal{L},\overline{\mathcal{L}}\big]\big],\,\,
\\
&
\big[\mathcal{L},\,
\big[\mathcal{L},\,
\big[\mathcal{L},\overline{\mathcal{L}}\big]\big]\big]
\Big),
\endaligned
\]
and one arrives at a fourth new {\sl general class} of CR-generic manifolds:
\[
\aligned
&
\boxed{\text{\sf General Class $\text{\sf III}_{\text{\sf 2}}$:}}
\\
&
\boxed{\,\,
\aligned
M^5\subset\C^4
\ \
&
\text{\rm with}\ \
\Big\{\mathcal{L},\,\overline{\mathcal{L}},\,\,
\big[\mathcal{L},\overline{\mathcal{L}}\big],\,\,
\big[\mathcal{L},\,
\big[\mathcal{L},\overline{\mathcal{L}}\big]\big],\,\,
\big[\mathcal{L},\,\big[\mathcal{L},\,
\big[\mathcal{L},\overline{\mathcal{L}}\big]\big]\big]
\Big\}\,\,
\\
&
\text{\rm constituting a frame for}\ \
\C\otimes_\R TM,
\\
&
\text{\rm while}\ \
{\bf 4}
=
\rank_\C
\Big(\mathcal{L},\overline{\mathcal{L}},\,
\big[\mathcal{L},\overline{\mathcal{L}}\big],\,\,
\big[\mathcal{L},\,\big[\mathcal{L},\overline{\mathcal{L}}\big]\big],\,\,
\big[\overline{\mathcal{L}},\,
\big[\mathcal{L},\overline{\mathcal{L}}\big]\big]\Big).
\endaligned}
\endaligned
\]

\medskip\noindent
{\bf Theorem.}
{\em Restricting to:}
\[
\dim_\R\,M
\,\leqslant\,
{\bf 5},
\]
{\em this is the last general class in CR dimension:}
\[
n
=
{\bf 1}.
\qed
\]

\medskip

The Class $\text{\sf III}_{\text{\sf 2}}$ model example in
$\C^4 \ni (z, w_1, w_2, w_3)$ is:
\[
\aligned
w_1
-
\overline{w}_1
&
=
2i\,z\overline{z},
\\
w_2
-
\overline{w}_2
&
=
2i\,z\overline{z}\big(z+\overline{z}\big),
\\
w_3
-
\overline{w}_3
&
=
2i\,z\overline{z}
\big(
z^2
+
{\textstyle{\frac{3}{2}}}\,z\overline{z}
+
\overline{z}^2
\big).
\endaligned
\]

%%%%%%%%%%%%%%%%%%%%%%%%%%%%%%%%%%%%%%%%%%%%%%%%%%%%%%%%%%%%%%%%%%%%%

\bigskip

\section{\sf Biholomorphic invariance of the Levi form}
\label{biholomorphic-invariance-LF}
\HEAD{\ref{biholomorphic-invariance-LF}.~Biholomorphic 
invariance of the Levi form}{
Jo\"el {\sc Merker} (Paris-Sud), 
Samuel {\sc Pocchiola} (Paris-Sud), 
Masoud {\sc Sabzevari} (Shahrekord)}

\medskip

Consider a $\mathcal{ C}^\omega$ hypersurface (hence CR-generic):
\[
M^{2n+1}
\subset
\C^{n+1},
\]
with:
\[
\aligned
n
&
=
\CRdim\,M,
\\
1
&
=
\codim\,M.
\endaligned
\]
Pick $p \in M$ and ${\sf U}_p \ni p$ a small open ball
or polydisc.

Suppose a local biholomorphism is given:
\[
h\colon\ \ \ 
{\sf U}_p
\overset{\sim}{\,\longrightarrow\,}
{\sf U}_{p'}'
\]
of ${\sf U}_p$ onto the image open set:
\[
{\sf U}_{p'}' 
\,:=\, 
h\big({\sf U}_p\big)
\subset
\C^{n+1}
\ \ \ \ \ \ \ \ \ \ \ \ \
{\scriptstyle{(p'\,=\,h(p))}}.
\]

\begin{center}
\begin{picture}(0,0)%
\includegraphics{M-hypersurface-polydiscs.pstex}%
\end{picture}%
\setlength{\unitlength}{4144sp}%
\begingroup\makeatletter\ifx\SetFigFont\undefined%
\gdef\SetFigFont#1#2#3#4#5{%
  \reset@font\fontsize{#1}{#2pt}%
  \fontfamily{#3}\fontseries{#4}\fontshape{#5}%
  \selectfont}%
\fi\endgroup%
\begin{picture}(4693,1052)(874,-2319)
\put(2002,-1922){\makebox(0,0)[lb]{\smash{{\SetFigFont{10}{12.0}{\familydefault}{\mddefault}{\updefault}{\color[rgb]{0,0,0}$p$}%
}}}}
\put(4844,-2009){\makebox(0,0)[lb]{\smash{{\SetFigFont{10}{12.0}{\familydefault}{\mddefault}{\updefault}{\color[rgb]{0,0,0}$p'$}%
}}}}
\put(1612,-1507){\makebox(0,0)[lb]{\smash{{\SetFigFont{10}{12.0}{\familydefault}{\mddefault}{\updefault}{\color[rgb]{0,0,0}${\sf U}_p$}%
}}}}
\put(943,-1939){\makebox(0,0)[lb]{\smash{{\SetFigFont{10}{12.0}{\familydefault}{\mddefault}{\updefault}{\color[rgb]{0,0,0}$M$}%
}}}}
\put(4459,-1546){\makebox(0,0)[lb]{\smash{{\SetFigFont{10}{12.0}{\familydefault}{\mddefault}{\updefault}{\color[rgb]{0,0,0}$h({\sf U}_p)$}%
}}}}
\put(3339,-1656){\makebox(0,0)[lb]{\smash{{\SetFigFont{10}{12.0}{\familydefault}{\mddefault}{\updefault}{\color[rgb]{0,0,0}$h$}%
}}}}
\put(5117,-1833){\makebox(0,0)[lb]{\smash{{\SetFigFont{10}{12.0}{\familydefault}{\mddefault}{\updefault}{\color[rgb]{0,0,0}$M'$}%
}}}}
\put(902,-1414){\makebox(0,0)[lb]{\smash{{\SetFigFont{10}{12.0}{\familydefault}{\mddefault}{\updefault}{\color[rgb]{0,0,0}$\C^{n+1}$}%
}}}}
\put(5552,-1470){\makebox(0,0)[lb]{\smash{{\SetFigFont{10}{12.0}{\familydefault}{\mddefault}{\updefault}{\color[rgb]{0,0,0}${\C'}^{n+1}$}%
}}}}
\end{picture}%

\end{center}

Denote the image hypersurface by:
\[
M':=h(M)
\,\subset\C^{n+c}.
\]

\medskip\noindent{\bf Current goal.}
{\em The objective is to compare the Levi forms of $M$ and of $M'$.}

\medskip

Choose local frames:
\[
\aligned
&
\big\{\mathcal{L}_1,\dots,\mathcal{L}_n\big\}
\ \ \ \ \ \ \ \ \ \ 
\text{\rm for}\ \
T^{1,0}M,
\\
&
\big\{\mathcal{L}_1',\dots,\mathcal{L}_n'\big\}
\ \ \ \ \ \ \ \ \ \ 
\text{\rm for}\ \
T^{1,0}M'.
\endaligned
\]
Because:
\[
h_*\big(T^{1,0}M\big)
=
T^{1,0}M',
\]
there must exist $\mathcal{ C}^\omega$ functions
$a_{jk}'$ defined on $M'$ so that (mind indices):
\[
\aligned
h_*\big(\mathcal{L}_1\big)
&
=
a_{11}'\,\mathcal{L}_1'
+\cdots+
a_{n1}'\,\mathcal{L}_n',
\\
\cdots\cdots\cdot
&
\cdots\cdots\cdots\cdots\cdots\cdots\cdots\cdots
\\
h_*\big(\mathcal{L}_n\big)
&
=
a_{1n'}\,\mathcal{L}_1'
+\cdots+
a_{nn}'\,\mathcal{L}_n'.
\endaligned
\]

Simultaneously:
\[
h_*\big(\overline{\mathcal{M}}\big)
=
\overline{h_*\big(\mathcal{M}\big)}
\]
yields:
\[
\aligned
h_*\big(\overline{\mathcal{L}}_1\big)
&
=
\overline{a}_{11}'\,\overline{\mathcal{L}}_1'
+\cdots+
\overline{a}_{n1}'\,\overline{\mathcal{L}}_n',
\\
\cdots\cdots\cdot
&
\cdots\cdots\cdots\cdots\cdots\cdots\cdots\cdots
\\
h_*\big(\overline{\mathcal{L}}_n\big)
&
=
\overline{a}_{1n}'\,\overline{\mathcal{L}}_1'
+\cdots+
\overline{a}_{nn}'\,\overline{\mathcal{L}}_n'.
\endaligned
\]

Of course, the $n \times n$ matrix function (notice the index
transposition):
\[
M'\,\ni\,q'
\,\,\longmapsto\,\,
\left(\!
\begin{array}{ccc}
a_{11}'(q') & \cdots & a_{1n}'(q')
\\
\vdots & \ddots & \vdots
\\
a_{n1}'(q') & \cdots & a_{nn}'(q')
\end{array}
\!\right)
\,\in\,
{\sf GL}_n(\C)
\]
is invertible.

\medskip

Hence at points $q \in M \cap {\sf U}_p$, the aim is to compare
(mind index places):
\[
\mathmotsf{Levi-Matrix}_{\mathcal{L},\overline{\mathcal{L}}}^M(q)
\,:=\,
\left(\!
\begin{array}{ccc}
\rho_0\big(\isqrt\big[\mathcal{L}_1,\overline{\mathcal{L}}_1\big]\big)
& \cdots &
\rho_0\big(\isqrt\big[\mathcal{L}_n,\overline{\mathcal{L}}_1\big]\big)
\\
\vdots & \ddots & \vdots
\\
\rho_0\big(\isqrt\big[\mathcal{L}_1,\overline{\mathcal{L}}_n\big]\big)
& \cdots &
\rho_0\big(\isqrt\big[\mathcal{L}_n,\overline{\mathcal{L}}_n\big]\big)
\end{array}
\!\right)
(q),
\]
where:
\[
\rho_0
\colon\ \ \
TM
\longrightarrow
\R
\]
is any nonzero (local) real $1$-form satisfying:
\[
\big\{\rho_0=0\big\}
=
TM\cap J(TM),
\]
with:
\[
\mathmotsf{Levi-Matrix}_{\mathcal{L}',\overline{\mathcal{L}}'}^{M'}\big(h(q)\big)
\,:=\,
\left(\!
\begin{array}{ccc}
\rho_0'\big(\isqrt\big[\mathcal{L}_1',\overline{\mathcal{L}}_1'\big]\big)
& \cdots &
\rho_0\big(\isqrt\big[\mathcal{L}_n',\overline{\mathcal{L}}_1'\big]\big)
\\
\vdots & \ddots & \vdots
\\
\rho_0\big(\isqrt\big[\mathcal{L}_1',\overline{\mathcal{L}}_n'\big]\big)
& \cdots &
\rho_0\big(\isqrt\big[\mathcal{L}_n',\overline{\mathcal{L}}_n'\big]\big)
\end{array}
\!\right)
\big(h(q)\big),
\]
where similarly:
\[
\rho_0'
\colon\ \ \
TM'
\longrightarrow
\R
\]
satisfies:
\[
\big\{\rho_0'=0\big\}
=
TM'\cap J'(TM').
\]

At first, there exists a nowhere zero function $b'$ with:
\[
\big(h^{-1}\big)^*
(\rho_0)
=
b'\,\rho_0'.
\]
For clarity, it is better to drop $h_*$ symbols:
\[
\left[
\aligned
\rho_0
&
=
b'\,\rho_0',
\\
\mathcal{L}_1
&
=
a_{11}'\,\mathcal{L}_1'
+\cdots+
a_{n1}'\,\mathcal{L}_n',
\\
\cdots
&
\cdots\cdots\cdots\cdots\cdots\cdots\cdots\cdots
\\
\mathcal{L}_n
&
=
a_{1n}'\,\mathcal{L}_1'
+\cdots+
a_{nn}'\,\mathcal{L}_n',
\endaligned\right.
\]
keeping in mind that such equalities are {\em truly satisfied}
after the replacement:
\[
q'
=
h(q)
\]
through which source points are linked to image points.

\smallskip

Hence the thing is to replace all this in the Levi matrix of $M$.

\smallskip

Abbreviating:
\[
\aligned
\mathcal{L}_1
&
=
\sum_{j=1}^n\,a_{j1}'\,\mathcal{L}_j',
\ \ \ \ \ \ \ \ \ \ \ \ \ \ \ \ \ \ \ \ \ \ \ \ \
\overline{\mathcal{L}}_1
=
\sum_{k=1}^n\,\overline{a}_{k1}'\,\overline{\mathcal{L}}_k',
\\
\cdots
&
\cdots\cdots\cdots\cdots\cdot\cdot
\ \ \ \ \ \ \ \ \ \ \ \ \ \ \ \ \ \ \ \ \ \ \ \ \
\cdots
\ \ \
\cdots\cdots\cdots\cdots\cdot\cdot
\\
\mathcal{L}_n
&
=
\sum_{j=1}^n\,a_{jn}'\,\mathcal{L}_j',
\ \ \ \ \ \ \ \ \ \ \ \ \ \ \ \ \ \ \ \ \ \ \ \ \
\overline{\mathcal{L}}_n
=
\sum_{k=1}^n\,\overline{a}_{kn}'\,\overline{\mathcal{L}}_k',
\endaligned
\]
One therefore has to expand:
\[
\aligned
\mathmotsf{Levi-Matrix}_{\mathcal{L},\overline{\mathcal{L}}}^M
\,:=\,
&
\rho_0
\left(\!
\begin{array}{ccc}
\isqrt\big[\mathcal{L}_1,\overline{\mathcal{L}}_1\big]
& \cdots &
\isqrt\big[\mathcal{L}_n,\overline{\mathcal{L}}_1\big]
\\
\vdots & \ddots & \vdots
\\
\isqrt\big[\mathcal{L}_1,\overline{\mathcal{L}}_n\big]
& \cdots &
\isqrt\big[\mathcal{L}_n,\overline{\mathcal{L}}_n\big]
\end{array}
\!\right)
\endaligned
\]
which is:
\[
\small
\aligned
=\,
&
b'\,\rho_0'
\left(\!
\begin{array}{ccc}
\isqrt
\left[
\sum_{j=1}^n\,a_{j1}'\mathcal{L}_j',\,
\sum_{k=1}^n\,\overline{a}_{k1}'\overline{\mathcal{L}}_k'
\right]
& \cdots\cdots &
\isqrt
\left[
\sum_{j=1}^n\,a_{jn}'\mathcal{L}_j',\,
\sum_{k=1}^n\,\overline{a}_{k1}'\overline{\mathcal{L}}_k'
\right]
\\
\vdots & \ddots & \vdots
\\
\isqrt
\left[
\sum_{j=1}^n\,a_{j1}'\mathcal{L}_j',\,
\sum_{k=1}^n\,\overline{a}_{kn}'\overline{\mathcal{L}}_k'
\right]
& \cdots\cdots &
\isqrt
\left[
\sum_{j=1}^n\,a_{jn}'\mathcal{L}_j',\,
\sum_{k=1}^n\,\overline{a}_{kn}'\overline{\mathcal{L}}_k'
\right]
\end{array}
\!\right).
\endaligned
\]
Here, when one expands any appearing Lie bracket:
\[
\bigg[
\sum_{j=1}^n\,a_{jl}'\,\mathcal{L}_j',\,\,
\sum_{k=1}^n\,\overline{a}_{km}'\,\overline{\mathcal{L}}_k'
\bigg]
\,\equiv\,
\sum_{j=1}^n\,\sum_{k=1}^n\,
a_{jl}'\,\overline{a}_{km}'\,
\big[\mathcal{L}_j',\,\overline{\mathcal{L}}_k'\big]
\ \ \ \ \
\mod\,\big(\mathcal{L}_\bullet',\overline{\mathcal{L}}_\bullet'\big),
\]
taking account of:
\[
\rho_0'
\Big(
\mathmotsf{\sf vector}
\ \ 
\mathmotsf{modulo}\ \
\big(\mathcal{L}_\bullet',\overline{\mathcal{L}}_\bullet'\big)
\Big)
\,=\,
\rho_0'
\big(
\mathmotsf{\sf vector}
\big),
\]
one gets as a continuation:
\[
\aligned
=
b'\,\rho_0'
\left(\!
\begin{array}{ccc}
\isqrt\sum_{j=1}^n\sum_{k=1}^na_{j1}'\overline{a}_{k1}'
\big[\mathcal{L}_j',\overline{\mathcal{L}}_k'\big]
& \cdots\cdots &
\isqrt\sum_{j=1}^n\sum_{k=1}^na_{jn}'\overline{a}_{k1}'
\big[\mathcal{L}_j',\overline{\mathcal{L}}_k'\big]
\\
\vdots & \ddots & \vdots
\\
\isqrt\sum_{j=1}^n\sum_{k=1}^na_{j1}'\overline{a}_{kn}'
\big[\mathcal{L}_j',\overline{\mathcal{L}}_k'\big]
& \cdots\cdots &
\isqrt\sum_{j=1}^n\sum_{k=1}^na_{jn}'\overline{a}_{kn}'
\big[\mathcal{L}_j',\overline{\mathcal{L}}_k'\big]
\end{array}
\!\right),
\endaligned
\]
and classically, one may reconstitute the product of 3 matrices:
\[
\!\!\!\!\!\!\!
\small
\aligned
=
b'\,
\underbrace{
\left(\!
\begin{array}{ccc}
\overline{a}_{11}' 
&\cdots &
\overline{a}_{n1}'
\\
\vdots & \ddots & \vdots
\\
\overline{a}_{1n}' 
& \cdots &
\overline{a}_{nn}'
\end{array}
\!\right)}_{
{\sf recognize}\,\,{}^{\tt t}\!\overline{A}'}
\,\,
\underbrace{
\left(\!
\begin{array}{ccc}
\rho_0'\big(\isqrt\big[\mathcal{L}_1',\overline{\mathcal{L}}_1'\big]\big)
& \cdots\cdots &
\rho_0'\big(\isqrt\big[\mathcal{L}_n',\overline{\mathcal{L}}_1'\big]\big)
\\
\vdots & \ddots & \vdots
\\
\rho_0'\big(\isqrt\big[\mathcal{L}_1',\overline{\mathcal{L}}_n'\big]\big)
& \cdots\cdots &
\rho_0'\big(\isqrt\big[\mathcal{L}_n',\overline{\mathcal{L}}_n'\big]\big)
\end{array}
\!\right)}_{
\text{\sf Levi-Matrix}_{\mathcal{L}',\overline{\mathcal{L}}'}^{M'}}
\,\,
\underbrace{
\left(\!
\begin{array}{ccc}
a_{11}' 
&\cdots &
a_{1n}'
\\
\vdots & \ddots & \vdots
\\
a_{n1}' 
&\cdots &
a_{nn}'
\end{array}
\!\right)}_{
{\sf recognize}\,\,A'}.
\endaligned
\]

\medskip\noindent{\bf Conclusion.}
{\em Through any local biholomorphism:}
\[
h
\colon\ \ \
M
\,\longrightarrow\,
M'
\]
{\em between hypersurfaces of $\C^{n+1}$, one has:}
\[
\boxed{\,\,
\mathmotsf{Levi-Matrix}_{\mathcal{L},\overline{\mathcal{L}}}^M(q)
\,=\,
\underbrace{b'\big(h(q)\big)}_{
{\sf nowhere}\,\,0
\atop
{\sf function}}
\,\cdot\,
\underbrace{{}^{\tt t}\!\overline{A}'
\big(h(q)\big)}_{
{\sf invertible}
\atop
{\sf matrix}}
\,\cdot\,
\mathmotsf{Levi-Matrix}_{\mathcal{L}',\overline{\mathcal{L}}'}^{M'}\big(h(q)\big)
\,\cdot\,
\underbrace{A'\big(h(q)\big)}_{
{\sf invertible}
\atop
{\sf matrix}},\,\,}
\]
{\em and moreover:}
\[
\boxed{\,\,
\rank_\C
\Big(
\mathmotsf{Levi-Matrix}_{\mathcal{L},\overline{\mathcal{L}}}^M(q)
\Big)
\,=\,
\rank_\C
\Big(
\mathmotsf{Levi-Matrix}_{\mathcal{L}',\overline{\mathcal{L}}'}^{M'}
\big(h(q)\big)
\Big),\,\,}
\]
{\em for every point $q \in M \cap {\sf U}_p$.\qed}

\medskip

Importantly, if follows from the latter conclusion that:

\medskip\noindent{\bf Scholium.}
{\em The rank of the Levi form of a hypersurface
$M^{ 2n+1} \subset \C^{n+1}$ at one
of its points is independent
both of the choice of local coordinates, and
of the choice of a local frame for $T^{1, 0}M$.\qed}

\medskip

This is why one will allow to employ the lightedned notation:
\[
\mathmotsf{Levi-Form}^M
\]
to emphasize the invariant features of the Levi form.

\medskip\noindent{\bf Yet a bit more about the Levi form.}
In the local frame:
\[
\big\{
\mathcal{L}_1,\dots,\mathcal{L}_n
\big\}
\]
for $T^{1, 0} M$, at a point $q \in M \cap {\sf U}_p$, 
pick two constant vectors:
\[
\aligned
\mathcal{M}_q
&
=
\mu_{1q}\,\mathcal{L}_1\big\vert_q
+\cdots+
\mu_{nq}\,\mathcal{L}_n\big\vert_q,
\\
\mathcal{N}_q
&
=
\nu_{1q}\,\mathcal{L}_1\big\vert_q
+\cdots+
\nu_{nq}\,\mathcal{L}_n\big\vert_q,
\endaligned
\]
and define as a matrix triple product:
\[
\small
\aligned
\!\!\!\!\!\!\!\!\!\!\!\!\!\!\!\!\!\!\!\!\!\!\!\!\!\!\!\!\!\!\!\!
\mathmotsf{Levi-Form}_{\mathcal{L},\overline{\mathcal{L}}}^{M,q}\!\!
\left(\!\!\!
\left(\!\!\!
\begin{array}{c}
\mu_{1q}
\\
\vdots
\\
\mu_{nq}
\end{array}
\!\!\!\right)\!,\!
\left(\!\!\!
\begin{array}{c}
\nu_{1q}
\\
\vdots
\\
\nu_{nq}
\end{array}
\!\!\!\right)
\!\!\!\right)
:=
\big(\overline{\nu}_{1q},\dots,\overline{\nu}_{nq}\big)\!\!
\left(\!\!\!
\begin{array}{ccc}
\rho_0\big(\isqrt\big[\mathcal{L}_1,\overline{\mathcal{L}}_1\big]\big)(q)
& \cdots &
\rho_0\big(\isqrt\big[\mathcal{L}_n,\overline{\mathcal{L}}_1\big]\big)(q)
\\
\vdots & \ddots & \vdots
\\
\rho_0\big(\isqrt\big[\mathcal{L}_1,\overline{\mathcal{L}}_n\big]\big)(q)
& \cdots &
\rho_0\big(\isqrt\big[\mathcal{L}_n,\overline{\mathcal{L}}_n\big]\big)(q)
\end{array}
\!\!\!\right)
\!\!\!
\left(\!\!\!
\begin{array}{c}
\mu_{1q}
\\
\vdots
\\
\mu_{nq}
\end{array}
\!\!\!\right)
\endaligned
\]
which matches up with all what precedes and which shows, once
more, that the value depends only on the two vectors
at $q$.

\medskip\noindent{\bf Kernel of the Levi form and its biholomorphic
invariance.}
At a point $q \in M \cap {\sf U}_p$, consider a vector:
\[
\mathcal{K}_q
=
\kappa_{1q}\,\mathcal{L}_1\big\vert_q
+\cdots+
\kappa_{nq}\,\mathcal{L}_n\big\vert_q,
\]
with constants $\kappa_{ 1q}, \dots, \kappa_{ nq} \in \C$.

\medskip\noindent{\bf Definition.}
Such a vector $\mathcal{ K}_q$ is said to {\sl belong to the
kernel of the Levi form} when:
\[
0
\,=\,
\mathmotsf{Levi-Form}_{\mathcal{L},\overline{\mathcal{L}}}^{M,q}
\left(\!
\left(\!
\begin{array}{c}
\kappa_{1q}
\\
\vdots
\\
\kappa_{nq}
\end{array}
\!\right),
\,\,
\left(\!
\begin{array}{c}
\nu_{1q}
\\
\vdots
\\
\nu_{nq}
\end{array}
\!\right)
\!\right),
\]
for every:
\[
\big(\nu_{1q},\dots,\nu_{nq}\big)
\,\in\,\C^n,
\]
that is to say equivalently, when:
\[
\left(\!
\begin{array}{c}
\kappa_{1q}
\\
\vdots
\\
\kappa_{nq}
\end{array}
\!\right)
\,\in\,
\mathmotsf{Kernel}\,
\Big(
\mathmotsf{Levi-Matrix}_{\mathcal{L},\overline{\mathcal{L}}}^M
(q)
\Big).
\]

\medskip

{\em Passim}, recall the known fact that:
\[
\mathmotsf{Kernel}
\,\subset\,
\mathmotsf{Isotropic cone},
\]
which means, chosing plainly:
\[
\big(\nu_{1q},\dots,\nu_{nq}\big)
:=
\big(\kappa_{1q},\dots,\kappa_{nq}\big),
\]
that:
\[
0
\,=\,
\mathmotsf{Levi-Form}_{\mathcal{L},\overline{\mathcal{L}}}^{M,q}
\left(\!
\left(\!
\begin{array}{c}
\kappa_{1q}
\\
\vdots
\\
\kappa_{nq}
\end{array}
\!\right),
\,\,
\left(\!
\begin{array}{c}
\kappa_{1q}
\\
\vdots
\\
\kappa_{nq}
\end{array}
\!\right)
\!\right).
\]

\medskip

Next, examine how Levi kernels transfer through biholomorphisms.

\smallskip

Thus, assume:
\[
0
=
\underbrace{
\mathmotsf{Levi-Matrix}_{\mathcal{L},\overline{\mathcal{L}}}^M(q)}_{
\sf replace}
\cdot
\left(\!
\begin{array}{c}
\kappa_{1q}
\\
\vdots
\\
\kappa_{nq}
\end{array}
\!\right),
\]
and replace, or {\em transfer}, the Levi matrix:
\[
0
=
\underbrace{b'\big(h(q)\big)}_{\text{\sf nowhere}\,\,0}
\cdot
\underbrace{{}^{\tt t}\!\overline{A}'\big(h(q)\big)}_{
\text{\sf invertible}}
\cdot\,
\mathmotsf{Levi-Matrix}_{\mathcal{L}',\overline{\mathcal{L}}'}^{M'}
\big(h(q)\big)
\cdot
A'\big(h(q)\big)
\cdot
\left(\!
\begin{array}{c}
\kappa_{1q}
\\
\vdots
\\
\kappa_{nq}
\end{array}
\!\right),
\]
that is to say after simplification:
\[
0
=
\mathmotsf{Levi-Matrix}_{\mathcal{L}',\overline{\mathcal{L}}'}^{M'}
\big(h(q)\big)
\cdot
\left(\!
\begin{array}{ccc}
a_{11}'(h(q)) 
& \cdots &
a_{1n}'(h(q))
\\
\vdots & \ddots & \vdots
\\
a_{n1}'(h(q)) 
& \cdots &
a_{nn}'(h(q))
\end{array}
\!\right)
\cdot
\left(\!
\begin{array}{c}
\kappa_{1q}
\\
\vdots
\\
\kappa_{nq}
\end{array}
\!\right).
\]

\medskip\noindent{\bf Natural proposition.}
{\em At an arbitrary point $q \in M \cap {\sf U}_p$, a vector:}
\[
\mathcal{K}_q
=
\kappa_{1q}\,\mathcal{L}_1\big\vert_q
+\cdots+
\kappa_{nq}\,\mathcal{L}_n\big\vert_q
\,\in\,
\mathmotsf{Kernel}
\Big(
\mathmotsf{Levi-Matrix}_{\mathcal{L},\overline{\mathcal{L}}}^M(q)
\Big)
\]
{\em belongs to the kernel of the source Levi form
{\em if and only if} its transferred image:}
\[
\aligned
h_*\big(\mathcal{K}_q\big)
&
=
\big(a_{11}'(h(q))\,\kappa_{1q}
+\cdots+
a_{1n}'(h(q))\,\kappa_{nq}
\big)\,
\mathcal{L}_1'\big\vert_{h(q)}
+
\\
&
\ \ \ \ \
+
\cdots\cdots\cdots\cdots\cdots\cdots\cdots\cdots\cdots\cdots\cdots
\cdots\cdots
+
\\
&
\ \ \ \ \
+
\big(a_{n1}'(h(q))\,\kappa_{1q}
+\cdots+
a_{nn}'(h(q))\,\kappa_{nq}
\big)\,
\mathcal{L}_n'\big\vert_{h(q)}
\endaligned
\]
{\em belongs to the kernel of the target Levi form:}
\[
h_*\big(\mathcal{K}_q\big)
\,\in\,
\mathmotsf{Kernel}
\Big(
\mathmotsf{Levi-Matrix}_{\mathcal{L}',\overline{\mathcal{L}}'}^{M'}
\big(h(q)\big)
\Big).
\qed
\]

%%%%%%%%%%%%%%%%%%%%%%%%%%%%%%%%%%%%%%%%%%%%%%%%%%%%%%%%%%%%%%%%%%%%%

\bigskip

\section{\sf Levi kernel and Freeman form 
\\
in CR dimension $n = {\bf 2}$}
\label{freeman-n-2}
\HEAD{\ref{freeman-n-2}.~Levi kernel and Freeman form
nondegeneracies in CR dimension $n = {\bf 2}$}{
Jo\"el {\sc Merker} (Paris-Sud), 
Samuel {\sc Pocchiola} (Paris-Sud), 
Masoud {\sc Sabzevari} (Shahrekord)}

\medskip

Now, consider a connected $\mathcal{ C}^\omega$ hypersurface:
\[
M^5
\subset
\C^3,
\]
hence:
\[
\aligned
c
&
=
1,
\\
n
&
=
2.
\endaligned
\]
Let $p \in M$, let ${\sf U}_p \ni p$ be a small open ball, 
and let:
\[
\big\{
\mathcal{L}_1,\,\mathcal{L}_2
\big\}
\]
be a local frame for $T^{ 1, 0} M$.

Also, in terms of a differential $1$-form:
\[
\rho_0
\colon\ \ \
TM
\,\longrightarrow\,
\R
\]
satisfying:
\[
\big\{\rho_0=0\big\}
=
TM\cap J(TM),
\]
at every point $q \in M \cap {\sf U}_p$, abbreviate:
\[
\aligned
\mathmotsf{Levi-Matrix}_{\mathcal{L},\overline{\mathcal{L}}}^M(q)
&
=
\left(\!
\begin{array}{cc}
\rho_0\big(\isqrt\,\big[\mathcal{L}_1,
\overline{\mathcal{L}}_1\big]\big)
&
\rho_0\big(\isqrt\,\big[\mathcal{L}_2,
\overline{\mathcal{L}}_1\big]\big)
\\
\rho_0\big(\isqrt\,\big[\mathcal{L}_1,
\overline{\mathcal{L}}_2\big]\big)
&
\rho_0\big(\isqrt\,\big[\mathcal{L}_2,
\overline{\mathcal{L}}_2\big]\big)
\end{array}
\!\right)
(q)
\\
&
=:
\left(\!
\begin{array}{cc}
\ell_{11}(q) & \ell_{12}(q)
\\
\ell_{21}(q) & \ell_{22}(q)
\end{array}
\!\right),
\endaligned
\]
in terms of a $2 \times 2$ matrix-valued $\mathcal{ C}^\omega$ function:
\[
q
\,\longmapsto\,
\left(\!
\begin{array}{cc}
\ell_{11}(q) & \ell_{12}(q)
\\
\ell_{21}(q) & \ell_{22}(q)
\end{array}
\!\right).
\]

At any point $q \in M \cap {\sf U}_p$:
\[
\mathmotsf{possible ranks}_\C
\,=\,
{\bf 0},\ \ 
{\bf 1},\ \ 
{\bf 2}.
\]

Assume temporarily that the:
\[
\aligned
\mathmotsf{Levi-Determinant}
(q)
:=
&\,
\ell_{11}(q)\,\ell_{22}(q)
-
\ell_{12}(q)\,\ell_{21}(q)
\\
\not\equiv
&\,
0
\endaligned
\]
is not identically zero as a $\mathcal{ C}^\omega$ function
of $q \in M \cap {\sf U}_p$.

Introducing then the {\em proper} real analytic subset:
\[
\Sigma_p
\,:=\,
\big\{
q\in M\cap{\sf U}_p\colon\,
0
=
\big(\ell_{11}\,\ell_{22}-\ell_{12}\,\ell_{21}\big)(q)
\big\},
\]
one by definition has then for every:
\[
q\in
\big(M\cap{\sf U}_p\big)
\big\backslash
\Sigma_p
\]
that:
\[
{\bf 2}
=
\rank_\C
\Big(
\mathmotsf{Levi-Matrix}_{\mathcal{L},\overline{\mathcal{L}}}^M(q)
\Big).
\]

Although the following fact is well known, it is advisable to spend
some energy in explaining it.

\medskip\noindent{\bf Assertion.}
{\em If the Levi form of
a \underline{connected} hypersurface $M^5 \subset \C^3$ 
is of rank ${\bf 2}$ at one of its points, then:}
\[
\Sigma
:=
\big\{
q\in M\colon\,
\rank_\C
\big(
\mathmotsf{Levi-Form}^M(q)
\big)
\leqslant
{\bf 1}
\big\}
\]
{\em is a \underline{global}, \underline{proper} $\mathcal{ C}^\omega$
subset of $M$, so that the Levi form
is of rank $2$ at almost every point of $M$.}

\proof
Let $p \in M$, let ${\sf U}_p \ni p$ be a small
ball in which the equation of $M$ is locally
expandable in convergent Taylor series.

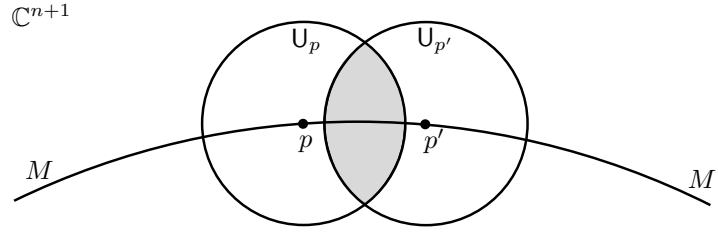
\begin{center}
\begin{picture}(0,0)%
\includegraphics{2-balls.pstex}%
\end{picture}%
\setlength{\unitlength}{4144sp}%
\begingroup\makeatletter\ifx\SetFigFont\undefined%
\gdef\SetFigFont#1#2#3#4#5{%
  \reset@font\fontsize{#1}{#2pt}%
  \fontfamily{#3}\fontseries{#4}\fontshape{#5}%
  \selectfont}%
\fi\endgroup%
\begin{picture}(4212,1336)(1181,-2510)
\put(1196,-1300){\makebox(0,0)[lb]{\smash{{\SetFigFont{10}{12.0}{\familydefault}{\mddefault}{\updefault}{\color[rgb]{0,0,0}$\C^{n+1}$}%
}}}}
\put(2918,-2026){\makebox(0,0)[lb]{\smash{{\SetFigFont{10}{12.0}{\familydefault}{\mddefault}{\updefault}{\color[rgb]{0,0,0}$p$}%
}}}}
\put(1277,-2219){\makebox(0,0)[lb]{\smash{{\SetFigFont{10}{12.0}{\familydefault}{\mddefault}{\updefault}{\color[rgb]{0,0,0}$M$}%
}}}}
\put(3618,-1423){\makebox(0,0)[lb]{\smash{{\SetFigFont{10}{12.0}{\familydefault}{\mddefault}{\updefault}{\color[rgb]{0,0,0}${\sf U}_{p'}$}%
}}}}
\put(5242,-2273){\makebox(0,0)[lb]{\smash{{\SetFigFont{10}{12.0}{\familydefault}{\mddefault}{\updefault}{\color[rgb]{0,0,0}$M$}%
}}}}
\put(2867,-1426){\makebox(0,0)[lb]{\smash{{\SetFigFont{10}{12.0}{\familydefault}{\mddefault}{\updefault}{\color[rgb]{0,0,0}${\sf U}_p$}%
}}}}
\put(3665,-2036){\makebox(0,0)[lb]{\smash{{\SetFigFont{10}{12.0}{\familydefault}{\mddefault}{\updefault}{\color[rgb]{0,0,0}$p'$}%
}}}}
\end{picture}%

\end{center}

Let $p' \in M$ be another point, let ${\sf U}_{ p'} \ni p'$ be another
small ball and suppose: 
\[
{\sf U}_p
\cap
{\sf U}_{p'}
\neq
\emptyset.
\]

Take affine coordinates:
\[
(z_1,z_2,w)
\]
centered at $p$ and affine coordinates:
\[
(z_1',z_2',w')
\]
centered at $p'$, so that:
\[
\aligned
(z_1',z_2',w')
&
=
\mathmotsf{affine}
\big(z_1,z_2,w\big)
\\
&
=:
h(z_1,z_2,w).
\endaligned
\]

Let $\big\{ \mathcal{ L}_1, \mathcal{ L}_2 \big\}$ be a local frame in
${\sf U}_p$ for $T^{ 1, 0} M$ having $\mathcal{ C}^\omega$
coefficients.

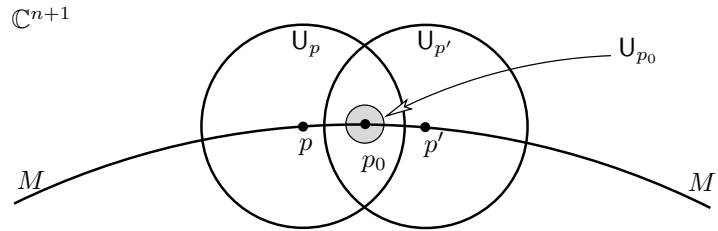
\begin{center}
\begin{picture}(0,0)%
\includegraphics{3-balls.pstex}%
\end{picture}%
\setlength{\unitlength}{4144sp}%
\begingroup\makeatletter\ifx\SetFigFont\undefined%
\gdef\SetFigFont#1#2#3#4#5{%
  \reset@font\fontsize{#1}{#2pt}%
  \fontfamily{#3}\fontseries{#4}\fontshape{#5}%
  \selectfont}%
\fi\endgroup%
\begin{picture}(4212,1357)(1181,-2510)
\put(1196,-1300){\makebox(0,0)[lb]{\smash{{\SetFigFont{10}{12.0}{\familydefault}{\mddefault}{\updefault}{\color[rgb]{0,0,0}$\C^{n+1}$}%
}}}}
\put(5242,-2273){\makebox(0,0)[lb]{\smash{{\SetFigFont{10}{12.0}{\familydefault}{\mddefault}{\updefault}{\color[rgb]{0,0,0}$M$}%
}}}}
\put(2867,-1426){\makebox(0,0)[lb]{\smash{{\SetFigFont{10}{12.0}{\familydefault}{\mddefault}{\updefault}{\color[rgb]{0,0,0}${\sf U}_p$}%
}}}}
\put(3665,-2036){\makebox(0,0)[lb]{\smash{{\SetFigFont{10}{12.0}{\familydefault}{\mddefault}{\updefault}{\color[rgb]{0,0,0}$p'$}%
}}}}
\put(2918,-2026){\makebox(0,0)[lb]{\smash{{\SetFigFont{10}{12.0}{\familydefault}{\mddefault}{\updefault}{\color[rgb]{0,0,0}$p$}%
}}}}
\put(3293,-2126){\makebox(0,0)[lb]{\smash{{\SetFigFont{10}{12.0}{\familydefault}{\mddefault}{\updefault}{\color[rgb]{0,0,0}$p_0$}%
}}}}
\put(3618,-1423){\makebox(0,0)[lb]{\smash{{\SetFigFont{10}{12.0}{\familydefault}{\mddefault}{\updefault}{\color[rgb]{0,0,0}${\sf U}_{p'}$}%
}}}}
\put(4827,-1468){\makebox(0,0)[lb]{\smash{{\SetFigFont{10}{12.0}{\familydefault}{\mddefault}{\updefault}{\color[rgb]{0,0,0}${\sf U}_{p_0}$}%
}}}}
\put(1227,-2259){\makebox(0,0)[lb]{\smash{{\SetFigFont{10}{12.0}{\familydefault}{\mddefault}{\updefault}{\color[rgb]{0,0,0}$M$}%
}}}}
\end{picture}%

\end{center}

Let also $\big\{ \mathcal{ L}_1', \mathcal{ L}_2' \big\}$
be a local frame in ${\sf U}_{p'}$ for $T^{ 1, 0} M$ having 
$\mathcal{ C}^\omega$
coefficients, so that:
\[
\aligned
\mathcal{L}_1
&
=
a_{11}'\,\mathcal{L}_1'
+
a_{21}'\,\mathcal{L}_2',
\\
\mathcal{L}_2
&
=
a_{12}'\,\mathcal{L}_1'
+
a_{22}'\,\mathcal{L}_2',
\endaligned
\]
on the intersection. Take in fact to fix ideas:
\[
p_0
\in
{\sf U}_p
\cap
{\sf U}_{p'},
\]
and a small open ball:
\[
{\sf U}_{p_0}
\subset
{\sf U}_p
\cap
{\sf U}_{p'}.
\]

Of course:
\[
0
\neq
\det
\left(\!
\begin{array}{cc}
a_{11}'(q') & a_{12}'(q')
\\
a_{21}'(q') & a_{22}'(q')
\end{array}
\!\right),
\]
for $q'$ in the intersection domain.

In terms of this invertible matrix and of two differential
$1$-forms $\rho_0$, $\rho_0'$, one already knows that:
\[
\left(\!
\begin{array}{cc}
\ell_{11}(q) & \ell_{12}(q)
\\
\ell_{21}(q) & \ell_{22}(q)
\end{array}
\!\right)
=
\left(\!
\begin{array}{c}
\mathmotsf{invertible}
\\
\mathmotsf{matrix}
\end{array}
\!\right)
\!\!\cdot\!\!
\left(\!
\begin{array}{cc}
\ell_{11}'(h(q)) & \ell_{12}'(h(q))
\\
\ell_{21}'h((q)) & \ell_{22}'h((q))
\end{array}
\!\right)
\!\!\cdot\!\!
\left(\!
\begin{array}{c}
\mathmotsf{invertible}
\\
\mathmotsf{matrix}
\end{array}
\!\right),
\]
which implies that the two real analytic subsets:
\[
\aligned
\Sigma_p
&
:=
\big\{
q\in M\cap{\sf U}_p\colon\,
0
=
(\ell_{11}\,\ell_{22}-\ell_{12}\,\ell_{21})(q)
\big\},
\\
\Sigma_{p'}
&
:=
\big\{
q'\in M\cap{\sf U}_{p'}\colon\,
0
=
(\ell_{11}'\,\ell_{22}'-\ell_{12}'\,\ell_{21}')(q')
\big\},
\endaligned
\]
{\em coincide} on ${\sf U}_{ p_0}$. Thus they glue
together, and from point to point $p, p', p'', \dots \in M$
(use connectedness),
all $\Sigma_p, \Sigma_{ p'}, \Sigma_{ p''}, \dots$
glue alltogether 
as a global real analytic subset $\Sigma \subset M$.

The uniqueness principle on {\em connected} open sets:
\[
\mathmotsf{analytic functions are}
\left\{
\aligned
&
\mathmotsf{either}\,\,\equiv 0,
\\
&
\mathmotsf{or almost everywhere}\,\,\neq 0,
\endaligned\right.
\]
then insures that:
\[
\aligned
\Big[
\big(\ell_{11}\,\ell_{22}
-
\ell_{12}\,\ell_{21}\big)(q)
\not\equiv
0
\ \
\text{\rm on}\,\,
M\cap{\sf U}_p
\Big]
&
\,\,\Longrightarrow\,\,
\Big[
\big(\ell_{11}'\,\ell_{22}'
-
\ell_{12}'\,\ell_{21}'\big)(q')
\not\equiv
0
\ \
\text{\rm on}\,\,
M\cap{\sf U}_{p_0}
\Big]
\\
&
\,\,\Longrightarrow\,\,
\Big[
\big(\ell_{11}'\,\ell_{22}'
-
\ell_{12}'\,\ell_{21}'\big)(q')
\not\equiv
0
\ \
\text{\rm on}\,\,
M\cap{\sf U}_{p'}
\Big]
\endaligned
\]
(the converse is also trivially true), so that:
\[
\aligned
\Sigma_p
\cap\big(M\cap{\sf U}_p\big)
\,\,\text{\rm is proper}
&
\,\,\Longrightarrow\,\,
\Sigma_p
\cap\big(M\cap{\sf U}_p\cap{\sf U}_{p'}\big)
\,\,\text{\rm is proper}
\\
&
\,\,\Longrightarrow\,\,
\Sigma_p
\cap\big(M\cap{\sf U}_{p_0}\big)
\,\,\text{\rm is proper}
\\
\explain{$\Sigma_{p'} = \Sigma_p$ inside ${\sf U}_{p_0}$}
\ \ \ \ \ \ \ \ \ \ \ \ \ \ \ \
&
\,\,\Longrightarrow\,\,
\Sigma_{p'}
\cap\big(M\cap{\sf U}_{p_0}\big)
\,\,\text{\rm is proper}
\\
&
\,\,\Longrightarrow\,\,
\Sigma_{p'}
\cap\big(M\cap{\sf U}_{p'}\big)
\,\,\text{\rm is proper},
\endaligned
\]
which shows by connectedness 
of $M$ that the glued global $\Sigma$ is nowhere dense
(proper) as soon as one $\Sigma_p$ is.
\endproof

{\em Passim}, a known generalized statement can be mentioned.

\medskip\noindent{\bf Theorem.}
{\em On a $\mathcal{ C}^\omega$ connected hypersurface:}
\[
M^{2n+1}
\subset
\C^{n+1},
\]
{\em there exists an integer:}
\[
\aligned
r_M
&
=
\mathmotsf{Zariski-generic (maximal possible) rank}
\\
&
\ \ \ \ \
\mathmotsf{of the Levi form of}\,\,M
\\
&
=:
\genrank_\C\big(\mathmotsf{Levi-Form}^M\big),
\endaligned
\]
{\em satisfying:}
\[
0
\leqslant
r_M
\leqslant
n
=
\CRdim\,M,
\]
{\em and there exists a global proper real analytic
($\mathcal{ C}^\omega$) subset:}
\[
\Sigma
\subset
M
\]
{\em such that:}
\[
M\backslash\Sigma
\ni
p
\,\,\Longleftrightarrow\,\,
r_M
=
\rank_\C
\big(
\mathmotsf{Levi-Form}^M(p)
\big),
\]
{\em which means equivalently:}
\[
\Sigma
\ni
p
\,\,\Longleftrightarrow\,\,
\rank_\C
\big(
\mathmotsf{Levi-Form}^M(p)
\big)
\leqslant
\genrank_\C
\big(\mathmotsf{Levi-Form}^M\big)
-
1.
\qed
\]

\medskip

For a hypersurface $M^5 \subset \C^3$:
\[
r_M
=
\left\{
\aligned
&
{\bf 0}\colon\ \ \
\text{\rm Levi-flat mostly degenerate case}\,\,
M\cong\C^2\times\R.
\\
&
{\bf 1}\colon\ \ \
\text{\rm Intermediate case to be examined}.
\\
&
{\bf 2}\colon\ \ \
\text{\rm Anciently known Levi nondegenerate case}.
\endaligned\right.
\]

\medskip

Admitting as before the:
\[
\text{\sl Lie-Cartan Principle of Relocalization},
\]
one arrives at a well known fifth 
{\sl general class} of CR-generic manifolds:
\[
\aligned
&
\boxed{\text{\sf General Class $\text{\sf IV}_{\text{\sf 1}}$:}}
\\
&
\boxed{\,\,
\aligned
M^5\subset\C^3
\ \
&
\text{\rm with}\ \
\Big\{\mathcal{L}_1,\,\mathcal{L}_2,\,
\overline{\mathcal{L}}_1,\,\overline{\mathcal{L}}_2,\,\,
\big[\mathcal{L}_1,\overline{\mathcal{L}}_1\big]
\Big\}\,\,
\\
&
\text{\rm constituting a frame for}\ \
\C\otimes_\R TM,
\\
&
\text{\rm and with the Levi-Form:}\ \
\\
&
\ \ \ \ \ \ \ \ \ \ \ \ \ \ \
\mathmotsf{Levi-Form}^M(p)
\\
&
\text{\rm being of rank}\,\,{\bf 2}\,\,
\text{\rm at every point}\,\,p\in M.
\endaligned}
\endaligned
\]

\medskip\noindent{\bf Hypersurfaces $M^5 \subset \C^3$ having
Levi Form everywhere of rank 1.}
Examine now the circumstance where:
\[
\rank_\C\big(
\mathmotsf{Levi-Form}^M(p)
\big)
=
{\bf 1},
\]
at every point $p$ of a connected hypersurface
$M^5 \subset \C^3$. 

\medskip\noindent{\bf Assertion.}
{\em If $M^5 \subset \C^3$ is of $\mathcal{ C}^\kappa$ smoothness
with $\kappa \geqslant 2$, or $\mathcal{ C}^\infty$, 
or $\mathcal{ C}^\omega$, then there exists a unique
$\C$-vector subbundle:}
\[
K^{1,0}M
\subset
T^{1,0}M
\]
{\em of $\mathcal{ C}^{\kappa - 2}$ smoothness, or $\mathcal{ C}^\infty$, 
or $\mathcal{ C}^\omega$, having:}
\[
\rank_\C\big(K^{1,0}M\big)
=
{\bf 1},
\]
{\em such that, at every point $p \in M$:}
\[
K_p^{1,0}M
\,\ni\,
\mathcal{K}_p
\,\,\Longleftrightarrow\,\,
\mathcal{K}_p
\,\in\,
\mathmotsf{Kernel}
\big(
\mathmotsf{Levi-Form}^M(p)
\big).
\]
{\em In other words, this means that the union of all the
$1$-dimensional kernels gathers coherently and
smoothly to constitute a true subbundle of
$T^{1, 0} M$.}

\proof
Let $p \in M$, let ${\sf U}_p \ni p$ be a small open ball or polydisc,
let:
\[
\big\{
\mathcal{L}_1,\mathcal{L}_2
\big\}
\]
be a local frame for $T^{1, 0} M$, and pick a local 
real differential $1$-form: 
\[
\rho_0
\colon\ \ \
TM
\,\longrightarrow\,
\R
\]
whose extension to $\C \otimes_\R TM$ satisfies:
satisfying:
\[
\big\{\rho_0=0\big\}
=
T^{0,1}M\oplus T^{0,1}M.
\]
One should look at the kernels of the Levi Matrices:
\[
\mathmotsf{Kernel}\,
\left(\!
\begin{array}{cc}
\rho_0\big(\isqrt\big[\mathcal{L}_1,\overline{\mathcal{L}}_1\big]\big)
& 
\rho_0\big(\isqrt\big[\mathcal{L}_2,\overline{\mathcal{L}}_1\big]\big)
\\
\rho_0\big(\isqrt\big[\mathcal{L}_1,\overline{\mathcal{L}}_2\big]\big)
& 
\rho_0\big(\isqrt\big[\mathcal{L}_2,\overline{\mathcal{L}}_2\big]\big)
\end{array}
\!\right)(q)
\]
at various points $q \in M \cap {\sf U}_p$, which one abbreviates as:
\[
\mathmotsf{Kernel}\,
\left(\!
\begin{array}{cc}
\ell_{11}(q) & \ell_{12}(q)
\\
\ell_{21}(q) & \ell_{22}(q)
\end{array}
\!\right),
\]
with entry functions $\ell_{ 11}$, $\ell_{ 12}$, $\ell_{21}$, $\ell_{22}$
being $\mathcal{ C}^{\kappa-2}$, or $\mathcal{ C}^\infty$, or
$\mathcal{ C}^\omega$.

By assumption, at $p$, one entry is nonzero:
\[
\ell_{11}(p)\neq 0,
\ \ \ \ \ \ \
\text{\rm or}
\ \ \ \ \ \ \
\ell_{12}(p)\neq 0,
\ \ \ \ \ \ \
\text{\rm or}
\ \ \ \ \ \ \
\ell_{21}(p)\neq 0,
\ \ \ \ \ \ \
\text{\rm or}
\ \ \ \ \ \ \
\ell_{22}(p)\neq 0,
\]
because at no point of $M$, the rank $1$ Levi Form can be zero.

It is easy to check that there exists a constant matrix:
\[
\left(\!
\begin{array}{cc}
\alpha & \beta
\\
\gamma & \delta
\end{array}
\!\right)
\,\in\,
{\sf GL}_2(\C)
\]
so that, replacing the $T^{1,0}M$-frame by:
\[
\aligned
\mathcal{L}_1^\sharp
&
:=
\alpha\,\mathcal{L}_1
+
\beta\,\mathcal{L}_2,
\\
\mathcal{L}_2^\sharp
&
:=
\gamma\,\mathcal{L}_1
+
\delta\,\mathcal{L}_2,
\endaligned
\]
one can assume (dropping ${}^\sharp$ symbols) that:
\[
0
\,\neq\,
\rho_0\big(\isqrt\big[\mathcal{L}_1,\overline{\mathcal{L}}_1\big](p)
=
\ell_{11}(p),
\]
so that, shrinking ${\sf U}_p$ if necessary, one has
by continuity:
\[
0
\neq
\ell_{11}(q)
\]
for every $q \in M \cap {\sf U}_p$. 

On the other hand, the rank of the Levi form being
nowhere equal to ${\bf 2}$ by hypothesis, one must have:
\[
\aligned
0
&
\equiv
\mathmotsf{Levi-Determinant}
\\
&
=
\ell_{11}(q)\,\ell_{22}(q)
-
\ell_{12}(q)\,\ell_{21}(q)
\ \ \ \ \ \ \ \ \ \ \ \ \
{\scriptstyle{(\forall\,q\,\in\,M\,\cap\,{\sf U}_p)}}.
\endaligned
\]

Consequently, by plain linear algebra:
\[
{\bf 1}
=
\dim_\C\Big(
\mathmotsf{Kernel}\big(
\mathmotsf{Levi-Form}^M(q)\big)
\Big)
\ \ \ \ \ \ \ \ \ \ \ \ \
{\scriptstyle{(\forall\,q\,\in\,M\,\cap\,{\sf U}_p)}}.
\]

Next, at one point $q \in M \cap {\sf U}_p$, consider
a {\em nonzero} vector:
\[
\mathcal{K}_q
=
\kappa_{1q}\,\mathcal{L}_1\big\vert_q
+
\kappa_{2q}\,\mathcal{L}_2\big\vert_q
\] 
which is assumed to belong to the Levi kernel, namely:
\[
0
=
\left(\!
\begin{array}{cc}
\ell_{11}(q) & \ell_{12}(q)
\\
\ell_{21}(q) & \ell_{22}(q)
\end{array}
\!\right)
\left(\!
\begin{array}{c}
\kappa_{1q}
\\
\kappa_{2q}
\end{array}
\!\right),
\]
that is to say:
\[
\aligned
0
&
=
\overbrace{\ell_{11}(q)}^{\neq\,0}\,\kappa_{1q}
+
\ell_{12}(q)\,\kappa_{2q},
\\
0
&
=
\ell_{21}(q)\,\kappa_{1q}
+
\ell_{22}(q)\,\kappa_{2q}.
\endaligned
\]
Thanks to the nowhere vanishing of $\ell_{11}$, 
one can solve the first line:
\[
\kappa_{1q}
=
-\,\frac{\ell_{12}(q)}{\ell_{11}(q)}\,
\kappa_{2q},
\]
while the second line is automatically satisfied
(mental exercise) thanks to the zeroness of the
Levi determinant.

One therefore gets that at $q \in M \cap {\sf U}_p$
arbitrary, the Levi kernel 
is generated, as a $1$-dimensional $\C$-vector
space, precisely by:
\[
\mathcal{K}_q
:=
-\,\frac{\ell_{12}(q)}{\ell_{11}(q)}\,
\mathcal{L}_1\big\vert_q
+
\mathcal{L}_2\big\vert_q.
\]

Final observation: here, because the coefficient-function:
\[
-\,\frac{\ell_{12}(q)}{\ell_{11}(q)}\,
\]
is of $\mathcal{ C}^{ \kappa - 2}$ smoothness, or $\mathcal{ C}^\infty$,
or $\mathcal{ C}^\omega$ on $M \cap {\sf U}_p$, since
$\ell_{ 11}$ is nowhere vanishing, then the collection of
all complex lines $\C \cdot \mathcal{ K}_q$ organizes coherently as a 
certain true line $\C$-subbundle $K^{1, 0} M
\subset T^{1, 0}M$. 

Of course, when one passes from one open set ${\sf U}_p$ to a nearby 
open set ${\sf U}_{ p'}$ with ${\sf U}_{ p'} \cap {\sf U}_p \neq
\emptyset$, the two definitions match up in the intersection,
because the Levi kernel exists independently of the choice
of local coordinates, and independently of the choice
of a local frame for $T^{1, 0}M$, as is
already known.
\endproof

Now, come back temporarily to hypersurfaces of any dimension $2 n + 1$.

\medskip\noindent{\bf Lemma.}
{\em On a connected hypersurface:}
\[
M^{2n+1}
\subset
\C^{n+1},
\]
{\em which is $\mathcal{ C}^\kappa$ $(\kappa \geqslant 2)$, 
or $\mathcal{ C}^\infty$, or $\mathcal{ C}^\omega$,
if the kernel of the Levi form is of
constant rank equal to a certain integer $e$
with:}
\[
0
\leqslant
e
\leqslant
n,
\]
{\em then the union of kernels gathers coherently to constitute
a certain complex vector subbundle:}
\[
K^{1,0}M
\subset
T^{1,0}M
\]
{\em of rank $e$ which in addition satisfies the three
involutiveness conditions:}
\[
\aligned
\big[K^{1,0}M,\,K^{1,0}M\big]
&
\,\subset\,
K^{1,0}M,
\\
\big[K^{0,1}M,\,K^{0,1}M\big]
&
\,\subset\,
K^{0,1}M,
\\
\big[K^{1,0}M,\,K^{0,1}M\big]
&
\,\subset\,
K^{1,0}M
\oplus
K^{0,1}M.
\endaligned
\]

\proof
Taking for granted (exercise) that the proof of the previous
Assertion can be elementarily generalized to yield
that $K^{ 1, 0}M$ is a rank $e$ vector subbundle
of $T^{1, 0}M$, it remains
to check the stated involutiveness conditions.

Recall how one expresses the hypotheses that a local vector field section:
\[
\mathcal{K}
\]
of $K^{1, 0}M$ belongs to the Levi-kernel 
at every point. In terms of any frame:
\[
\big\{
\mathcal{L}_1,\dots,\mathcal{L}_n
\big\}
\]
for $T^{1, 0} M$, and in terms of any local differential $1$-form:
\[
\rho_0\colon\ \ \
TM
\,\longrightarrow\,
\C
\]
satisfying:
\[
\big\{\rho_0=0\big\}
=
T^{1,0}M
\oplus
T^{0,1}M,
\]
one must have by hypothesis:
\[
0
=
\rho_0\big(\big[\mathcal{K},\,\overline{\mathcal{L}}_1\big]\big)
=\cdots=
\rho_0\big(\big[\mathcal{K},\,\overline{\mathcal{L}}_n\big]\big).
\]

Since the goal is to prove that:
\[
\aligned
\big[K^{1,0}M,\,K^{1,0}M\big]
&
\,\subset\,
K^{1,0}M,
\\
\big[K^{0,1}M,\,K^{0,1}M\big]
&
\,\subset\,
K^{0,1}M,
\\
\big[K^{1,0}M,\,K^{0,1}M\big]
&
\,\subset\,
K^{1,0}M
\oplus
K^{0,1}M,
\endaligned
\]
the second condition being the conjugate of the first,
take two local sections:
\[
\mathcal{K}_1
\ \ \ \ \ \ \ \ \ \ \ \ \
\text{\rm and}
\ \ \ \ \ \ \ \ \ \ \ \ \
\mathcal{K}_2
\]
which both satisfy such a hypothesis:
\[
\aligned
\big[\mathcal{K}_1,\,\overline{\mathcal{L}}_j\big]
&
=
\sum_{k=1}^n\,\function\cdot
\mathcal{L}_k
+
\sum_{l=1}^n\,\function\cdot\overline{\mathcal{L}}_l,
\\
\big[\mathcal{K}_2,\,\overline{\mathcal{L}}_j\big]
&
=
\sum_{k=1}^n\,\function\cdot
\mathcal{L}_k
+
\sum_{l=1}^n\,\function\cdot\overline{\mathcal{L}}_l.
\endaligned
\]
Simultaneously, reminding the automatic involutiveness:
\[
\big[T^{1,0}M,\,T^{1,0}M\big]
\,\subset\,
T^{1,0}M
\,\,\,\Longrightarrow\,\,\,
\big[K^{1,0}M,\,T^{1,0}M\big]
\,\subset\,
T^{1,0}M,
\]
one also has for free:
\[
\aligned
\big[\mathcal{K}_1,\,\mathcal{L}_j\big]
&
=
\sum_{k=1}^n\,\function\cdot\mathcal{L}_k,
\\
\big[\mathcal{K}_2,\,\mathcal{L}_j\big]
&
=
\sum_{k=1}^n\,\function\cdot\mathcal{L}_k.
\endaligned
\]

It is now time to examine whether the bracket:
\[
\big[\mathcal{K}_1,\,\mathcal{K}_2\big]
\]
still belongs to the kernel of the Levi form, namely to compute:
\[
\rho_0\big(
\big[\big[\mathcal{K}_1,\mathcal{K}_2\big],\,
\overline{\mathcal{L}}_j\big]
\big)
\overset{?}{\,=\,}
0
\ \ \ \ \ \ \ \ \ \ \ \ \
{\scriptstyle{(j\,=\,1\,\cdots\,n)}}.
\]
But the Jacobi identity:
\[
\big[\big[\mathcal{K}_1,\mathcal{K}_2\big],\,
\overline{\mathcal{L}}_j\big]
=
-\,
\big[\big[\overline{\mathcal{L}}_j,\mathcal{K}_1\big],\,
\mathcal{K}_2\big]
-
\big[\big[\mathcal{K}_2,\overline{\mathcal{L}}_j\big],\,\mathcal{K}_1
\big]
\]
yields thanks to two applications of the hypothesis:
\[
\footnotesize
\aligned
\rho_0\Big(
\big[\big[\mathcal{K}_1,\mathcal{K}_2\big],\,
\overline{\mathcal{L}}_j\big]\Big)
&
=
-\,
\rho_0\Big(
\big[\big[\overline{\mathcal{L}}_j,\mathcal{K}_1\big],\,
\mathcal{K}_2\big]
\Big)
-
\rho_0
\Big(
\big[\big[\mathcal{K}_2,\overline{\mathcal{L}}_j\big],\,\mathcal{K}_1
\big]
\Big)
\\
&
=
-\,\rho_0\Big(
\big[
\function_\bullet\,\mathcal{L}_\bullet
+
\function_\bullet\,\overline{\mathcal{L}}_\bullet,\,
\mathcal{K}_2
\big]
\Big)
-
\\
&
\ \ \ \ \ \ \ \ \ \ \ \ \ \ \ \
-\,\rho_0\Big(
\big[
\function_\bullet\,\mathcal{L}_\bullet
+
\function_\bullet\,\overline{\mathcal{L}}_\bullet,\,
\mathcal{K}_1
\big]
\Big)
\\
&
=
\rho_0
\Big(
\function_\bullet\,\mathcal{L}_\bullet
+
\function_\bullet\,\overline{\mathcal{L}}_\bullet
\Big)
\\
&
=
0,
\endaligned
\]
as desired.

The second involutiveness condition is proved similarly (or by conjugating).

To prove the third condition, one has to compute:
\[
\big[\mathcal{K}_1,\,\overline{\mathcal{K}}_2\big].
\]
Since this bracket is at least a local section of:
\[
T^{1,0}M\oplus T^{0,1}M,
\]
one decomposes it as:
\[
\big[\mathcal{K}_1,\,\overline{\mathcal{K}}_2\big]
\,=\,
\mathcal{M}
+
\overline{\mathcal{N}},
\]
where $\mathcal{ M}$ and $\mathcal{ N}$ are local sections
of $T^{1, 0}M$.
The goal is to prove that both are local sections
of the Levi kernel subbundle $K^{1, 0} M \subset T^{1, 0} M$.

At first, one abbreviates what precedes as:
\[
\aligned
\big[\mathcal{K}_1,\mathcal{L}_j\big]
&
=
\mathcal{R}_{1j},
\\
\big[\mathcal{K}_2,\mathcal{L}_j\big]
&
=
\mathcal{R}_{2j},
\\
\big[\mathcal{K}_1,\overline{\mathcal{L}}_j\big]
&
=
\mathcal{S}_{1j}
+
\overline{\mathcal{T}}_{1j},
\\
\big[\mathcal{K}_2,\overline{\mathcal{L}}_j\big]
&
=
\mathcal{S}_{2j}
+
\overline{\mathcal{T}}_{2j}.
\endaligned
\]

Then the Jacobi identity yields on one hand:
\[
\aligned
\Big[\big[\mathcal{K}_1,\overline{\mathcal{K}}_2\big],\,
\overline{\mathcal{L}}_j\Big]
&
=
-\,\Big[
\big[\overline{\mathcal{L}}_j,\mathcal{K}_1\big],\,
\overline{\mathcal{K}}_2
\Big]
-
\Big[
\big[\overline{\mathcal{K}}_2,\overline{\mathcal{L}}_j\big],\,
\mathcal{K}_1
\Big]
\\
&
=
\big[
\mathcal{S}_{1j}+\overline{\mathcal{T}}_{1j},\,\,
\overline{\mathcal{K}}_2
\big]
+
\big[\overline{\mathcal{R}}_{2j},\,\mathcal{K}_1\big]
\\
&
=
\underbrace{\big[
\mathcal{S}_{1j},\,\overline{\mathcal{K}}_2
\big]}_{
{\sf section of}
\atop
T^{1,0}M\oplus T^{0,1}M}
+
\underbrace{
\big[\overline{\mathcal{T}}_{1j},\,\overline{\mathcal{K}}_2\big]}_{
{\sf section of}
\atop
T^{0,1}M}
+
\underbrace{
\big[
\overline{\mathcal{R}}_{2j},\,
\mathcal{K}_1
\big]}_{
{\sf section of}
\atop
T^{1,0}M\oplus T^{0,1}M}
\\
&
=
\mathmotsf{vector field section of}\,\,
T^{1,0}M\oplus T^{0,1}M,
\endaligned
\]
while on another hand the same length $3$ bracket has value:
\[
\aligned
\Big[\big[\mathcal{K}_1,\overline{\mathcal{K}}_2\big],\,
\overline{\mathcal{L}}_j\Big]
&
=
\big[\mathcal{M}+\overline{\mathcal{N}},\,\,\overline{\mathcal{L}}_j\big]
\\
&
=
\big[\mathcal{M},\,\overline{\mathcal{L}}_j\big]
+
\underbrace{
\big[\overline{\mathcal{N}},\,\overline{\mathcal{L}}_j\big]}_{
{\sf section of}
\atop
T^{0,1}M},
\endaligned
\]
and then a final subtraction provides:
\[
\big[\mathcal{M},\,\overline{\mathcal{L}}_j\big]\,\,
\mathmotsf{is a local section of}\,\,
T^{1,0}M\oplus T^{0,1}M,
\]
{\em which means precisely that:}
\[
\mathcal{M}\,\,
\mathmotsf{is a section of the Levi kernel subbundle}\,\,
K^{1,0}M.
\]

One verifies (exercise) that a similar reasoning starting from:
\[
\Big[
\big[\mathcal{K}_1,\overline{\mathcal{K}}_2\big],\,
\mathcal{L}_j
\Big]
\]
proves that:
\[
\big[\overline{\mathcal{N}},\,\mathcal{L}_j\big]\,\,
\mathmotsf{is a local section of}\,\,
T^{1,0}M\oplus T^{0,1}M,
\]
which means, after conjugation, that:
\[
\mathcal{N}\,\,
\mathmotsf{is a section of the Levi kernel subbundle}\,\,
K^{1,0}M,
\]
and this concludes the proof.
\endproof

Introduce the {\em real} subbundle:
\[
\aligned
K^cM
:=
&\,
\Re\,K^{1,0}M
\\
\subset
&\,
\Re\,T^{1,0}M
=
T^cM.
\endaligned
\]

\medskip\noindent{\bf Corollary.}
{\em One has the involutiveness condition:}
\[
\big[K^cM,K^cM\big]
\subset
K^cM.
\]

\proof
Indeed, two general local sections of $K^cM$ write:
\[
\aligned
K_1
&
=
\mathcal{K}_1
+
\overline{\mathcal{K}}_1,
\\
K_2
&
=
\mathcal{K}_2
+
\overline{\mathcal{K}}_2,
\endaligned
\]
whence:
\[
\aligned
\big[K_1,K_2\big]
&
=
\big[
\mathcal{K}_1+\overline{\mathcal{K}}_1,\,\,
\mathcal{K}_2+\overline{\mathcal{K}}_2
\big]
\\
&
=
\underbrace{\big[\mathcal{K}_1,\,\mathcal{K}_2\big]}_{
{\sf section}\,\,\mathcal{K}_3
\atop
{\sf of}\,\,K^{1,0}M}
+
\underbrace{
\big[\overline{\mathcal{K}}_1,\,\overline{\mathcal{K}}_2\big]}_{
\overline{\mathcal{K}}_3}
+
\underbrace{
\big[\mathcal{K}_1,\,\overline{\mathcal{K}}_2\big]}_{
{\sf section}\,\,\mathcal{K}_4+\overline{\mathcal{K}}_5
\atop
{\sf of}\,\,K^{1,0}M\oplus K^{0,1}M}
+
\underbrace{
\big[\overline{\mathcal{K}}_1,\,\mathcal{K}_2\big]}_{
\overline{\mathcal{K}}_4+\mathcal{K}_5}
\\
&
=
\underbrace{
\mathcal{K}_3+\overline{\mathcal{K}}_3}_{
{\sf section}
\atop
{\sf of}\,\,K^cM}
+
\underbrace{
\mathcal{K}_4+\overline{\mathcal{K}}_4}_{
{\sf section}
\atop
{\sf of}\,\,K^cM}
+
\underbrace{
\mathcal{K}_5+\overline{\mathcal{K}}_5}_{
{\sf section}
\atop
{\sf of}\,\,K^cM},
\endaligned
\]
which concludes.
\endproof

\noindent{\bf Explicit expression of the function $k$ for $M^5 \subset
\C^3$.}
Now, come back to the case under study of a hypersurface:
\[
M^5
\subset
\C^3
\]
whose Levi-kernel bundle:
\[
K^{1,0}M
\subset
T^{1,0}M
\]
if of rank:
\[
{\bf 1}
=
\rank_\C\big(K^{1,0}M\big)
\,<\,
{\bf 2}
=
\rank_\C\big(T^{1,0}M\big).
\]
Take as usual a local differential $1$-form:
\[
\rho_0\colon\ \ \
TM
\,\longrightarrow\,
\R
\]
whose extension to $\C \otimes_\R TM$ satisfies:
\[
\big\{\rho_0=0\big\}
=
T^{1,0}M\oplus T^{0,1}M.
\]
Of course, such a $\rho_0$ is far from unique: it can be changed
by multiplying it by any nowhere vanishing real function.

Here, one must make more explicit all the data for $K^{1, 0}M$
in terms of a graphing function $\varphi$ for $M$:
\[
v
=
\varphi\big(x_1,x_2,y_1,y_2,u\big),
\]
in some local affine holomorphic coordinates:
\[
\big(z_1,z_2,w\big)
=
\big(x_1+\isqrt\,y_1,\,
x_2+\isqrt\,y_2,\,
u+\isqrt\,v\big)
\]
centered at some point $p \in M$.

Concerning smoothness, assume that:
\[
\varphi\in\mathcal{C}^\kappa\,\,
(\kappa\geqslant 3),
\ \ \ \ \ \ \ \ \ \ \ \ \ \
\text{\rm or}
\ \ \ \ \ \ \ \ \ \ \ \ \ \
\varphi\in\mathcal{C}^\infty,
\ \ \ \ \ \ \ \ \ \ \ \ \ \
\text{\rm or}
\ \ \ \ \ \ \ \ \ \ \ \ \ \
\varphi\in\mathcal{C}^\omega,
\]
becaues $\varphi$ will be differentiated thrice.

First of all, a natural local frame for $T^{1,0}M$:
\[
\big\{\mathcal{L}_1,\mathcal{L}_2\big\}
\]
is constituted by the two $(1,0)$ vector fields:
\[
\aligned
\mathcal{L}_1
&
=
\frac{\partial}{\partial z_1}
-
\frac{\varphi_{z_1}}{\isqrt+\varphi_u}\,
\frac{\partial}{\partial u},
\\
\mathcal{L}_2
&
=
\frac{\partial}{\partial z_2}
-
\frac{\varphi_{z_2}}{\isqrt+\varphi_u}\,
\frac{\partial}{\partial u},
\endaligned
\]
together with their conjugates:
\[
\aligned
\overline{\mathcal{L}}_1
&
=
\frac{\partial}{\partial\overline{z}_1}
-
\frac{\varphi_{\overline{z}_1}}{-\isqrt+\varphi_u}\,
\frac{\partial}{\partial u},
\\
\overline{\mathcal{L}}_2
&
=
\frac{\partial}{\partial\overline{z}_2}
-
\frac{\varphi_{\overline{z}_2}}{-\isqrt+\varphi_u}\,
\frac{\partial}{\partial u};
\endaligned
\]
recall indeed that $\varphi$ being real:
\[
\overline{\varphi(x_1,x_2,y_1,y_2,u)}
=
\varphi(x_1,x_2,y_1,y_2,u),
\]
one has:
\[
\overline{\varphi_{z_1}}
=
\varphi_{\overline{z}_1},
\ \ \ \ \ \ \ \ \ \ \ \ \ \ \ \ \ \ 
\overline{\varphi_{z_2}}
=
\varphi_{\overline{z}_2}.
\]

To begin with, compute the entries of:
\[
{\sf Levi-Matrix}_{\mathcal{L},\overline{\mathcal{L}}}^M(q)
=
\left(\!
\begin{array}{cc}
\rho_0\big(\isqrt\big[\mathcal{L}_1,\overline{\mathcal{L}}_1\big]\big)
&
\rho_0\big(\isqrt\big[\mathcal{L}_2,\overline{\mathcal{L}}_1\big]\big)
\\
\rho_0\big(\isqrt\big[\mathcal{L}_1,\overline{\mathcal{L}}_2\big]\big)
&
\rho_0\big(\isqrt\big[\mathcal{L}_2,\overline{\mathcal{L}}_2\big]\big)
\end{array}
\!\right)
(q).
\]
As an intermediation, if one abbreviates:
\[
\aligned
\mathcal{L}_1
&
=
\frac{\partial}{\partial z_1}
+
A_1\,\frac{\partial}{\partial u},
\ \ \ \ \ \ \ \ \ \ \ \ \ \ \ \ \ \ \ \ \ \ \ \ 
\overline{\mathcal{L}}_1
=
\frac{\partial}{\partial\overline{z}_1}
+
\overline{A_1}\,\frac{\partial}{\partial u},
\\
\mathcal{L}_2
&
=
\frac{\partial}{\partial z_2}
+
A_2\,\frac{\partial}{\partial u},
\ \ \ \ \ \ \ \ \ \ \ \ \ \ \ \ \ \ \ \ \ \ \ \ 
\overline{\mathcal{L}}_2
=
\frac{\partial}{\partial\overline{z}_2}
+
\overline{A_2}\,\frac{\partial}{\partial u},
\endaligned
\]
with of course:
\[
\aligned
A_1
:=
-\,
\frac{\varphi_{z_1}}{\isqrt+\varphi_u},
\\
A_2
:=
-\,
\frac{\varphi_{z_2}}{\isqrt+\varphi_u},
\endaligned
\]
when one computes the $4$ Lie brackets:
\[
\aligned
\isqrt\,\big[\mathcal{L}_1,\overline{\mathcal{L}}_1\big]
&
=
\isqrt
\bigg[
\frac{\partial}{\partial z_1}
+
A_1\,\frac{\partial}{\partial u},\,\,
\frac{\partial}{\partial\overline{z}_1}
+
\overline{A}_1\,\frac{\partial}{\partial u}
\bigg]
\\
&
=
\isqrt\,
\Big(\mathcal{L}_1\big(\overline{A_1}\big)
-
\overline{\mathcal{L}}_1\big(A_1\big)
\Big)\,
\frac{\partial}{\partial u},
\\
\isqrt\,\big[\mathcal{L}_2,\overline{\mathcal{L}}_1\big]
&
=
\isqrt\,
\Big(\mathcal{L}_2\big(\overline{A_1}\big)
-
\overline{\mathcal{L}}_1\big(A_2\big)
\Big)\,
\frac{\partial}{\partial u},
\\
\isqrt\,\big[\mathcal{L}_1,\overline{\mathcal{L}}_2\big]
&
=
\isqrt\,
\Big(\mathcal{L}_1\big(\overline{A_2}\big)
-
\overline{\mathcal{L}}_2\big(A_1\big)
\Big)\,
\frac{\partial}{\partial u},
\\
\isqrt\,\big[\mathcal{L}_2,\overline{\mathcal{L}}_2\big]
&
=
\isqrt\,
\Big(\mathcal{L}_2\big(\overline{A_2}\big)
-
\overline{\mathcal{L}}_2\big(A_2\big)
\Big)\,
\frac{\partial}{\partial u},
\endaligned
\]
$8$ functions appear:
\[
\aligned
&
\mathcal{L}_1\big(\overline{A_1}\big),
\ \ \ \ \ \ \ \ \ \ \ \ \ \ \ \ \ \ \ \ \ \ \ \ 
\overline{\mathcal{L}}_1\big(A_1\big),
\\
&
\mathcal{L}_2\big(\overline{A_1}\big),
\ \ \ \ \ \ \ \ \ \ \ \ \ \ \ \ \ \ \ \ \ \ \ \ 
\overline{\mathcal{L}}_1\big(A_2\big),
\\
&
\mathcal{L}_1\big(\overline{A_2}\big),
\ \ \ \ \ \ \ \ \ \ \ \ \ \ \ \ \ \ \ \ \ \ \ \ 
\overline{\mathcal{L}}_2\big(A_1\big),
\\
&
\mathcal{L}_2\big(\overline{A_2}\big),
\ \ \ \ \ \ \ \ \ \ \ \ \ \ \ \ \ \ \ \ \ \ \ \ 
\overline{\mathcal{L}}_2\big(A_2\big),
\endaligned
\]
that one should express in terms of $\varphi$.

Here is detailed computation for the first one:
\[
\aligned
\mathcal{L}_1\big(\overline{A_1}\big)
&
=
\bigg(
\frac{\partial}{\partial z_1}
-
\frac{\varphi_{z_1}}{\isqrt+\varphi_u}\,
\frac{\partial}{\partial u}
\bigg)
\bigg[
-\,\frac{\varphi_{\overline{z}_1}}{-\isqrt+\varphi_u}
\bigg]
\\
&
=
-\,\frac{\varphi_{z_1\overline{z}_1}}{-\isqrt+\varphi_u}
+
\frac{\varphi_{\overline{z}_1}\,\varphi_{z_1u}}{
(-\isqrt+\varphi_u)^2}
-
\frac{\varphi_{z_1}}{\isqrt+\varphi_u}\,
\bigg[
-\,\frac{\varphi_{\overline{z}_1u}}{-\isqrt+\varphi_u}
+
\frac{\varphi_{\overline{z}_1}\,\varphi_{uu}}{
(-\isqrt+\varphi_u)^2}
\bigg]
\\
&
=
\frac{-\,\varphi_{z_1\overline{z}_1}(1+\varphi_u^2)
+
\varphi_{\overline{z}_1}\varphi_{z_1u}(\isqrt+\varphi_u)
+
\varphi_{z_1}\varphi_{\overline{z}_1u}(-\isqrt+\varphi_u)
-
\varphi_{z_1}\varphi_{\overline{z}_1}\varphi_{uu}}{
(\isqrt+\varphi_u)\,(-\isqrt+\varphi_u)^2}.
\endaligned
\]
The conjugate is:
\[
\overline{\mathcal{L}}_1\big(A_1\big)
=
\frac{-\,\varphi_{z_1\overline{z}_1}(1+\varphi_u^2)
+
\varphi_{z_1}\varphi_{\overline{z}_1u}(-\isqrt+\varphi_u)
+
\varphi_{\overline{z}_1}\varphi_{z_1u}(\isqrt+\varphi_u)
-
\varphi_{z_1}\varphi_{\overline{z}_1}\varphi_{uu}}{
(\isqrt+\varphi_u)^2\,(-\isqrt+\varphi_u)}.
\]

Similarly, one obtains:
\[
\aligned
\mathcal{L}_2\big(\overline{A_1}\big)
&
=
\frac{-\,\varphi_{z_2\overline{z}_1}(1+\varphi_u^2)
+
\varphi_{\overline{z}_1}\varphi_{z_2u}(\isqrt+\varphi_u)
+
\varphi_{z_2}\varphi_{\overline{z}_1u}(-\isqrt+\varphi_u)
-
\varphi_{z_2}\varphi_{\overline{z}_1}\varphi_{uu}}{
(\isqrt+\varphi_u)\,(-\isqrt+\varphi_u)^2},
\\
\mathcal{L}_1\big(\overline{A_2}\big)
&
=
\frac{-\,\varphi_{z_1\overline{z}_2}(1+\varphi_u^2)
+
\varphi_{\overline{z}_2}\varphi_{z_1u}(\isqrt+\varphi_u)
+
\varphi_{z_1}\varphi_{\overline{z}_2u}(-\isqrt+\varphi_u)
-
\varphi_{z_1}\varphi_{\overline{z}_2}\varphi_{uu}}{
(\isqrt+\varphi_u)\,(-\isqrt+\varphi_u)^2},
\\
\mathcal{L}_2\big(\overline{A_2}\big)
&
=
\frac{-\,\varphi_{z_2\overline{z}_2}(1+\varphi_u^2)
+
\varphi_{\overline{z}_2}\varphi_{z_2u}(\isqrt+\varphi_u)
+
\varphi_{z_2}\varphi_{\overline{z}_2u}(-\isqrt+\varphi_u)
-
\varphi_{z_2}\varphi_{\overline{z}_2}\varphi_{uu}}{
(\isqrt+\varphi_u)\,(-\isqrt+\varphi_u)^2},
\endaligned
\]
while the conjugates are:
\[
\aligned
\overline{\mathcal{L}}_2\big(A_1\big)
&
=
\frac{-\,\varphi_{z_1\overline{z}_2}(1+\varphi_u^2)
+
\varphi_{z_1}\varphi_{\overline{z}_2u}(-\isqrt+\varphi_u)
+
\varphi_{\overline{z}_2}\varphi_{z_1u}(\isqrt+\varphi_u)
-
\varphi_{z_1}\varphi_{\overline{z}_2}\varphi_{uu}}{
(\isqrt+\varphi_u)^2\,(-\isqrt+\varphi_u)},
\\
\overline{\mathcal{L}}_1\big(A_2\big)
&
=
\frac{-\,\varphi_{z_2\overline{z}_1}(1+\varphi_u^2)
+
\varphi_{z_2}\varphi_{\overline{z}_1u}(-\isqrt+\varphi_u)
+
\varphi_{\overline{z}_1}\varphi_{z_2u}(\isqrt+\varphi_u)
-
\varphi_{z_2}\varphi_{\overline{z}_1}\varphi_{uu}}{
(\isqrt+\varphi_u)^2\,(-\isqrt+\varphi_u)},
\\
\overline{\mathcal{L}}_2\big(A_2\big)
&
=
\frac{-\,\varphi_{z_2\overline{z}_2}(1+\varphi_u^2)
+
\varphi_{z_2}\varphi_{\overline{z}_2u}(-\isqrt+\varphi_u)
+
\varphi_{\overline{z}_2}\varphi_{z_2u}(\isqrt+\varphi_u)
-
\varphi_{z_2}\varphi_{\overline{z}_2}\varphi_{uu}}{
(\isqrt+\varphi_u)^2\,(-\isqrt+\varphi_u)}.
\endaligned
\]

While computing the entries of the Levi matrix, in the
subtractions after reduction 
a common denominator, a number of terms disappear.
For instance, in:
\[
\aligned
\isqrt\,\big[\mathcal{L}_1,\overline{\mathcal{L}}_1\big]
&
=
\isqrt\,
\Big(\mathcal{L}_1(\overline{A_1}\big)
-
\overline{\mathcal{L}}_1\big(A_1\big)
\Big)\,
\frac{\partial}{\partial u}
\\
&
=
\isqrt\,
\Big(
\overline{A_1}_{z_1}
+
A_1\,\overline{A_1}_u
-
{A_1}_{\overline{z_1}}
-
\overline{A_1}\,
{A_1}_u
\Big)\,
\frac{\partial}{\partial u}
\endaligned
\]
one obtains after simplifications:
\[
\aligned
\isqrt\,
\big(\mathcal{L}_1(\overline{A_1}\big)
-
\overline{\mathcal{L}}_1\big(A_1\big)\big)
&
=
\frac{1}{(\isqrt+\varphi_u)^2(-\isqrt+\varphi_u)^2}
\bigg\{
2\,\varphi_{z_1\overline{z}_1}
+
2\,\varphi_{z_1\overline{z}_1}\varphi_u\varphi_u
-
\\
&
\ \ \ \ \
-\,
2\isqrt\,\varphi_{\overline{z}_1}\varphi_{z_1u}
-
2\,\varphi_{\overline{z}_1}\varphi_{z_1u}\varphi_u
+
2\isqrt\,\varphi_{z_1}\varphi_{\overline{z}_1u}
+
\\
&
\ \ \ \ \ \
+
2\,\varphi_{z_1}\varphi_{\overline{z}_1}\varphi_{uu}
-
2\,\varphi_{z_1}\varphi_{\overline{z}_1u}\varphi_u
\bigg\}.
\endaligned
\]
Of course, with the natural choice of:
\[
\rho_0
:=
-\,A_1\,dz_1
-
A_2\,dz_2
-
\overline{A_1}\,d\overline{z}_1
-
\overline{A_2}\,d\overline{z}_2
+
du,
\]
one has:
\[
\aligned
\rho_0\big(\isqrt\big[\mathcal{L}_1,\overline{\mathcal{L}}_1\big]\big)
&
=
\mathmotsf{this last expression}
\\
&
=
\isqrt\,
\Big(\mathcal{L}_1(\overline{A_1}\big)
-
\overline{\mathcal{L}}_1\big(A_1\big)
\Big)
\\
&
=
\isqrt\,
\Big(
\overline{A_1}_{z_1}
+
A_1\,\overline{A_1}_u
-
{A_1}_{\overline{z_1}}
-
\overline{A_1}\,
{A_1}_u
\Big).
\endaligned
\]

Recall that, introducing:
\[
\mathcal{T}
:=
\isqrt\,\big[\mathcal{L}_1,\overline{\mathcal{L}}_1\big],
\]
one is currently working under the assumption 
that the Levi form is everywhere of rank $1$, so
that after a possible ${\sf GL}_2 (\C)$ change of
coordinates, one can assume that the $5$ fields:
\[
\big\{
\mathcal{L}_1,\mathcal{L}_2,\overline{\mathcal{L}}_1,
\overline{\mathcal{L}}_2,\mathcal{T}\big\}
\]
constitute a frame for $\C \otimes_\R TM$.

Next, one computes similarly the remaining three entries
of the Levi matrix:
\[
\aligned
\isqrt\,
\big(\mathcal{L}_2(\overline{A_1}\big)
-
\overline{\mathcal{L}}_1\big(A_2\big)\big)
&
=
\frac{1}{(\isqrt+\varphi_u)^2(-\isqrt+\varphi_u)^2}
\bigg\{
2\,\varphi_{z_2\overline{z}_1}
+
2\,\varphi_{z_2\overline{z}_1}\varphi_u\varphi_u
-
\\
&\ \ \ \ \
-\,
2\isqrt\,\varphi_{\overline{z}_1}\varphi_{z_2u}
-
2\,\varphi_{\overline{z}_1}\varphi_{z_2u}\varphi_u
+
2\isqrt\,\varphi_{z_2}\varphi_{\overline{z}_1u}
+
\\
&
\ \ \ \ \ 
+
2\,\varphi_{z_2}\varphi_{\overline{z}_1}\varphi_{uu}
-
2\,\varphi_{z_2}\varphi_{\overline{z}_1u}\varphi_u
\bigg\},
\endaligned
\]
\[
\aligned
\isqrt\,
\big(\mathcal{L}_1(\overline{A_2}\big)
-
\overline{\mathcal{L}}_2\big(A_1\big)\big)
&
=
\frac{1}{(\isqrt+\varphi_u)^2(-\isqrt+\varphi_u)^2}
\bigg\{
2\,\varphi_{z_1\overline{z}_2}
+
2\,\varphi_{z_1\overline{z}_2}\varphi_u\varphi_u
-
\\
&
\ \ \ \ \
-\,
2\isqrt\,\varphi_{\overline{z}_2}\varphi_{z_1u}
-
2\,\varphi_{\overline{z}_2}\varphi_{z_1u}\varphi_u
+
2\isqrt\,\varphi_{z_1}\varphi_{\overline{z}_2u}
+
\\
&
\ \ \ \ \
+
2\,\varphi_{z_1}\varphi_{\overline{z}_2}\varphi_{uu}
-
2\,\varphi_{z_1}\varphi_{\overline{z}_2u}\varphi_u
\bigg\},
\endaligned
\]
\[
\aligned
\isqrt\,
\big(\mathcal{L}_2(\overline{A_2}\big)
-
\overline{\mathcal{L}}_2\big(A_2\big)\big)
&
=
\frac{1}{(\isqrt+\varphi_u)^2(-\isqrt+\varphi_u)^2}
\bigg\{
2\,\varphi_{z_2\overline{z}_2}
+
2\,\varphi_{z_2\overline{z}_2}\varphi_u\varphi_u
-
\\
&
\ \ \ \ \
-\,
2\isqrt\,\varphi_{\overline{z}_2}\varphi_{z_2u}
-
2\,\varphi_{\overline{z}_2}\varphi_{z_2u}\varphi_u
+
2\isqrt\,\varphi_{z_2}\varphi_{\overline{z}_2u}
+
\\
&
\ \ \ \ \
+
2\,\varphi_{z_2}\varphi_{\overline{z}_2}\varphi_{uu}
-
2\,\varphi_{z_2}\varphi_{\overline{z}_2u}\varphi_u
\bigg\}.
\endaligned
\]

Thanks to all these expressions, the:
\[
\aligned
\mathmotsf{Levi-Determinant}
&
=
\det\,
\left(\!
\begin{array}{cc}
\rho_0\big(\isqrt\big[\mathcal{L}_1,\overline{\mathcal{L}}_1\big]\big)
&
\rho_0\big(\isqrt\big[\mathcal{L}_2,\overline{\mathcal{L}}_1\big]\big)
\\
\rho_0\big(\isqrt\big[\mathcal{L}_1,\overline{\mathcal{L}}_2\big]\big)
&
\rho_0\big(\isqrt\big[\mathcal{L}_2,\overline{\mathcal{L}}_2\big]\big)
\end{array}
\!\right)
\\
&
=
\det\,
\left(\!
\begin{array}{cc}
\isqrt\big(
\mathcal{L}_1\big(\overline{A_1}\big)
-
\overline{\mathcal{L}}_1\big(A_1\big)
\big)
&
\isqrt\big(
\mathcal{L}_2\big(\overline{A_1}\big)
-
\overline{\mathcal{L}}_1\big(A_2\big)
\big)
\\
\isqrt\big(
\mathcal{L}_1\big(\overline{A_2}\big)
-
\overline{\mathcal{L}}_2\big(A_1\big)
\big)
&
\isqrt\big(
\mathcal{L}_2\big(\overline{A_2}\big)
-
\overline{\mathcal{L}}_2\big(A_2\big)
\big)
\end{array}
\!\right),
\endaligned
\]
equal to:
\[
\aligned
&
=
\det\,
\left(\!
\begin{array}{cc}
\isqrt
\big(
\overline{A_1}_{z_1}
+
A_1\,\overline{A_1}_u
-
{A_1}_{\overline{z_1}}
-
\overline{A_1}\,
{A_1}_u
\big)
&
\isqrt
\big(
\overline{A_1}_{z_2}
+
A_2\,\overline{A_1}_u
-
{A_2}_{\overline{z_1}}
-
\overline{A_1}\,
{A_2}_u
\big)
\\
\isqrt
\big(
\overline{A_2}_{z_1}
+
A_1\,\overline{A_2}_u
-
{A_1}_{\overline{z_2}}
-
\overline{A_2}\,
{A_1}_u
\big)
&
\isqrt
\big(
\overline{A_2}_{z_2}
+
A_2\,\overline{A_2}_u
-
{A_2}_{\overline{z_2}}
-
\overline{A_2}\,
{A_2}_u
\big)
\end{array}
\!\right)
\endaligned
\]
can be computed in terms of $\varphi$, 
and one obtains after simplifications:
\[
\boxed{\,
\aligned
&
\mathmotsf{Levi-Determinant}
=
\frac{4}{(\isqrt+\varphi_u)^3(-\isqrt+\varphi_u)^3}
\bigg\{
\varphi_{z_2\overline{z}_2}\varphi_{z_1\overline{z}_1}
-
\varphi_{z_2\overline{z}_1}\varphi_{z_1\overline{z}_2}
+
\\
&
+
\varphi_{z_2\overline{z}_1}\varphi_{\overline{z}_2}\varphi_{z_1u}
\varphi_u
-
\varphi_{z_2\overline{z}_1}\varphi_{\overline{z}_2}\varphi_{z_1}
\varphi_{uu}
-
\varphi_{\overline{z}_1}\varphi_{z_2u}\varphi_{z_1}
\varphi_{\overline{z}_2u}
+
\varphi_{\overline{z}_1}\varphi_{z_2u}\varphi_u
\varphi_{z_1\overline{z}_2}
-
\\
&
-\,
\varphi_{z_2}\varphi_{\overline{z}_1u}\varphi_{\overline{z}_2}
\varphi_{z_1u}
-
\varphi_{z_2}\varphi_{\overline{z}_1}\varphi_{uu}
\varphi_{z_1\overline{z}_2}
+
\varphi_{z_2}\varphi_{\overline{z}_1u}
\varphi_u\varphi_{z_1\overline{z}_2}
-
\varphi_{z_2\overline{z}_2}\varphi_{\overline{z}_1}
\varphi_{z_1u}\varphi_u
+
\\
&
+
\varphi_{z_2\overline{z}_2}\varphi_{z_1}
\varphi_{\overline{z}_1}\varphi_{uu}
-
\varphi_{z_2\overline{z}_2}\varphi_{z_1}
\varphi_{\overline{z}_1u}\varphi_u
+
\varphi_{z_2\overline{z}_1}\varphi_{z_1}
\varphi_{\overline{z}_2u}\varphi_u
+
\varphi_{z_2}\varphi_{\overline{z}_2u}
\varphi_{\overline{z}_1}\varphi_{z_1u}
-
\\
&
-\,
\varphi_{z_2}\varphi_{\overline{z}_2u}
\varphi_{z_1\overline{z}_1}\varphi_u
+
\varphi_{\overline{z}_2}\varphi_{z_2u}
\varphi_{z_1}\varphi_{\overline{z}_1u}
-
\varphi_{\overline{z}_2}\varphi_{z_2u}
\varphi_u\varphi_{z_1\overline{z}_1}
+
\varphi_{\overline{z}_2}\varphi_{z_2}
\varphi_{uu}\varphi_{z_1\overline{z}_1}
+
\\
&
+
\isqrt
\big(
\varphi_{z_2\overline{z}_2}\varphi_{z_1}
\varphi_{\overline{z}_1u}
+
\varphi_{\overline{z}_1}\varphi_{z_2u}
\varphi_{z_1\overline{z}_2}
+
\varphi_{z_2\overline{z}_1}\varphi_{\overline{z}_2}
\varphi_{z_1u}
+
\varphi_{z_2}\varphi_{\overline{z}_2u}
\varphi_{z_1\overline{z}_1}
\big)
-
\\
&
-\,\isqrt\big(
\varphi_{\overline{z}_2}\varphi_{z_2u}
\varphi_{z_1\overline{z}_1}
+
\varphi_{z_2\overline{z}_1}\varphi_{z_1}
\varphi_{\overline{z}_2u}
+
\varphi_{z_2}\varphi_{\overline{z}_1u}
\varphi_{z_1\overline{z}_2}
+
\varphi_{z_2\overline{z}_2}\varphi_{\overline{z}_1}
\varphi_{z_1u}
\big)
-
\\
&
-\,
\varphi_{z_2\overline{z}_1}\varphi_{z_1\overline{z}_2}
\varphi_u\varphi_u
+
\varphi_{z_2\overline{z}_2}\varphi_{z_1\overline{z}_1}
\varphi_u\varphi_u
\bigg\}.\,\,
\endaligned}
\]
So, this Levi determinant is assumed to be identically zero:
\[
\aligned
0
&
\equiv
\varphi_{z_2\overline{z}_2}\varphi_{z_1\overline{z}_1}
-
\varphi_{z_2\overline{z}_1}\varphi_{z_1\overline{z}_2}
+
\\
&\ \ \ \ \
+
\varphi_{z_2\overline{z}_1}\varphi_{\overline{z}_2}\varphi_{z_1u}
\varphi_u
-
\varphi_{z_2\overline{z}_1}\varphi_{\overline{z}_2}\varphi_{z_1}
\varphi_{uu}
-
\varphi_{\overline{z}_1}\varphi_{z_2u}\varphi_{z_1}
\varphi_{\overline{z}_2u}
+
\varphi_{\overline{z}_1}\varphi_{z_2u}\varphi_u
\varphi_{z_1\overline{z}_2}
-
\\
&\ \ \ \ \
-\,
\varphi_{z_2}\varphi_{\overline{z}_1u}\varphi_{\overline{z}_2}
\varphi_{z_1u}
-
\varphi_{z_2}\varphi_{\overline{z}_1}\varphi_{uu}
\varphi_{z_1\overline{z}_2}
+
\varphi_{z_2}\varphi_{\overline{z}_1u}
\varphi_u\varphi_{z_1\overline{z}_2}
-
\varphi_{z_2\overline{z}_2}\varphi_{\overline{z}_1}
\varphi_{z_1u}\varphi_u
+
\\
&\ \ \ \ \
+
\varphi_{z_2\overline{z}_2}\varphi_{z_1}
\varphi_{\overline{z}_1}\varphi_{uu}
-
\varphi_{z_2\overline{z}_2}\varphi_{z_1}
\varphi_{\overline{z}_1u}\varphi_u
+
\varphi_{z_2\overline{z}_1}\varphi_{z_1}
\varphi_{\overline{z}_2u}\varphi_u
+
\varphi_{z_2}\varphi_{\overline{z}_2u}
\varphi_{\overline{z}_1}\varphi_{z_1u}
-
\\
&\ \ \ \ \
-\,
\varphi_{z_2}\varphi_{\overline{z}_2u}
\varphi_{z_1\overline{z}_1}\varphi_u
+
\varphi_{\overline{z}_2}\varphi_{z_2u}
\varphi_{z_1}\varphi_{\overline{z}_1u}
-
\varphi_{\overline{z}_2}\varphi_{z_2u}
\varphi_u\varphi_{z_1\overline{z}_1}
+
\varphi_{\overline{z}_2}\varphi_{z_2}
\varphi_{uu}\varphi_{z_1\overline{z}_1}
+
\\
&\ \ \ \ \
+
\isqrt
\big(
\varphi_{z_2\overline{z}_2}\varphi_{z_1}
\varphi_{\overline{z}_1u}
+
\varphi_{\overline{z}_1}\varphi_{z_2u}
\varphi_{z_1\overline{z}_2}
+
\varphi_{z_2\overline{z}_1}\varphi_{\overline{z}_2}
\varphi_{z_1u}
+
\varphi_{z_2}\varphi_{\overline{z}_2u}
\varphi_{z_1\overline{z}_1}
\big)
-
\\
&\ \ \ \ \
-\,\isqrt\big(
\varphi_{\overline{z}_2}\varphi_{z_2u}
\varphi_{z_1\overline{z}_1}
+
\varphi_{z_2\overline{z}_1}\varphi_{z_1}
\varphi_{\overline{z}_2u}
+
\varphi_{z_2}\varphi_{\overline{z}_1u}
\varphi_{z_1\overline{z}_2}
+
\varphi_{z_2\overline{z}_2}\varphi_{\overline{z}_1}
\varphi_{z_1u}
\big)
-
\\
&\ \ \ \ \
-\,
\varphi_{z_2\overline{z}_1}\varphi_{z_1\overline{z}_2}
\varphi_u\varphi_u
+
\varphi_{z_2\overline{z}_2}\varphi_{z_1\overline{z}_1}
\varphi_u\varphi_u.
\endaligned
\]

Now, remind that the local generator:
\[
\mathcal{K}
=
k\,\mathcal{L}_1
+
\mathcal{L}_2
\]
of the Levi-kernel bundle $K^{1, 0}M$ has as its coefficient-function:
\[
\aligned
k
&
=
-\,
\frac{
\rho_0\big(\isqrt\,
\big[\mathcal{L}_2,\overline{\mathcal{L}}_1\big]\big)}{
\rho_0\big(\isqrt\,
\big[\mathcal{L}_1,\overline{\mathcal{L}}_1\big]\big)}
\\
&
=
-\,
\frac{
\mathcal{L}_2\big(\overline{A_1}\big)
-
\overline{\mathcal{L}}_1\big(A_2\big)
}{
\mathcal{L}_1\big(\overline{A_1}\big)
-
\overline{\mathcal{L}}_1\big(A_1\big)
},
\endaligned
\]
namely:
\[
\!\!\!\!\!\!\!\!\!\!\!\!\!\!\!\!\!\!
\boxed{\,\,
\aligned
k
&
=
\frac{\varphi_{z_2\overline{z}_1}
+
\varphi_{z_2\overline{z}_1}\varphi_u\varphi_u
-
\isqrt\varphi_{\overline{z}_1}\varphi_{z_2u}
-
\varphi_{\overline{z}_1}\varphi_{z_2u}\varphi_u
+
\isqrt\varphi_{z_2}\varphi_{\overline{z}_1u}
+
\varphi_{z_2}\varphi_{\overline{z}_1}\varphi_{uu}
-
\varphi_{z_2}\varphi_{\overline{z}_1u}\varphi_u}{
-\,\varphi_{z_1\overline{z}_1}
-
\varphi_{z_1\overline{z}_1}\varphi_u\varphi_u
+
\isqrt\varphi_{\overline{z}_1}\varphi_{z_1u}
+
\varphi_{\overline{z}_1}\varphi_{z_1u}\varphi_u
-
\isqrt\varphi_{z_1}\varphi_{\overline{z}_1u}
-
\varphi_{z_1}\varphi_{\overline{z}_1}\varphi_{uu}
+
\varphi_{z_1}\varphi_{\overline{z}_1u}\varphi_u
},
\endaligned}
\]
and there is a:

\medskip\noindent{\bf Surprising Computational fact.}
{\em This function:}
\[
k
=
-\,
\frac{
\mathcal{L}_2\big(\overline{A_1}\big)
-
\overline{\mathcal{L}}_1\big(A_2\big)
}{
\mathcal{L}_1\big(\overline{A_1}\big)
-
\overline{\mathcal{L}}_1\big(A_1\big)
},
\]
{\em when expressed back in terms
of the graphing function for $M$:}
\[
\varphi
=
\varphi\big(x_1,x_2,y_1,y_2,u\big),
\]
{\em happens to be also equal to the other two quotients:}
\[
\aligned
k
&
=
-\,
\frac{\mathcal{L}_2\big(\overline{A_1}\big)}{
\mathcal{L}_1\big(\overline{A_1}\big)}
\\
&
=
-\,
\frac{-\,\overline{\mathcal{L}}_1\big(A_2\big)}{
-\,\overline{\mathcal{L}}_1\big(A_1\big)}.
\endaligned
\]

\proof
By coming back to the expressions of these
numerators:
\[
\mathcal{L}_2\big(\overline{A_1}\big),
\ \ \ \ \ \ \ \ \ \ \ \ \ \ \ \ \ \ \ \ \ \ \ \
-\,\overline{\mathcal{L}}_1\big(A_2\big),
\] 
and of these denominators:
\[
\mathcal{L}_1\big(\overline{A_1}\big),
\ \ \ \ \ \ \ \ \ \ \ \ \ \ \ \ \ \ \ \ \ \ \ \
\overline{\mathcal{L}}_1\big(A_1\big)
\]
provided explicitly above, the fact becomes visible (eyes exercise).
\endproof

One can provide a partial enlightement of why this
fact is true by sticking to the particular
so-called {\sl rigid case} in which the
graphing function $\varphi$ is {\em independent
of $u = \Re\, w$}, namely when the equation of
$M$ is under the form:
\[
v
=
\varphi\big(x_1,x_2,y_1,y_2\big).
\]

In this case, the computations greatly simplify. Indeed:
\[
\aligned
\mathcal{L}_1
&
=
\frac{\partial}{\partial z_1}
+
\underbrace{\isqrt\,\varphi_{z_1}}_{=\,A_1}\,
\frac{\partial}{\partial u},
\ \ \ \ \ \ \ \ \ \ \ \ \ \ \ \ \ \ \ \ \ \ \ \ 
\overline{\mathcal{L}}_1
=
\frac{\partial}{\partial\overline{z}_1}
+
\underbrace{-\isqrt\,\varphi_{\overline{z}_1}}_{
=\,
\overline{A_1}}\,
\frac{\partial}{\partial u},
\\
\mathcal{L}_2
&
=
\frac{\partial}{\partial z_2}
+
\underbrace{\isqrt\,\varphi_{z_2}}_{
=\,A_2}\,\frac{\partial}{\partial u},
\ \ \ \ \ \ \ \ \ \ \ \ \ \ \ \ \ \ \ \ \ \ \ \ 
\overline{\mathcal{L}}_2
=
\frac{\partial}{\partial\overline{z}_2}
+
\underbrace{
-\isqrt\,\varphi_{\overline{z}_2}}_{
=\,\overline{A_2}}\,
\frac{\partial}{\partial u},
\endaligned
\]
so that the Levi matrix becomes:
\[
\left(\!
\begin{array}{cc}
2\,\varphi_{z_1\overline{z}_1} 
&
2\,\varphi_{z_2\overline{z}_1}
\\
2\,\varphi_{z_1\overline{z}_2}
&
2\,\varphi_{z_2\overline{z}_2}
\end{array}
\!\right),
\]
whence:
\[
\aligned
k
&
=
-\,\frac{2\,\varphi_{z_2\overline{z}_1}}{
2\,\varphi_{z_1\overline{z}_1}}
\\
&
=
-\,\frac{\varphi_{z_2\overline{z}_1}}{
\varphi_{z_1\overline{z}_1}}
\\
&
=
-\,\frac{\mathcal{L}_2(\varphi_{\overline{z}_1})}{
\mathcal{L}_1(\varphi_{\overline{z}_1})}
\\
&
=
-\,
\frac{-\,\overline{\mathcal{L}}_1(\varphi_{z_2})}{
-\,\overline{\mathcal{L}}_1(\varphi_{z_1})}.
\endaligned
\]

\medskip\noindent{\bf Transfer through local biholomorphisms
and Freeman form.} With $p \in M^5 \subset \C^3$ and ${\sf U}_p \ni p$ open, 
when a local biholomorphism is given:
\[
h\colon\ \ \ 
{\sf U}_p
\overset{\sim}{\,\longrightarrow\,}
{\sf U}_{p'}'
\]
of ${\sf U}_p$ onto an image open set:
\[
{\sf U}_{p'}' 
\,:=\, 
h\big({\sf U}_p\big)
\subset
\C^{n+1}
\ \ \ \ \ \ \ \ \ \ \ \ \
{\scriptstyle{(p'\,=\,h(p))}},
\]
so that:
\[
M'
\subset
{\sf U}_{p'}'
\]
is also a hypersurface of ${\C'}^3$, it has already been proved above that
the rank of the Levi form of $M'$ is also equal to 
$1$ at every point $q' \in M'$, and that 
the Levi-kernel bundles transfer properly through $h$:
\[
\aligned
h_*(K^{1,0}M)
&
=
K^{1,0}M',
\\
h_*(K^{0,1}M)
&
=
K^{0,1}M'.
\endaligned
\]

In terms of two local vector field generators:
\[
\mathcal{K},
\ \ \ \ \ \ \ \ \ \ \ \ \ \ \ \ \ \ \ \ \ \ \ \ 
\mathcal{K}',
\]
for $K^{1,0}M$ and for $K^{1,0}M'$, this
means that:
\[
h_*(\mathcal{K})
=
c'\,\mathcal{K'},
\]
for some $\mathcal{ C}^{\kappa - 1}$ ($\kappa \geqslant 3$),
or $\mathcal{ C}^\infty$, or $\mathcal{ C}^\omega$ nowhere
vanishing function: 
\[
c'
\colon\ \ \
M'
\,\longrightarrow\,
\C\backslash\{0\}.
\]

\medskip

{\em It is therefore absolutely natural not to take 
$\{ \mathcal{ L}_1, \mathcal{ L}_2\}$ but:}
\[
\big\{
\mathcal{L}_1,\,\mathcal{K}
\big\}
\]
{\em as a frame for $T^{1,0}M$ and similarly:}
\[
\big\{
\mathcal{L}_1',\,\mathcal{K}'
\big\}
\]
{\em as a frame for $T^{ 1, 0} M'$.}

\medskip

Through a local bihololomorphism $h$, one also has:
\[
h_*\big(\mathcal{L}_1\big)
=
a'\,\mathcal{L}_1'
+
b'\,\mathcal{K}',
\]
for two certain functions:
\[
\aligned
a'
\colon\ \ \
M'
&
\,\longrightarrow\,
\C\backslash\{0\},
\\
b'
\colon\ \ \
M'
&
\,\longrightarrow\,
\C,
\endaligned
\]
with $a'$ (only) also vanishing nowhere, since linear independency
is preserved.

\medskip

Now, introduce as before a local differential $1$-form:
\[
\rho_0
\colon\ \ \ 
TM
\,\longrightarrow\,
\R,
\]
whose extension to $\C \otimes_\R TM$ satisfies:
\[
\big\{\rho_0=0\big\}
=
T^{1,0}M\cap T^{0,1}M,
\]
and do the same for the $M'$-side:
\[
\big\{\rho_0'=0\big\}
=
T^{1,0}M'\cap T^{0,1}M'.
\]
As was seen above, a possible, natural choice is
\[
\aligned
\rho_0
&
:=
-\,A_1\,dz_1
-
A_2\,dz_2
-
\overline{A_1}\,d\overline{z}_1
-
\overline{A_2}\,d\overline{z}_2
+
du,
\\
\rho_0'
&
:=
-\,A_1'\,dz_1'
-
A_2'\,dz_2'
-
\overline{A_1}'\,d\overline{z}_1'
-
\overline{A_2}'\,d\overline{z}_2'
+
du',
\endaligned
\]
if, as understood, the equation of $M'$ writes similarly:
\[
v'
=
\varphi'\big(x_1',x_2',y_1',y_2',u'\big),
\]
with:
\[
\aligned
\mathcal{L}_1'
&
=
\frac{\partial}{\partial z_1'}
\underbrace{-
\frac{\varphi_{z_1'}'}{\isqrt+\varphi_{u'}'}}_{
=:\,\,A_1'}\,
\frac{\partial}{\partial u'}
=
\frac{\partial}{\partial z_1'}
+
A_1'
\frac{\partial}{\partial u'},
\\
\mathcal{L}_2'
&
=
\frac{\partial}{\partial z_2'}
\underbrace{-
\frac{\varphi_{z_2'}'}{\isqrt+\varphi_{u'}'}}_{
=:\,\,A_2'}\,
\frac{\partial}{\partial u'}
=
\frac{\partial}{\partial z_2'}
+
A_2'
\frac{\partial}{\partial u'}.
\endaligned
\]

One therefore has:
\[
\aligned
0
&
=
\rho_0\big(\mathcal{L}_1\big)
=
\rho_0\big(\mathcal{K}\big)
=
\rho_0\big(\overline{\mathcal{L}}_1\big)
=
\rho_0\big(\overline{\mathcal{K}}\big),
\\
0
&
=
\rho_0'\big(\mathcal{L}_1'\big)
=
\rho_0'\big(\mathcal{K}'\big)
=
\rho_0'\big(\overline{\mathcal{L}}_1'\big)
=
\rho_0'\big(\overline{\mathcal{K}}'\big).
\endaligned
\]

Also, since the two Levi determinants:
\[
\aligned
0
&
\equiv
\det
\left(\!
\begin{array}{cc}
\rho_0\big(\isqrt\big[\mathcal{L}_1,\overline{\mathcal{L}}_1\big]\big)
&
\rho_0\big(\isqrt\big[\mathcal{L}_2,\overline{\mathcal{L}}_1\big]\big)
\\
\rho_0\big(\isqrt\big[\mathcal{L}_1,\overline{\mathcal{L}}_2\big]\big)
&
\rho_0\big(\isqrt\big[\mathcal{L}_2,\overline{\mathcal{L}}_2\big]\big)
\end{array}
\!\right),
\\
0
&
\equiv
\det
\left(\!
\begin{array}{cc}
\rho_0'\big(\isqrt\big[\mathcal{L}_1',\overline{\mathcal{L}}_1'\big]\big)
&
\rho_0'\big(\isqrt\big[\mathcal{L}_2',\overline{\mathcal{L}}_1'\big]\big)
\\
\rho_0'\big(\isqrt\big[\mathcal{L}_1',\overline{\mathcal{L}}_2'\big]\big)
&
\rho_0'\big(\isqrt\big[\mathcal{L}_2',\overline{\mathcal{L}}_2'\big]\big)
\end{array}
\!\right),
\endaligned
\]
vanish identically, one introduces two {\sl slant-functions}:
\[
\aligned
k
&
:=
-\,
\frac{
\mathcal{L}_2\big(\overline{A_1}\big)
-
\overline{\mathcal{L}}_1\big(A_2\big)
}{
\mathcal{L}_1\big(\overline{A_1}\big)
-
\overline{\mathcal{L}}_1\big(A_1\big)
},
\\
k'
&
:=
-\,
\frac{
\mathcal{L}_2'\big(\overline{A_1}'\big)
-
\overline{\mathcal{L}}_1'\big(A_2'\big)
}{
\mathcal{L}_1'\big(\overline{A_1}'\big)
-
\overline{\mathcal{L}}_1'\big(A_1'\big)
},
\endaligned
\]
in terms of which two local generators of $K^{1, 0}M$ and
of $K^{1,0}M'$ are:
\[
\aligned
\mathcal{K}
&
=
k\,\mathcal{L}_1
+
\mathcal{L}_2,
\\
\mathcal{K}'
&
=
k'\,\mathcal{L}_1'
+
\mathcal{L}_2'.
\endaligned
\]

Introduce also the two differential $1$-forms:
\[
\aligned
\kappa_0
&
:=
dz_1
-
k\,dz_2,
\\
\kappa_0'
&
:=
dz_1'
-
k'\,dz_2',
\endaligned
\]
that are local sections of $T^{*(1,0)}M$ and of $T^{*(1,0)}M'$
which visibly satisfy:
\[
\aligned
0
&
=
\kappa_0\big(\mathcal{K}\big)
=
\kappa_0\big(\overline{\mathcal{L}}_1\big)
=
\kappa_0\big(\overline{\mathcal{K}}\big),
\\
0
&
=
\kappa_0'\big(\mathcal{K}'\big)
=
\kappa_0'\big(\overline{\mathcal{L}}_1'\big)
=
\kappa_0'\big(\overline{\mathcal{K}}'\big).
\endaligned
\]
Lastly, introduce the two $(1, 0)$-differential $1$-forms:
\[
\aligned
\zeta_0
&
:=
dz_1,
\\
\zeta_0'
&
:=
dz_1',
\endaligned
\]
which complete a coframe:
\[
\aligned
&
\big\{\zeta_0,\,\kappa_0\big\},
\\
&
\big\{\zeta_0',\,\kappa_0'\big\},
\endaligned
\]
for $T^{*(1,0)}M$ and for $T^{*(1,0)}M'$. By conjugating:
\[
\aligned
&
\big\{\overline{\zeta}_0,\,\overline{\kappa}_0\big\},
\\
&
\big\{\overline{\zeta}_0',\,\overline{\kappa}_0'\big\},
\endaligned
\]
constitute a coframe for $T^{*(0,1)}M$ and for 
$T^{*(0,1)}M'$.

\medskip

Now, since $h$ is a local biholomorphism, it transfers $(1, 0)$-differential
$1$-forms to $(1, 0)$-differential $1$-forms, 
in the sense that:
\[
\aligned
h_*\big(
\big\{
0=\rho_0
\big\}
\big)
&
\,=\,
\big\{
0=\rho_0'
\big\},
\\
h_*\big(
\big\{
0=\rho_0=\overline{\kappa}_0=\overline{\zeta}_0
\big\}
\big)
&
\,=\,
\big\{
0=\rho_0'=\overline{\kappa}_0'=\overline{\zeta}_0'
\big\},
\\
h_*\big(
\big\{
0=\rho_0=\kappa_0=\zeta_0
\big\}
\big)
&
\,=\,
\big\{
0=\rho_0'=\kappa_0'=\zeta_0'
\big\}.
\endaligned
\]
Dropping now any symbolic mention of $h_*$, one therefore has:
\[
\aligned
\rho_0
&
=
d'\,\rho_0',
\\
\kappa_0
&
=
e'\,\rho_0'
+
f'\,\kappa_0'
+
g'\,\zeta_0',
\endaligned
\]
with four certain functions:
\[
\aligned
d'\colon\ \ \ 
M'
&
\,\longrightarrow\,
\C,
\\
e'\colon\ \ \
M'
&
\,\longrightarrow\,
\C,
\ \ \ \ \ \ \ \ \ \ \ \ \ \ \ \ \ \ \ \ 
f'\colon\ \ \
M'
&
\,\longrightarrow\,
\C,
\ \ \ \ \ \ \ \ \ \ \ \ \ \ \ \ \ \ \ \ 
g'\colon\ \ \
M'
&
\,\longrightarrow\,
\C.
\endaligned
\]

But since furthermore:
\[
h_*
\big(
K^{1,0}M
\big)
=
K^{1,0}M'
\]
a condition which reads in terms of the coframe:
\[
h_*\big(
\big\{
0=\rho_0=\kappa_0
\big\}
\big)
\,=\,
\big\{
0=\rho_0'=\kappa_0'
\big\},
\]
one must have:
\[
g'
=
0,
\]
and hence:
\[
\aligned
\rho_0
&
=
d'\,\rho_0',
\\
\kappa_0
&
=
e'\,\rho_0'
+
f'\,\kappa_0',
\endaligned
\]
with nowhere vanishing:
\[
f'\colon\ \ \
M'
\,\longrightarrow\,
\C\backslash\{0\},
\]
to preserve independency.

\medskip

To motivate the concept of {\sl Freeman form}, pick 
two functions:
\[
\mu\colon\ \ \
M
\,\longrightarrow\,
\C,
\ \ \ \ \ \ \ \ \ \ \ \ \ \ \ \ \ \ \ \ 
\nu\colon\ \ \
M
\,\longrightarrow\,
\C,
\]
and compute:
\[
\aligned
\kappa_0\Big(
\big[\mu\,\mathcal{K},\,\,
\overline{\nu}\,\overline{\mathcal{L}}_1\big]
\Big)
&
=
\kappa_0
\Big(
\mu\overline{\nu}\,\big[\mathcal{K},\overline{\mathcal{L}}_1\big]
+
\mu\,\mathcal{K}\big(\overline{\nu}\big)
\cdot
\overline{\mathcal{L}}_1
-
\overline{\nu}\,\overline{\mathcal{L}}_1\big(\mu\big)
\cdot
\mathcal{K}
\Big)
\\
&
=
\kappa_0
\Big(
\mu\overline{\nu}\,\big[\mathcal{K},\overline{\mathcal{L}}_1\big]
\Big)
+
\mu\,\mathcal{K}\big(\overline{\nu}\big)
\cdot
\zero{\kappa_0
\big(
\overline{\mathcal{L}}_1
\big)}
-
\overline{\nu}\,\overline{\mathcal{L}}_1\big(\mu\big)
\cdot
\zero{
\kappa_0
\big(
\mathcal{K}
\big)}
\\
&
=
\kappa_0
\Big(
\mu\overline{\nu}\,\big[\mathcal{K},\overline{\mathcal{L}}_1\big]
\Big).
\endaligned
\]

The introduction of a natural generalization
of the Levi form, the so-called {\sl Freeman form},
will be justified by the:

\medskip\noindent{\bf Claim.}
{\em Up to a nonzero function-factor, the result is the 
same in the right-hand side hypersurface $M'$:}
\[
\kappa_0\Big(
\big[\mu\,\mathcal{K},\,\,
\overline{\nu}\,\overline{\mathcal{L}}_1\big]
\Big)
=
\mathmotsf{nonzero function}\cdot
\kappa_0'\Big(
\big[\mu'\,\mathcal{K}',\,\,
\overline{\nu}'\,\overline{\mathcal{L}}_1'\big]
\Big).
\]

\medskip

Indeed, setting:
\[
\aligned
\mu'
&
:=
\mu\circ h^{-1},
\\
\nu'
&
:=
\nu\circ h^{-1},
\endaligned
\]
one starts to compute how this expression transfers:
\[
\aligned
\kappa_0
\Big(
\big[\mu\,\mathcal{K},\,\,
\overline{\nu}\,\overline{\mathcal{L}}_1\big]
\Big)
&
=
\big(e'\rho_0'+f'\kappa_0'\big)
\Big(
\big[\mu'c'\,\mathcal{K}',\,\,
\overline{\nu}'\overline{a}'\,\overline{\mathcal{L}}_1'
+
\overline{\nu}'\overline{b}'\,\overline{\mathcal{K}}'
\big]
\Big)
\\
&
=
\big(e'\rho_0'+f'\kappa_0'\big)
\Big(
\mu'\overline{\nu}'c'\overline{a}'\,
\big[\mathcal{K}',\overline{\mathcal{L}}_1'\big]
+
\mu'\overline{\nu}'c'\overline{b}'\,
\big[\mathcal{K}',\overline{\mathcal{K}}'\big]
+
\\
&
\ \ \ \ \ \ \ \ \ \ \ \ \ \ \ \ \ \ \ \ \ \ \ \ \ \ \ \ \ \
+
\mu'c'\mathcal{K}'\big(\overline{\nu}'\overline{a}'\big)
\cdot
\overline{\mathcal{L}}_1'
+
\mu'c'\mathcal{K}'\big(\overline{\nu}'\overline{b}'\big)
\cdot
\overline{\mathcal{K}}'
-
\\
&
\ \ \ \ \ \ \ \ \ \ \ \ \ \ \ \ \ \ \ \ \ \ \ \ \ \ \ \ \ \
-
\overline{\nu'}\overline{a}'\overline{\mathcal{L}}_1'
\big(\mu'c'\big)
\cdot
\mathcal{K}'
-
\overline{\nu}'\overline{b}'\overline{\mathcal{K}}'
\big(\mu'c'\big)
\cdot
\mathcal{K}'
\Big),
\endaligned
\]
and further, by distributing the actions of the
two differential $1$-forms $e' \rho_0'$ and 
$f'\kappa_0'$, starting with the second one:
\[
\aligned
\kappa_0
\Big(
\big[\mu\,\mathcal{K},\,\,
\overline{\nu}\,\overline{\mathcal{L}}_1\big]
\Big)
&
=
f'c'\overline{a}'\,
\mu'\overline{\nu}'\,
\kappa_0'
\big(\big[\mathcal{K}',\overline{\mathcal{L}}_1'\big]\big)
+
f'\mu'\overline{\nu}'c'\overline{b}'\,
\kappa_0'\big(\big[\mathcal{K}',\overline{\mathcal{K}}'\big]\big)
+
\\
&
\ \ \ \ \ 
+
f'\mu'c'\mathcal{K}'\big(\overline{\nu}'\overline{a}'\big)\,
\zero{\kappa_0'\big(\overline{\mathcal{L}}_1'\big)}
+
f'\mu'c'\mathcal{K}'\big(\overline{\nu}'\overline{b}'\big)\,
\zero{\kappa_0'\big(\overline{\mathcal{K}}'\big)}
-
\\
&
\ \ \ \ \
-
f'\overline{\nu}'\overline{a}'\overline{\mathcal{L}}_1'
\big(\mu'c'\big)\,
\zero{\kappa_0'\big(\mathcal{K}'\big)}
-
f'\overline{\nu}'\overline{b}'\overline{\mathcal{K}}'
\big(\mu'c'\big)\,
\zero{\kappa_0'\big(\mathcal{K}'\big)}
\\
&
\ \ \ \ \
+
e'\mu'\overline{\nu}'c'\overline{a}'\,
\rho_0'\big(\big[\mathcal{K}',\overline{\mathcal{L}}_1'\big]\big)
+
e'\mu'\overline{\nu}'c'\overline{b}'\,
\rho_0'\big(\big[\mathcal{K}',\overline{\mathcal{K}}'\big]\big)
+
\\
&
\ \ \ \ \ 
+
e'\mu'c'\mathcal{K}'\big(\overline{\nu}'\overline{a}'\big)\,
\zero{\rho_0'\big(\overline{\mathcal{L}}_1'\big)}
+
e'\mu'c'\mathcal{K}'\big(\overline{\nu}'\overline{b}'\big)\,
\zero{\rho_0'\big(\overline{\mathcal{K}}'\big)}
-
\\
&
\ \ \ \ \
-
e'\overline{\nu}'\overline{a}'\overline{\mathcal{L}}_1'
\big(\mu'c'\big)\,
\zero{\rho_0'\big(\mathcal{K}'\big)}
-
e'\overline{\nu}'\overline{b}'\overline{\mathcal{K}}'
\big(\mu'c'\big)\,
\zero{\rho_0'\big(\mathcal{K}'\big)},
\endaligned
\]
so that after eight direct zero-ifications:
\[
\aligned
\kappa_0
\Big(
\big[\mu\,\mathcal{K},\,\,
\overline{\nu}\,\overline{\mathcal{L}}_1\big]
\Big)
&
=
f'c'\overline{a}'\,
\mu'\overline{\nu}'\,
\kappa_0'
\big(\big[\mathcal{K}',\overline{\mathcal{L}}_1'\big]\big)
+
f'\mu'\overline{\nu}'c'\overline{b}'\,
\kappa_0'\big(\big[\mathcal{K}',\overline{\mathcal{K}}'\big]\big)
+
\\
&
\ \ \ \ \
+
e'\mu'\overline{\nu}'c'\overline{a}'\,
\rho_0'\big(\big[\mathcal{K}',\overline{\mathcal{L}}_1'\big]\big)
+
e'\mu'\overline{\nu}'c'\overline{b}'\,
\rho_0'\big(\big[\mathcal{K}',\overline{\mathcal{K}}'\big]\big).
\endaligned
\]

But at this point, reminding that:
\[
\big[K^{1,0}M',K^{1,0}M'\big]
\,\subset\,
K^{1,0}M'
\oplus
K^{0,1}M',
\]
one has with two certain functions $g', h'$:
\[
\big[\mathcal{K}',\,\overline{\mathcal{K}}'\big]
=
g'\,\mathcal{K}'
+
h'\,\overline{\mathcal{K}}',
\]
whence the second term in the right-hand side above vanishes:
\[
f'\mu'\overline{\nu}'c'\overline{b}'\,\kappa_0'
\big(\big[\mathcal{K}',\,\overline{\mathcal{K}}'\big]\big)
=
g'f'\mu'\overline{\nu}'c'\overline{b}'\,
\zero{\kappa_0'\big(\mathcal{K}'\big)}
+
h'f'\mu'\overline{\nu}'c'\overline{b}'\,
\zero{\kappa_0'\big(\overline{\mathcal{K}}'\big)},
\]
and similarly, the fourth term does also vanish:
\[
e'\mu'\overline{\nu}'c'\overline{b}'\rho_0'
\big(\big[\mathcal{K}',\,\overline{\mathcal{K}}'\big]\big)
=
g'e'\mu'\overline{\nu}'c'\overline{b}'\,
\zero{\rho_0'\big(\mathcal{K}'\big)}
+
h'e'\mu'\overline{\nu}'c'\overline{b}'\,
\zero{\rho_0'\big(\overline{\mathcal{K}}'\big)}.
\]
Lastly, the third term vanishes too, because
$\mathcal{ K}'$ being in the Levi kernel, one
has:
\[
\big[\mathcal{K}',\,\overline{\mathcal{L}}_1'\big]
=
p'\,\mathcal{L}_1'
+
q'\,\mathcal{K}'
+
r'\,\overline{\mathcal{L}}_1'
+
s'\,\overline{\mathcal{K}}',
\]
for certain four other functions.

\medskip

Thus, just the first term remains:
\[
\aligned
\kappa_0
\Big(
\big[\mu\,\mathcal{K},\,\,
\overline{\nu}\,\overline{\mathcal{L}}_1\big]
\Big)
&
=
f'c'\overline{a}'\,
\mu'\overline{\nu}'\,
\kappa_0'
\big(\big[\mathcal{K}',\overline{\mathcal{L}}_1'\big]\big),
\endaligned
\]
and since one easily verifies, by applying an already
seen reasoning, that:
\[
\mu'\overline{\nu}'\,
\kappa_0'
\big(\big[\mathcal{K}',\overline{\mathcal{L}}_1'\big]\big)
=
\kappa_0'
\big(\big[\mu'\,\mathcal{K}',\,\,\overline{\nu}'\,
\overline{\mathcal{L}}_1'\big]\big)
\]
one concludes that:
\[
\boxed{\,\,
\kappa_0
\Big(
\big[\mu\,\mathcal{K},\,\,
\overline{\nu}\,\overline{\mathcal{L}}_1\big]
\Big)
=
\underbrace{f'c'\overline{a}'}_{{\sf nonzero}
\atop
{\sf factor}}\,
\kappa_0'
\big(\big[\mu'\,\mathcal{K}',\,\,\overline{\nu}'\,
\overline{\mathcal{L}}_1'\big]\big).\,\,
}
\]
so that this quantity is a biholomorphic invariant\,\,---\,\,in the
sense of \'Elie Cartan\,\,---\,\,for hypersurfaces $M^5 \subset \C^3$
whose Levi form is everywhere of rank $1$.
\endproof

Having reached this corner-point, it is advisable to state 
precisely and in a synthetically summarized manner a:

\medskip\noindent{\bf Proposition.}
{\em In any system of holomorphic coordinates, 
for any choice of Levi-kernel adapted local $T^{1,0}M$-frame:}
\[
\big\{\mathcal{L}_1,\mathcal{K}\big\}
\]
{\em satisfying:}
\[
K^{1,0}M
=
\C\,\mathcal{K},
\]
{\em and for any choice of differential $1$-forms:}
\[
\big\{
\rho_0,\kappa_0,\zeta_0
\big\}
\]
{\em satisfying:}
\[
\aligned
\big\{0=\rho_0\big\}
&
=
T^{1,0}M\oplus T^{0,1}M,
\\
\big\{0=\rho_0=\kappa_0=\zeta_0=\overline{\zeta}_0\big\}
&
=
K^{1,0}M,
\endaligned
\]
{\em the quantity:}
\[
\kappa_0
\big(
\big[\mathcal{K},\,\overline{\mathcal{L}}_1\big]
\big),
\]
{\em is, at one fixed point $p \in M$, either $0$ 
or nonzero, independently of any choice.}

\proof
In what precedes, specific choices have been made for $\mathcal{ K}$,
$\mathcal{ L}_1$, $\rho_0$, $\kappa_0$, $\zeta_0$, but the transfer
formula:
\[
\kappa_0\Big(
\big[\mu\,\mathcal{K},\,\,
\overline{\nu}\,\overline{\mathcal{L}}_1\big]
\Big)
=
\mathmotsf{nonzero function}\cdot
\kappa_0'\Big(
\big[\mu'\,\mathcal{K}',\,\,
\overline{\nu}'\,\overline{\mathcal{L}}_1'\big]
\Big),
\]
and the reasonings made there included the fact that when one does
other choices, any change of choice has the same general form as when
dealing with a transfer through a local biholomorphism.  Hence the
zeroness or the nonzeroness of the interesting quantity is definitely
invariant.
\endproof

\medskip

One should notice the strong similarity of these reasonings 
with the introduction of the concept of Levi form.
Hence it is advisable to conceptualize in an 
analogous way what was obtained.

\medskip\noindent{\bf Definition.}
At any point:
\[
p
\in
M^5
\subset
\C^3
\]
of a $\mathcal{ C}^\kappa$ ($\kappa \geqslant 3$), or
$\mathcal{ C}^\infty$, or $\mathcal{ C}^\omega$
hypersurface whose Levi form is everywhere of $\C$-rank $1$,
so that the Levi-kernel subbundle:
\[
K^{1,0}M
\subset
T^{1,0}M
\]
is also of $\C$-rank $1 = 2 - 1$, the {\sl Freeman form}
is the $\C$-skew bilinear for on:
\[
\underbrace{K_p^{1,0}M}_{\cong\,\C}
\times
\big(
\underbrace{T_p^{1,0}M\,\,
\mod\,K_p^{1,0}M}_{\cong\,\C}
\big)
\,\longrightarrow\,
\C
\]
defined as follows: given any two constant vectors:
\[
\mathcal{K}_p
\,\in\,
K_p^{1,0}M
\ \ \ \ \ \ \ \ \ \ \ \ \ 
\text{\rm and}
\ \ \ \ \ \ \ \ \ \ \ \ \ 
\mathcal{L}_{1p}^\sim
\,\in\,
T_p^{1,0}M\,\,
\mod\,K_p^{1,0}M,
\]
take any two local vector field extensions:
\[
\mathcal{K}
\ \ \ \ \ \ \ \ \ \ \ \ \ 
\text{\rm and}
\ \ \ \ \ \ \ \ \ \ \ \ \ 
\mathcal{L}_1
\]
of $K^{1, 0}M$ and of $T^{1, 0}M$ satisfying hence:
\[
\mathcal{K}\big\vert_p
=
\mathcal{K}_p
\ \ \ \ \ \ \ \ \ \ \ \ \ 
\text{\rm and}
\ \ \ \ \ \ \ \ \ \ \ \ \ 
\mathcal{L}_1\big\vert_p
=
\mathcal{L}_{1p}^\sim,
\]
and define:
\[
\mathmotsf{Freeman-form}^{M,p}
\big(
\mathcal{K}_p,\,\mathcal{L}_{1p}
\big)
\overset{\sf def}{\,:=\,}
\big[\mathcal{K},\,\overline{\mathcal{L}}_1\big](p)
\ \ \ \ \
\mod\,\,
\big(
K^{1,0}M
\oplus
T^{0,1}M
\big).
\]

\medskip
One can show directly that the result is independent of the
choice of vector field extensions $\mathcal{ K}$ and
$\mathcal{ L}_1$ (exercise), but this property
will also be clarified with frames in coordinates
just below. 

\medskip

To be more precise about the `mod out' of the right-hand side,
is is better to introduce as previously a $(1, 0)$-form $\kappa_0$ 
satisfying:
\[
\big\{0=\kappa_0\big\}
=
K^{1,0}M
\ \ \ \ \
\mathmotsf{inside}\,\,
T^{1,0}M,
\]
and to define in accordance to what precedes:
\[
\mathmotsf{Freeman-form}^{M,p}
\big(
\mathcal{K}_p,\,\mathcal{L}_{1p}
\big)
\overset{\sf def}{\,:=\,}
\kappa_0
\Big(
\big[\mathcal{K},\,\overline{\mathcal{L}}_1\big]
\Big)(p).
\]

The computations motivating this definition yielded
that, within a Levi-kernel adapted $T^{1, 0}M$-frame:
\[
\big\{
\mathcal{L}_1,\mathcal{K}
\big\},
\]
one may decompose the two constant vectors:
\[
\aligned
\mathcal{K}_p
&
=
\mu_p\,\mathcal{K}\big\vert_p,
\\
\mathcal{L}_{1p}
&
=
\nu_{1p}\,\mathcal{L}_1\big\vert_p,
\endaligned
\]
with two constants $\mu_p, \nu_{ 1p} \in \C$,
and extend them both as:
\[
\aligned
&
\mu\,
\mathcal{K}
\\
&
\nu_1\,
\mathcal{L}_1,
\endaligned
\]
by means of two local functions $\mu$ and $\nu_1$:
\[
\aligned
\mu\colon\ \ \ 
&
M\,\longrightarrow\,\C,
\\
\nu_1\colon\ \ \ 
&
M\,\longrightarrow\,\C,
\endaligned
\]
satisfying:
\[
\aligned
\mu(p)
&
=
\mu_p,
\\
\nu_1(p)
&
=
\nu_{1p},
\endaligned
\]
and one realizes that: 
\[
\mathmotsf{Freeman-form}_{\mathcal{K},\overline{\mathcal{L}},
\kappa_0}^{M,p}
\big(
\mathcal{K}_p,\,\mathcal{L}_{1p}
\big)
\,=\,
\mu_p\,\overline{\nu}_{1p}\,\,
\underbrace{\kappa_0
\big(
\big[\mathcal{K},\,\overline{\mathcal{L}}_1\big]
\big)(p)}_{\constant},
\]
an expression which indeed shows (again) that the result
is independent of vector field extensions.

Notice that the Freeman form is $\C$-skew bilinear, but not
necessarily Hermitian, for the constant above
needs not belong either to $\R$ or to $\isqrt\, \R$.

\medskip\noindent{\bf Analysis of everywhere degeneracy of the Freeman
form.}
Applying then the:
\[
\text{\sl Lie-Cartan Principle of Relocalization},
\]
one is lead to a natural dichotomy.
Either:
\[
\mathmotsf{Freeman-Form}^M(p)
\equiv
0,
\]
or:
\[
\mathmotsf{Freeman-Form}^M(p)
\neq
0,
\]
at every point $p \in M$.

\medskip

Examine the first possibility. This is a degenerate situation.

\medskip\noindent{\bf Proposition.}
{\em A $\mathcal{ C}^\omega$ hypersurface:}
\[
M^5
\subset
\C^3
\]
{\em having at every point $p$:}
\[
\rank_\C\big(\mathmotsf{Levi-Form}^M(p)\big)
=
{\bf 1}
\]
{\em has an identically vanishing:}
\[
\mathmotsf{Freeman-Form}^M(p)
\equiv
{\bf 0},
\]
{\em if and only if it is biholomorphic, locally
in some neighborhood of every point, to a product:}
\[
M^5
\,\cong\,
M^3
\times
\C
\]
{\em with a $\mathcal{ C}^\omega$ hypersurface:}
\[
M^3
\subset
\C^2.
\]

\medskip

In local coordinates, the graphing function:
\[
v
=
\varphi\big(x_1,y_1,u\big)
\]
happens then to be completely independent of $x_2, y_2$.

\medskip\noindent{\bf Interpretation.} 
{\em One sets aside such an exceptional supposition, because the
equivalence problem reduces to that of an:}
\[
M^3
\subset
\C^2
\]
{\em in smaller dimension, plus $1$ complex parameter coming
from $(\cdot) \times \C$.\qed}

\proof[Proof of the Proposition]
One centers affine holomorphic coordinates:
\[
(z_1,z_2,w)
\]
at some point $p \in M$ and one represents $M$ as usual:
\[
v
=
\varphi\big(x_1,x_2,y_1,y_2,u\big).
\]
One introduces the two generators of $T^{1, 0}M$:
\[
\aligned
\mathcal{L}_1
&
=
\frac{\partial}{\partial z_1}
\underbrace{-
\frac{\varphi_{z_1}}{\isqrt+\varphi_u}}_{A_1}\,
\frac{\partial}{\partial u},
\\
\mathcal{L}_2
&
=
\frac{\partial}{\partial z_2}
\underbrace{-
\frac{\varphi_{z_2}}{\isqrt+\varphi_u}}_{A_2}\,
\frac{\partial}{\partial u},
\endaligned
\]
and, with the hypothesis that the Levi determinant
vanishes, the generator:
\[
\aligned
\mathcal{K}
&
=
k\,\mathcal{L}_1
+
\mathcal{L}_2
\\
&
=
k\,\frac{\partial}{\partial z_1}
+
\frac{\partial}{\partial z_2}
+
\big(k\,A_1+A_2\big)\,
\frac{\partial}{\partial u},
\endaligned
\]
of the Levi-kernel bundle $K^{1, 0}M$, where the
slanting function $k$ happens, according to
an already seen surprising computational fact, to receive
{\em three equal expressions:}
\[
\aligned
k
&
=
-\,
\frac{
\mathcal{L}_2\big(\overline{A_1}\big)
-
\overline{\mathcal{L}}_1\big(A_2\big)
}{
\mathcal{L}_1\big(\overline{A_1}\big)
-
\overline{\mathcal{L}}_1\big(A_1\big)
}
\\
&
=
-\,
\frac{\mathcal{L}_2\big(\overline{A_1}\big)}{
\mathcal{L}_1\big(\overline{A_1}\big)}
\\
&
=
-\,
\frac{-\,\overline{\mathcal{L}}_1\big(A_2\big)}{
-\,\overline{\mathcal{L}}_1\big(A_1\big)}.
\endaligned
\]

First of all, the known involutiveness:
\[
\big[\mathcal{K},\,\overline{\mathcal{K}}\big]
=
\function\cdot
\mathcal{K}
+
\function\cdot
\overline{\mathcal{K}},
\]
and the fact that this bracket does not contain
either $\frac{ \partial}{ \partial z_2}$ or
$\frac{ \partial}{ \partial \overline{ z}_2}$ entail that:
\[
\aligned
0
&
=
\big[
\mathcal{K},\,
\overline{\mathcal{K}}
\big]
\\
&
=
\bigg[
k\,\frac{\partial}{\partial z_1}
+
\frac{\partial}{\partial z_2}
+
\big(k\,A_1+A_2\big)\,
\frac{\partial}{\partial u},\,\,\,
\overline{k}\,\frac{\partial}{\partial\overline{z}_1}
+
\frac{\partial}{\partial\overline{z}_2}
+
\big(\overline{k}\,\overline{A}_1+\overline{A}_2\big)\,
\frac{\partial}{\partial u}
\bigg]
\\
&
=
\mathcal{K}\big(\overline{k}\big)\,
\frac{\partial}{\partial\overline{z}_1}
-
\overline{\mathcal{K}}\big(k\big)\,
\frac{\partial}{\partial z_1}
+
\Big(
\mathcal{K}\big(\overline{k}\,\overline{A}_1+\overline{A}_2\big)
-
\overline{\mathcal{K}}\big(k\,A_1+A_2\big)
\Big)\,
\frac{\partial}{\partial u},
\endaligned
\]
whence one deduces at least that:
\[
\boxed{\,
0
\equiv
\overline{\mathcal{K}}(k).\,}
\]

Next, the (assumed) zeroness of the Freeman form means that:
\[
\big[\mathcal{K},\,\overline{\mathcal{L}}_1\big]
\,\equiv\,
0
\ \ \ \ \ 
\mod\,\,
\big(
\mathcal{K},\,\overline{\mathcal{K}},\,\overline{\mathcal{L}}_1
\big).
\]
But when one indeed computes this bracket:
\[
\aligned
\big[\mathcal{K},\,\overline{\mathcal{L}}_1\big]
&
=
\bigg[
k\,\frac{\partial}{\partial z_1}
+
\frac{\partial}{\partial z_2}
+
\big(k\,A_1+A_2\big)\,
\frac{\partial}{\partial u},\,\,\,
\frac{\partial}{\partial\overline{z}_1}
+
\overline{A}_1\,\frac{\partial}{\partial u}
\bigg]
\\
&
=
\overline{\mathcal{L}}_1(k)\,
\frac{\partial}{\partial z_1}
+
\something\,
\frac{\partial}{\partial u},
\endaligned
\]
one obtains a $\frac{ \partial}{\partial z_1}$-component
which is {\em nonzero} modulo $\big\{ 
\mathcal{ K}, \, \overline{ \mathcal{ K}},\,
\overline{ \mathcal{L}}_1 \big\}$, whence one obtains also:
\[
\boxed{\,
0
\,\equiv\,
\overline{\mathcal{L}}_1\big(k\big).
\,}
\]

Next, because:
\[
\big\{
\overline{\mathcal{K}},\,
\overline{\mathcal{L}}_1\big\}
\]
constitute a frame for $T^{ 0, 1} M$, the two latter boxed equations
yield that:
\[
\boxed{\,\text{\sl The $\mathcal{ C}^\omega$ slanting function $k$ is
Cauchy-Riemann!}\,}
\]

Hence, as is known and was reproved earlier on, there exists
a {\em holomorphic} function $K$ locally defined in some
open neighborhood of $M$ which extends $k$:
\[
K
\big\vert_M
=
k.
\]

The first coefficient of the $(1, 0)$-vector field:
\[
\mathcal{K}
=
K\,\frac{\partial}{\partial z_1}
+
\frac{\partial}{\partial z_2}
+
\big(K\,A_1+A_2\big)\,
\frac{\partial}{\partial u}
\]
is hence {\em holomorphic}.

What about the second and last coefficient:
\[
k\,A_1+A_2
=
K\,A_1+A_2\big\vert_M\,?
\]
Is it also CR? 

Yes, mainly thanks to the surprising computational fact
recalled above, firstly:
\[
\aligned
0
&
\overset{?}{\,=\,}
\overline{\mathcal{L}}_1
\big(
k\,A_1+A_2\big)
\\
&
\,=\,
k\,\overline{\mathcal{L}}_1(A_1)
+
\overline{\mathcal{L}}_1(A_2)
\\
&
\,=\,
0,
\endaligned
\]
while secondly a direct painful computation gives:
\[
\aligned
\!\!\!\!\!\!\!\!\!\!\!\!\!\!\!\!\!\!\!\!\!\!\!\!\!\!\!\!\!\!
0
&
\overset{?}{\,=\,}
\overline{\mathcal{L}}_2
\big(
k\,A_1+A_2\big)
\\
\!\!\!\!\!\!\!\!\!\!\!\!\!\!\!\!\!\!\!\!\!\!\!\!\!\!\!\!\!\!
&
\,=\,
k\,\overline{\mathcal{L}}_2(A_1)
+
\overline{\mathcal{L}}_2(A_2)
\\
\!\!\!\!\!\!\!\!\!\!\!\!\!\!\!\!\!\!\!\!\!\!\!\!\!\!\!\!\!\!
&
\,=\,
\frac{-\,\mathmotsf{numerator-LD}}{
(\isqrt+\varphi_u)\,
[\varphi_{z_1\overline{z}_1}+\varphi_{z_1\overline{z}_1}\varphi_u\varphi_u
-
\isqrt\,\varphi_{\overline{z}_1}\varphi_{z_1u}
-
\varphi_{\overline{z}_1}\varphi_{z_1u}\varphi_u
+
\isqrt\,\varphi_{z_1}\varphi_{\overline{z}_1u}\varphi_u
+
\varphi_{z_1}\varphi_{\overline{z}_1}\varphi_{uu}]}
\endaligned
\]
where:
\[
\aligned
\!\!\!\!\!\!\!\!\!\!\!\!\!\!\!\!\!\!\!\!\!\!\!\!\!\!\!\!\!\!
\mathmotsf{numerator-LD}
&
\equiv
\varphi_{z_2\overline{z}_2}\varphi_{z_1\overline{z}_1}
-
\varphi_{z_2\overline{z}_1}\varphi_{z_1\overline{z}_2}
+
\\
&\ \ \ \ \
+
\varphi_{z_2\overline{z}_1}\varphi_{\overline{z}_2}\varphi_{z_1u}
\varphi_u
-
\varphi_{z_2\overline{z}_1}\varphi_{\overline{z}_2}\varphi_{z_1}
\varphi_{uu}
-
\varphi_{\overline{z}_1}\varphi_{z_2u}\varphi_{z_1}
\varphi_{\overline{z}_2u}
+
\varphi_{\overline{z}_1}\varphi_{z_2u}\varphi_u
\varphi_{z_1\overline{z}_2}
-
\\
&\ \ \ \ \
-\,
\varphi_{z_2}\varphi_{\overline{z}_1u}\varphi_{\overline{z}_2}
\varphi_{z_1u}
-
\varphi_{z_2}\varphi_{\overline{z}_1}\varphi_{uu}
\varphi_{z_1\overline{z}_2}
+
\varphi_{z_2}\varphi_{\overline{z}_1u}
\varphi_u\varphi_{z_1\overline{z}_2}
-
\varphi_{z_2\overline{z}_2}\varphi_{\overline{z}_1}
\varphi_{z_1u}\varphi_u
+
\\
&\ \ \ \ \
+
\varphi_{z_2\overline{z}_2}\varphi_{z_1}
\varphi_{\overline{z}_1}\varphi_{uu}
-
\varphi_{z_2\overline{z}_2}\varphi_{z_1}
\varphi_{\overline{z}_1u}\varphi_u
+
\varphi_{z_2\overline{z}_1}\varphi_{z_1}
\varphi_{\overline{z}_2u}\varphi_u
+
\varphi_{z_2}\varphi_{\overline{z}_2u}
\varphi_{\overline{z}_1}\varphi_{z_1u}
-
\\
&\ \ \ \ \
-\,
\varphi_{z_2}\varphi_{\overline{z}_2u}
\varphi_{z_1\overline{z}_1}\varphi_u
+
\varphi_{\overline{z}_2}\varphi_{z_2u}
\varphi_{z_1}\varphi_{\overline{z}_1u}
-
\varphi_{\overline{z}_2}\varphi_{z_2u}
\varphi_u\varphi_{z_1\overline{z}_1}
+
\varphi_{\overline{z}_2}\varphi_{z_2}
\varphi_{uu}\varphi_{z_1\overline{z}_1}
+
\\
&\ \ \ \ \
+
\isqrt
\big(
\varphi_{z_2\overline{z}_2}\varphi_{z_1}
\varphi_{\overline{z}_1u}
+
\varphi_{\overline{z}_1}\varphi_{z_2u}
\varphi_{z_1\overline{z}_2}
+
\varphi_{z_2\overline{z}_1}\varphi_{\overline{z}_2}
\varphi_{z_1u}
+
\varphi_{z_2}\varphi_{\overline{z}_2u}
\varphi_{z_1\overline{z}_1}
\big)
-
\\
&\ \ \ \ \
-\,\isqrt\big(
\varphi_{\overline{z}_2}\varphi_{z_2u}
\varphi_{z_1\overline{z}_1}
+
\varphi_{z_2\overline{z}_1}\varphi_{z_1}
\varphi_{\overline{z}_2u}
+
\varphi_{z_2}\varphi_{\overline{z}_1u}
\varphi_{z_1\overline{z}_2}
+
\varphi_{z_2\overline{z}_2}\varphi_{\overline{z}_1}
\varphi_{z_1u}
\big)
-
\\
&\ \ \ \ \
-\,
\varphi_{z_2\overline{z}_1}\varphi_{z_1\overline{z}_2}
\varphi_u\varphi_u
+
\varphi_{z_2\overline{z}_2}\varphi_{z_1\overline{z}_1}
\varphi_u\varphi_u
\\
&
\equiv
0
\endaligned
\]
was, up to a constant factor $4$, the numerator
of the Levi determinant already shown above,
and assumed throughout to be identically zero!

\medskip

Once again, this function:
\[
k\,A_1+A_2
\]
being CR, it also happens to be the restriction, to $M$,
of a certain holomorphic function.

{\em Consequently, $\mathcal{ K}$ is not
only a $(1, 0)$-vector field, it is a perfect
holomorphic vector field with holomorphic coefficients}.

To conclude, one uses a local biholomorphism:
\[
(z_1,z_2,w)
\,\longmapsto\,
(z_1',z_2',w')
\]
which straightens out:
\[
\mathcal{K}'
=
\frac{\partial}{\partial z_2'}.
\]
Since $\mathcal{ K}'$ is again tangent to the
image $M'$, this means (exercise), dropping then primes
on coordinates, that $M$ becomes a product:
\[
M^3\times
\C_{z_2},
\]
with graphing function being independent of $(x_2, y_2)$.
\endproof

After all these detailed considerations, one therefore arrives
at the very last, sixth {\sl general class}
of CR-generic manifolds:
\[
\aligned
&
\boxed{\text{\sf General Class $\text{\sf IV}_{\text{\sf 2}}$:}}
\\
&
\boxed{\,\,
\aligned
M^5\subset\C^3
\ \
&
\text{\rm with}\ \
\Big\{\mathcal{L}_1,\,\mathcal{L}_2,\,
\overline{\mathcal{L}}_1,\,\overline{\mathcal{L}}_2,\,\,
\big[\mathcal{L}_1,\overline{\mathcal{L}}_1\big]
\Big\}\,\,
\\
&
\text{\rm constituting a frame for}\ \
\C\otimes_\R TM,
\\
&
\text{\rm with the Levi-Form:}\ \
\\
&
\ \ \ \ \ \ \ \ \ \ \ \ \ \ \
\mathmotsf{Levi-Form}^M(p)
\\
&
\text{\rm being of rank}\,\,{\bf 1}\,\,
\text{\rm at every point}\,\,p\in M,
\\
&
\text{\rm while the Freeman-Form:}\ \
\\
&
\ \ \ \ \ \ \ \ \ \ \ \ \ \ \
\mathmotsf{Freeman-Form}^M(p)
\\
&
\text{\rm is nondegenerate at every point.}\ \
\endaligned}
\endaligned
\]

%%%%%%%%%%%%%%%%%%%%%%%%%%%%%%%%%%%%%%%%%%%%%%%%%%%%%%%%%%%%%%%%%%%%%

\bigskip

\section{\sf General classes 
$\text{\sf I}$, 
$\text{\sf II}$, 
$\text{\sf III}_{\text{\sf 1}}$,
$\text{\sf III}_{\text{\sf 2}}$,
$\text{\sf IV}_{\text{\sf 1}}$,
$\text{\sf IV}_{\text{\sf 2}}$,
\\
of $M^3 \subset \C^2$, of $M^4 \subset \C^3$, of $M^5 \subset \C^4$,
of $M^5 \subset \C^3$}
\label{summary-general-classes}
\HEAD{\ref{summary-general-classes}.~General classes 
$\text{\sf I}$, 
$\text{\sf II}$, 
$\text{\sf III}_{\text{\sf 1}}$,
$\text{\sf III}_{\text{\sf 2}}$,
$\text{\sf IV}_{\text{\sf 1}}$,
$\text{\sf IV}_{\text{\sf 2}}$}{
Jo\"el {\sc Merker} (Paris-Sud), 
Samuel {\sc Pocchiola} (Paris-Sud), 
Masoud {\sc Sabzevari} (Shahrekord)}

\medskip

In conclusion, there are precisely {\em six} general
classes of real analytic CR-generic manifolds up to dimension 
{\bf 5}:
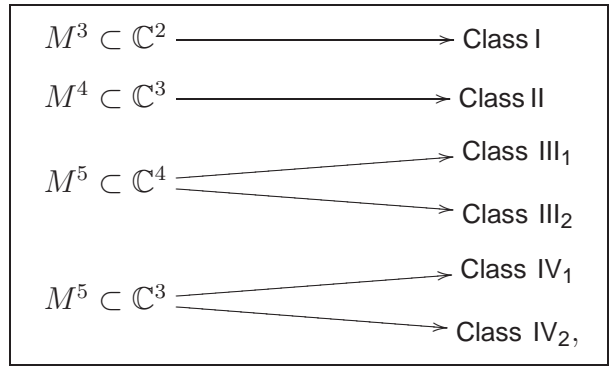
\[
\boxed{\,\,
\aligned
&
\xymatrix @!0 @R=1em @C=12.5pc {
M^3
\subset
\C^2
\ar[r]
&
\text{\footnotesize\sf Class}\,
\text{\footnotesize\sf I}
}
\\
&
\xymatrix @!0 @R=1em @C=12.5pc {
M^4
\subset
\C^3
\ar[r]
&
\text{\footnotesize\sf Class}\,
\text{\footnotesize\sf II}
}
\\
&
\xymatrix @!0 @R=1em @C=13pc {
& 
\text{\footnotesize\sf Class}\,\,
\text{\footnotesize\sf III}_{\text{\sf 1}}
\\ 
M^5\subset\C^4
\ar[ru] 
\ar[rd] 
\\ 
& 
\text{\footnotesize\sf Class}\,\,
\text{\footnotesize\sf III}_{\text{\sf 2}}
}
\\
&
\xymatrix @!0 @R=1em @C=13pc {
& 
\text{\footnotesize\sf Class}\,\,
\text{\footnotesize\sf IV}_{\text{\sf 1}}
\\ 
M^5\subset\C^3
\ar[ru] 
\ar[rd] 
\\ 
& 
\text{\footnotesize\sf Class}\,\,
\text{\footnotesize\sf IV}_{\text{\sf 2}},
}\,\,
\endaligned}
\]
when one disregards the degenerate classes.

%%%%%%%%%%%%%%%%%%%%%%%%%%%%%%%%%%%%%%%%%%%%%%%%%%%%%%%%%%%%%%%%%%%%%

\vfill\end{document}